\documentclass[12pt, reqno]{amsart}
\usepackage{amsmath,amssymb,amsthm}
\usepackage{amsmath, amssymb}
\usepackage{amsthm, amsfonts, mathrsfs}
\usepackage{mathptmx}
\usepackage{fullpage}
\usepackage{amsfonts,graphicx}
\usepackage{xcolor}
\colorlet{RED}{red}% amsart upper-cases some arguments

\numberwithin{equation}{section}
\usepackage[colorlinks=true, pdfstartview=FitV, linkcolor=blue, citecolor=blue, urlcolor=blue]{hyperref}

\numberwithin{equation}{section}
\newtheorem{theorem}{Theorem}[section]
\newtheorem{proposition}[theorem]{Proposition}
\newtheorem{lemma}[theorem]{Lemma}

\numberwithin{equation}{section}

\theoremstyle{remark}
\newtheorem{remark}[theorem]{Remark}

\newcommand{\beq}{\begin{equation}}
\newcommand{\eeq}{\end{equation}}

\begin{document}
\title[Global solutions  to the full compressible    non-resistive MHD system]{Global well-posedness for the 2D  full compressible viscous   non-resistive   magnetohydrodynamic
system}

\author[Qiao, Sun, Wu, Xu]{Liening Qiao$^{1}$, Juntao Sun$^{1}$,  Jiahong Wu$^{2}$ and Fuyi Xu$^{1\dag}$  }

\address{$^{1}$ School of Mathematics and  Statistics, Shandong University of
Technology,  Zibo 255049,  Shandong Province,    P.R. China}

\email{lnqiao2026@163.com}

\address{$^{1}$ School of Mathematics and  Statistics, Shandong University of
Technology,  Zibo 255049,  Shandong Province,    P.R. China}

\email{jtsun@sdut.edu.cn}

\address{$^2$ Department of Mathematics, University of Notre
	Dame, Notre Dame, IN 46556, USA}

\email{jwu29@nd.edu}

\address{$^{1\dag}$  School of Mathematics and  Statistics, Shandong University of
Technology,  Zibo 255049,  Shandong Province,    P.R. China}

\email{zbxufuyi@163.com (Corresponding author)}



\begin{abstract}{This paper establishes the global well-posedness and nonlinear stability of the 2D full (non-isentropic, heat-conducting) compressible viscous magnetohydrodynamic system without magnetic diffusion in the periodic domain $\mathbb T^2$, for smooth initial data close to the equilibrium with constant density and temperature, zero velocity and background magnetic field $\mathbf e_2=(0,1)$. Neither the density nor the magnetic field is damped or diffused in this system. Moreover, in contrast with the existing results on the non-isentropic non-resistive system, we impose no Diophantine condition on the background magnetic field, no symmetry on the data and no special geometry on the domain; in fact $\mathbf e_2$ violates every Diophantine condition, so that the stabilizing effect of the magnetic field degenerates completely in the $x_1$-direction. Our proof rests on two intrinsic stabilizing mechanisms hidden in the structure of the system. First, the background magnetic field generates damped wave structures which provide dissipation for the $x_2$-derivatives of the density and of the magnetic field. Second, we introduce a new good unknown $\Phi$, an effective total pressure gradient that combines the density, the temperature with the magnetic pressure. Although none of its constituents is dissipative by itself, we show that the coupling of $\Phi$ with the divergence of the velocity and with the temperature produces dissipation for $\Phi$ in all directions, together with an extra regularization of the divergence of the velocity. Combining these two mechanisms with new anisotropic commutator estimates and a hierarchy of time-weighted energy estimates, we prove that the system admits a unique global solution in $H^N(\mathbb T^2)$ for every $N\geq5$, and obtain a uniform-in-time bound on the highest-order energy as well as explicit algebraic decay rates. Even in comparison with the known result for the simpler isentropic counterpart, where the highest-order energy was only shown to grow at most algebraically in time for data in $H^{N}$ with $N\geq10$, our result removes the possible growth and lowers the regularity requirement to $N\geq5$.}
\end{abstract}
\maketitle
\tableofcontents

\section{Introduction and main results}
\subsection{Introduction }
{This paper is devoted to the global well-posedness and nonlinear stability, with a uniform-in-time bound on the highest-order Sobolev norm, of the two-dimensional (2D) full compressible viscous magnetohydrodynamic (MHD) system without magnetic diffusion near a background magnetic field. In the periodic domain $\mathbb T^2$, the system reads as follows:}
\begin{align} \label{heat-mhd1}
\left\{
\begin{aligned}
&\partial_{t}\rho+\mathrm{div}(\rho\mathbf u)=0,\\
&\partial_{t}(\rho\mathbf u)+\mathrm{div}(\rho\mathbf u\otimes\mathbf u)-\mu\Delta \mathbf u-(\mu+\lambda)\nabla\mathrm{div}\mathbf u+\nabla(P+\frac12|\mathbf
H|^2)=(\mathbf H\cdot\nabla)\mathbf H, \\
&c_\nu\big(\partial_{t}(\rho\Theta)+\mathrm{div}(\rho\Theta\mathbf u)\big)+P\mathrm{div}\mathbf u-\kappa\Delta\Theta=\mathcal Q(\mathbf
u), \\
&\partial_t\mathbf H+(\mathbf u\cdot \nabla)\mathbf H
=(\mathbf H\cdot \nabla)\mathbf u-\mathbf H\mathrm{div}\ \mathbf u,\\
&\mathrm{div}\mathbf H=0,
\end{aligned}
\right.
\end{align}
where  time $t\geq0$, the spatial coordinate $x=(x_1,x_2)\in \mathbb T^2$,     $\rho(x,t)\in\mathbb R^+$,  $\mathbf u(x,t)\in\mathbb R^2$, $\Theta(x,t)\in\mathbb R^+$, $\mathbf H(x,t)\in\mathbb R^2$ and $P(x,t)\in\mathbb R$ denote the
density, velocity field, absolute temperature, magnetic field and pressure, respectively.
 The positive
constants $c_\nu$ and $\kappa$  stand for the specific heat at
constant volume and the coefficient of heat conduction. The constants $\mu$ and $\lambda$ are the shear viscosity and
the second viscosity coefficient, and are subject to the physical
constraints
$$\mu>0,\ \ \mu+\lambda\geq0.$$
 The symbol $\otimes$ denotes the tensor product, $\mathcal Q(\mathbf u)$ is  given by
$$\mathcal Q(\mathbf u):=\frac{\mu}{2}\left|\nabla\mathbf u+(\nabla\mathbf u)^{T}\right|^2+\lambda(\mathrm{div}\ \mathbf u)^2,$$
where $(\nabla\mathbf u)^{T}$ is the transpose of $\nabla\mathbf u.$ Throughout the paper the fluid is a perfect gas, so that the pressure obeys
the equation of state
\begin{equation}\label{state-law}
	P=R\rho\Theta,\quad R>0 .
\end{equation}

{Magnetohydrodynamics describes the mutual interaction between electrically conducting fluids and electromagnetic fields, and its governing equations are obtained by coupling the compressible Navier--Stokes equations with Maxwell's equations. Owing to their physical relevance, rich dynamics and mathematical complexity, the multidimensional compressible MHD equations have attracted extensive attention. When viscosity, heat conduction and magnetic resistivity are all present, the analysis shares many features with that of the compressible Navier--Stokes equations, and a variety of results on the global existence and large-time behavior of solutions are available, both for the full (non-isentropic) system \cite{Hu-Wang1,Pu-Guo,Xu-Zhong} and for the isentropic system \cite{Gao-Wu-Xu,Hu-Wang2,KawashimaS,A.S,SuenA}; see also the references therein.

The system \eqref{heat-mhd1} is of a very different nature. Since the magnetic diffusion is absent, the induction equation $\eqref{heat-mhd1}_4$ is, like the continuity equation $\eqref{heat-mhd1}_1$, of transport type and contains neither dissipation nor damping. Thus two of the four unknowns, the density and the magnetic field, are not regularized or damped by any explicit mechanism, while the nonlinear terms couple them to the velocity and the temperature at the highest order of derivatives. For this reason the global well-posedness and stability problem for \eqref{heat-mhd1} is extremely challenging. On the other hand, physical experiments and numerical simulations on electrically conducting fluids such as liquid metals have revealed a remarkable phenomenon: a background magnetic field can smooth and stabilize turbulent conducting fluids (see \cite{AMS, Alex, Alf, Davi0, Davi1, DA}). A natural way to capture this stabilizing effect rigorously, which has guided much of the recent work, is to study the compressible MHD equations near an equilibrium carrying a nontrivial constant magnetic field, and to show that the coupling with this field generates enough dissipation for the unknowns that are not dissipated directly.

\vskip .05in
\noindent\textbf{Known results.} We first review the isentropic case. For the Cauchy problem, Hu \cite{HuX} established the global well-posedness of classical solutions to the 2D compressible non-resistive MHD equations with initial data close to the constant state $(1, \mathbf 0, \mathbf e_2)$ with $\mathbf e_2=(0,1)$. Wu and Wu \cite{Wu-Wu} subsequently developed a systematic approach to the global well-posedness and stability of the 2D system near the equilibrium $(1, \mathbf 0, \mathbf e_1)$ with $\mathbf e_1=(1,0)$. For strip domains, Tan and Wang \cite{Tan-Wang} established the global well-posedness and stability of smooth solutions with small initial perturbations around the background magnetic field $\mathbf e_3=(0,0,1)$. In periodic domains, Wu and Zhu \cite{Wu-Zhu} proved the global well-posedness and stability of smooth solutions to the 2D compressible non-resistive MHD equations with small initial perturbations around the background magnetic field $\mathbf e_2=(0,1)$. In the 3D periodic setting, Wu and Zhai \cite{Wu-Zhai} obtained a corresponding result when the background magnetic field satisfies a Diophantine condition. Recall that a vector \(\mathbf n\in\mathbb R^n\) is said to satisfy the Diophantine condition if there exist constants \(c>0\) and \(r>n-1\) such that
\begin{equation*}
|\mathbf n\cdot\mathbf k|
\geq \frac{c}{|\mathbf k|^r},
\qquad
\forall\,\mathbf k\in\mathbb Z^n\setminus\{0\}.
\end{equation*}
Jiu, Liu and Xie \cite{Jiu-Liu-Xie} further relaxed the regularity requirements on the initial data in the framework of \cite{Wu-Zhai}. It should be emphasized that the Diophantine condition excludes every background magnetic field whose components are all rational numbers. More recently, Li, Li, Wu and Xu \cite{Li-Li-Wu-Xu} established the global well-posedness and nonlinear stability of classical solutions to the 3D compressible viscous non-resistive MHD equations near the equilibrium $(1,\mathbf 0,\mathbf e_3)$ without a Diophantine assumption, for small initial perturbations satisfying a suitable $x_3$-symmetry condition. Further developments for the isentropic system can be found in \cite{Jiang-Jang2, Li-Qiao, Qiao-Zhai, Zhai-Wu-Xu,ZhaoY2} and the references therein.

For the full compressible system, the global analysis encounters additional difficulties stemming from the coupling among the temperature, the density, the velocity and the magnetic field, and far fewer results are available. For strip domains, Li \cite{LiY} extended the result of Tan and Wang \cite{Tan-Wang} to the non-isentropic setting. For periodic domains, Li, Xu and Zhai \cite{Li-Xu-Zhai} proved the global well-posedness and stability of smooth solutions to the 3D full compressible non-resistive MHD equations with small initial perturbations around a background magnetic field satisfying a Diophantine condition. When heat conduction and resistivity are present but viscosity is absent, the underlying inviscid flow may form shocks and results are scarce: Wang and Xin \cite{Wang-Xin} proved the global well-posedness for small data near a non-horizontal constant background field in a 3D strip domain, and Wu, Xu and Zhai \cite{Wu-Xu-Zhai} treated the periodic case under a Diophantine condition. Recently, Qiao, Sun, Wu and Xu \cite{Qiao-Xu} obtained global smooth solutions with explicit algebraic decay rates for the 3D full compressible viscous non-resistive MHD system without heat conductivity, provided that the initial velocity, pressure and magnetic field are small perturbations of the equilibrium and the background magnetic field satisfies a Diophantine condition; in particular, the initial density is allowed to have arbitrarily large variations.

\vskip .05in
\noindent\textbf{Motivation.} Two features of the existing theory should be emphasized. First, all the aforementioned results on the non-isentropic non-resistive system are confined to the three-dimensional case and rely on additional structure: either the special geometry of a strip domain, or the requirement that the background magnetic field satisfy a Diophantine condition. This additional structure is not a technical convenience; it is precisely what yields the coercive estimates that compensate for the lack of magnetic diffusion. Second, even for the simpler isentropic system in $\mathbb T^2$, the global solutions constructed in \cite{Wu-Zhu} require the initial data to be in $H^N$ with $N\geq10$, and their highest-order energy is not known to remain bounded: it is only controlled by a growing factor $(1+t)^\sigma$ with $\sigma>0$.

The aim of this paper is to solve the global well-posedness and stability problem for the 2D full compressible system \eqref{heat-mhd1} in $\mathbb T^2$ near the constant equilibrium $(\rho,\mathbf u,\Theta,\mathbf H)=(1,\mathbf 0,1,\mathbf e_2)$ \emph{without any additional assumption}. None of the results quoted above applies to this setting. Indeed, the background magnetic field $\mathbf e_2=(0,1)$ violates every Diophantine condition, since $\mathbf e_2\cdot\mathbf k=0$ for the infinitely many frequencies $\mathbf k=(k_1,0)$. As a consequence, the estimates that rely on the Diophantine inequality are not available, and the stabilizing effect of the magnetic field is genuinely degenerate: it acts only on the $x_2$-derivatives of the density and of the magnetic field, and it gives no control at all on the modes that do not depend on $x_2$. At the same time, compared with the isentropic case, the temperature equation creates new nonlinear interactions at the highest order of derivatives.

\vskip .05in
\noindent\textbf{Main contributions.} The main result of this paper, Theorem \ref{th:main1} below, states that \eqref{heat-mhd1} admits a unique global solution for any initial data that is a small perturbation of $(1,\mathbf 0,1,\mathbf e_2)$ in $H^N(\mathbb T^2)$ with $N\geq5$; the solution remains uniformly bounded in $H^N$ for all time and decays at explicit algebraic rates. The novelty of this work can be summarized as follows.
\begin{itemize}
\item[(1)] \emph{No Diophantine condition, no symmetry and no special geometry.} Theorem \ref{th:main1} holds for the background magnetic field $\mathbf e_2$, all of whose components are rational, in a fully periodic two-dimensional domain; the only structural condition on the data, besides smallness, is the natural zero-mean condition \eqref{junzhi} on the first component of the magnetic perturbation, which is preserved by the flow.
\item[(2)] \emph{A new good unknown and a hidden dissipative structure.} The key new ingredient is the quantity
$$\Phi:=R(1+\theta)\nabla a+\nabla \mathbf h_2+\nabla\frac{|\mathbf h|^2}{2},$$
where $(a,\theta,\mathbf h)$ denotes the perturbation of the density, the temperature and the magnetic field. It is an effective total pressure gradient built from the density, the temperature and the magnetic pressure. None of its constituents is dissipative by itself. Nevertheless, we discover that $\Phi$, $\mathrm{div}\mathbf u$ and $\theta$ form a partially dissipative subsystem which is stable at the linear level (see \eqref{linear-phi-divu-theta}), and we exploit this structure at the nonlinear level to obtain dissipation for $\Phi$ in \emph{all} directions, together with an extra regularization of $\mathrm{div}\mathbf u$. This mechanism, which is distinct from the wave-type stabilization of $\partial_2 a$ and $\partial_2\mathbf h$ induced by the background field, supplies exactly the control in the $x_1$-direction that the magnetic stabilization cannot provide.
\item[(3)] \emph{A uniform-in-time bound on the highest-order energy.} The derivatives of the density and of the magnetic field in the $x_1$-direction are not dissipated at the highest order. We establish new anisotropic commutator estimates (Lemma \ref{commutator111}) which, together with the constraint $\mathrm{div}\mathbf h=0$, guarantee that every undissipated top-order term is multiplied by a factor controlled by the dissipation at a lower order. A hierarchy of time-weighted energy estimates then shows that this factor is integrable in time, and the highest-order energy is bounded uniformly in time. This is in sharp contrast with \cite{Wu-Zhu}, where the highest-order energy is allowed to grow in time.
\item[(4)] \emph{Lower regularity and explicit decay rates.} Although the non-isentropic system is more complicated than the isentropic one, the mechanisms above allow us to reduce the regularity requirement from $N\geq10$ in \cite{Wu-Zhu} to $N\geq5$, and they yield the explicit decay rates stated in Theorem \ref{th:main1} for $\partial_2a$, $\partial_2\mathbf u$, $\mathrm{div}\mathbf u$, $D\theta$, $\partial_2\mathbf h$ and $\Phi$.
\end{itemize}
A more detailed description of the stabilizing mechanisms and  the proof strategy  is given in Subsections 1.4 and 1.5.}
 Here, we consider the  following initial data
$$(\rho,\mathbf u,\Theta,\mathbf H)\big|_{t=0}=(\rho_0,\mathbf u_0,\Theta_0,\mathbf H_0),$$ and  periodic boundary conditions
$$x\in\mathbb T^2:=[0,1]^2.$$
%% ------------------------------------------------------------------
\subsection{Notations}
%% ------------------------------------------------------------------
Let  $C$ denote  a generic  positive constant throughout this paper that may vary at
different places and denote  $"f\leq Cg"$ by  $"f\lesssim g"$. Here $C$ denotes a generic constant that may vary from place to place and depends only on the constants $|\mathbb T|,N,\mu,\lambda,R,c_{\nu},\kappa$ and $s$, where $s$ will be defined in Theorem \ref{th:main1}. We take $L^q(\mathbb T^2), q\geq1$ and $H^k(\mathbb T^2), k\geq0$ for the usual Sobolev spaces
on $\mathbb T^2$ with norms $\|\cdot\|_{L^q}, \|\cdot\|_k$. Here, we denote $\|\cdot\|_0:=\|\cdot\|_{L^2}.$ Without further statements, $\int\cdot dx$ will be
denoted as $\int_{\mathbb T^2}dx$; $\|h(f,g)\|_{X},\|hf,hg\|_{X}$ will be denoted as
$(\|hf\|^2_{X}+\|hg\|^2_{X})^\frac12$; $[U,V]W$ will be denoted as $U(VW)-V(UW)$; $\partial_i$ will be
denoted as $\partial_{x_i}$. For a multi-index $\alpha\in\mathbb N^2,$ we set
$$|\alpha|=\alpha_1+\alpha_2,\ \text{and}\ \partial^{\alpha}=\partial^{\alpha_1}_{1}\partial^{\alpha_2}_2,$$
\subsection{Main Results}
Let $(a,\theta,\mathbf h):=(\rho-1,\Theta-1,\mathbf H-\mathbf e_2)$  with $(a_0,\theta_0,\mathbf h_0):=(\rho_0-1,\Theta_0-1,\mathbf H_0-\mathbf e_2)$, then the system
\eqref{heat-mhd1} can be rewritten as
\begin{align} \label{heat-mhd}
\left\{
\begin{aligned}
&\partial_{t}a+\mathrm{div}(a\mathbf u)+\mathrm{div}u=0,\\
&(1+a)\left(\partial_{t}\mathbf u+(\mathbf u\cdot\nabla)\mathbf u\right)-\mu\Delta \mathbf u-(\mu+\lambda)\nabla\mathrm{div}\mathbf u+\nabla(R\theta+Ra+Ra\theta+\frac12|\mathbf
h|^2+\mathbf h_2)=(\mathbf h\cdot\nabla)\mathbf h+\partial_2\mathbf h, \\
&c_\nu(1+a)(\partial_{t}\theta+(\mathbf u\cdot\nabla)\theta)+R(1+a)(1+\theta)\mathrm{div}\mathbf u-\kappa\Delta\theta=\mathcal Q(\mathbf u), \\
&\partial_t\mathbf h+(\mathbf u\cdot \nabla)\mathbf h
=(\mathbf h\cdot \nabla)\mathbf u+\partial_2\mathbf u-\mathbf h\mathrm{div} \mathbf u-\mathrm{div} \mathbf u\mathbf e_2,\\
&\mathrm{div}\mathbf h=0.
\end{aligned}
\right.
\end{align}

Our  main result on the global well-posedness of  the  system \eqref{heat-mhd}  then reads  as follows.
\begin{theorem}\label{th:main1} Let $N\geq5$ and $s\in(0,1)$. Assume the initial data $(a_0,\mathbf u_0,\theta_0,\mathbf h_0)\in H^{N}(\mathbb T^2)$ with $\mathrm{div}\mathbf h_0=0$ satisfying
\begin{equation}\label{a-theta-bound}
|a_0|,|\theta_0|\leq\frac12,
\end{equation}
\begin{equation}\label{junzhi}
\int\mathbf (\mathbf h_{0})_1dx=0.
\end{equation}
There exists a small constant $\epsilon>0$ depending only on $|\mathbb T|,N,\mu,\lambda,\kappa, c_\nu, R$ and $s$ such that, if
\begin{equation}\label{initial-date-1}
\|a_0,\mathbf u_0,\theta_0,\mathbf h_0\|_{N}\leq\epsilon,
\end{equation}
then  the system \eqref{heat-mhd} admits a unique global solution  $(a,\mathbf u,\theta,\mathbf h)\in C([0,\infty);H^N(\mathbb T^2))$ satisfying
\begin{equation*}
\sup_{t\in[0,\infty)}\|a,\mathbf u,\theta,\mathbf h\|_{N}\lesssim\epsilon,
\end{equation*}
and
\begin{align*}
&\sup_{t\in[0,\infty)}(1+t)^{s+l}\|\partial_2a,\partial_2\mathbf u,D\theta,\partial_2\mathbf h,\mathrm{div}\mathbf u,\Phi\|_{N-l-1}^2\\
&\quad+\int_0^{\infty}(1+t)^{s+l}\big(\|\partial_2\mathbf u,\mathrm{div}\mathbf u,D\theta\|^2_{N-{l}}+\|D\mathbf u,\Phi\|^2_{N-l-1}+\|\partial_2^2a,\partial_2^2\mathbf h\|^2_{N-l-2}\big)dt\lesssim\epsilon^2,
\end{align*}
where $l=0,1,2,\cdots,N-3,$  $\Phi:=R\nabla a(1+\theta)+\nabla \mathbf h_2+\nabla\frac{|\mathbf h|^2}{2}.$
\end{theorem}
\begin{remark} {When the temperature is frozen ($\theta\equiv$ constant), the system \eqref{heat-mhd} formally reduces to the 2D viscous and non-resistive compressible isentropic MHD system studied in \cite{Wu-Zhu}. Compared with \cite{Wu-Zhu}, we exploit new dissipative mechanisms here, namely the good unknown $\Phi$ and its coupling with $\mathrm{div}\mathbf u$ and $\theta$. They enable us to treat the more complicated non-isentropic system and, at the same time, to reduce the regularity requirement on the initial data for the global well-posedness and the decay estimates from $H^{N}(\mathbb T^2)$ with $N\geq10$ to $H^{N}(\mathbb T^2)$ with $N\geq5$.}
\end{remark}

\begin{remark}{In \cite{Wu-Zhu}, the highest-order energy of smooth solutions to the isentropic system is allowed to grow in time. Here, we establish a uniform-in-time bound for the highest-order energy, under lower regularity assumptions and for the full system. In this sense, our result essentially improves upon the one available for the isentropic case.}
\end{remark}

\begin{remark}{The restriction $s\in(0,1)$ and the condition $N\geq5$ come from the proof: $s<1$ is used to close the time-weighted estimate at the level $N-1$ (see Step 2 in the proof of Lemma \ref{decay}), while $N\geq5$ ensures both the time integrability of the decay rate of the low-order energy $E_2$ and the bound $\|D^2\partial_2^2(a,\mathbf h)\|_0^2\lesssim\mathcal D_{N-2}$ needed in the highest-order estimate.}
\end{remark}

\subsection{The stabilizing mechanism at the linear level}
At the linear level, the stabilizing mechanism is determined by the interplay between viscous dissipation, thermal diffusion, and the anisotropic coupling induced by the background magnetic field. For the two-dimensional system with the constant background magnetic field $\mathbf e_2$, the velocity and temperature equations provide direct dissipation, whereas the density and magnetic field equations contain no explicit dissipative terms. Consequently, the main difficulty lies in recovering the missing control of the density and magnetic field through their coupling with the dissipative variables. Normalising $c_\nu=\mu=R=\kappa=1$ and $\mu+\lambda=0$ and for readability, the linearised system reads
\begin{align} \label{heat-mhd-linear}
\left\{
\begin{aligned}
&\partial_{t}a+\mathrm{div}u=0,\\
&\partial_{t}\mathbf u-\Delta \mathbf u+\nabla(\theta+a+\mathbf h_2)=\partial_2\mathbf h, \\
&\partial_{t}\theta+\mathrm{div}\mathbf u-\Delta\theta=0, \\
&\partial_t\mathbf h+\mathrm{div} \mathbf u\mathbf e_2
=\partial_2\mathbf u,\\
&\mathrm{div}\mathbf h=0.
\end{aligned}
\right.
\end{align}
\textbf{Effective dissipation for $\partial_2h$ and $\partial_2a$.}  Firstly, from $\eqref{heat-mhd-linear}_{{2,4}},$ we get a degenerate damped {wave} equation:
\begin{align}\label{linear-h}\left\{\begin{aligned}
&\partial_t^2\mathrm{curl}\mathbf h-\Delta\partial_t\mathrm{curl}\mathbf h-\partial_2^2\mathrm{curl}\mathbf h=\partial_1(\Delta\mathrm{div}\mathbf u-\partial_t\mathrm{div}\mathbf u),\\
&\mathrm{div} \mathbf h=0.\end{aligned}\right.
\end{align}
Then, the left-hand side of \eqref{linear-h} is a degenerate damped wave operator, while the right-hand side is a source term driven by \(\mathbf u\). Hence the dissipation intrinsic to \(\mathbf u\) itself allows us to exploit the dissipation of \(\mathbf h\) in the \(x_2\)-direction. Since the magnetic dissipation is degenerate in the \(x_1\)-direction, we may further examine the properties of the \(x_2\)-derivative of the magnetic field. %Its left-hand side is a degenerate damped wave operator; its right-hand side is a source driven by $\mathbf u$. Thus, by means of the dissipation inherent in $\mathbf u$ itself, we exploit the dissipation of $\mathbf h$ in the $x_2$-direction.
 %Since the magnetic dissipation is degenerate in the $x_1$-direction, we may further consider the properties of the derivative of the magnetic field in the $x_2$-direction. 
Thus,  differentiating the above equation with respect to $x_2$ yields the {same} damped wave equation for $\partial_2\mathbf h:$
\begin{align}\label{linear-h-1}\left\{\begin{aligned}
&\partial_t^2\mathrm{curl}\partial_2\mathbf h-\Delta\partial_t\mathrm{curl}\partial_2\mathbf h-\partial_2^2\mathrm{curl}\partial_2\mathbf h=\partial_1(\Delta\mathrm{div}\partial_2\mathbf u-\partial_t\mathrm{div}\partial_2\mathbf u),\\
&\mathrm{div}\partial_2\mathbf h=0.\end{aligned}\right.
\end{align}
Because $\partial_2\mathbf h$ has zero mean in $x_2$, Poincar\'{e}'s inequality in the $x_2$-variable makes the restoring term $-\partial_2^2$ coercive on $\partial_2\mathbf h$. In this sense, the wave structure is nondegenerate on \(\partial_2\mathbf h\), whereas it yields no control over the modes independent of \(x_2\). It follows from the above equation that the dissipation of $\partial_2\mathbf u$ determines the dissipation of $\partial_2\mathbf h$. Secondly, the equations $\eqref{heat-mhd-linear}_{1,2}$ yield the following damped wave equation for $\partial_2a:$
\begin{align}\label{linear-rho}
\partial_t^2\partial_2a-\Delta\partial_t\partial_2a-\Delta\partial_2a=\Delta(\partial_2\theta+\partial_2\mathbf h_2).
\end{align}
Obviously, the left-hand side of \eqref{linear-rho} is a damped wave operator, while the right-hand side is a source driven by $\partial_2\theta$ and $\partial_2\mathbf h$. This further shows that the dissipation of \(\partial_2 a\) can be recovered from the intrinsic dissipation of \(\partial_2\theta\) in \eqref{heat-mhd-linear}, together with that of \(\partial_2\mathbf h\) derived from \eqref{linear-h-1}.\\
\textbf{A hidden stabilizing mechanism.} Let $d{:=}\mathrm{div}\mathbf u$, $\phi:=a+\mathbf h_2$.  In view of \eqref{heat-mhd-linear} {(note that $\partial_t\mathbf h_2=\partial_2\mathbf u_2-\mathrm{div}\mathbf u$)}, we consider the following equations for $d, \phi$ and $\theta${, in which $\partial_2\mathbf u_2$ is regarded as a source term}:
\begin{align} \label{linear-phi-divu-theta}
\left\{
\begin{aligned}
&\partial_{t}d-\Delta d+\Delta \theta+\Delta\phi=0,\\
&\partial_{t}\theta-\Delta \theta+d=0, \\
&\partial_t\phi+2d={\partial_2\mathbf u_2}.
\end{aligned}
\right.
\end{align}
To investigate the stability structure of the linearized system, we perform a Fourier analysis of the coupled equations for \(d\), \(\theta\), and \(\phi\). For each nonzero spatial frequency \({\mathbf k}\in\mathbb Z^2\), the corresponding homogeneous system can be written as
$$
\frac{d}{dt}
\begin{pmatrix}
\widehat d\\
\widehat\theta\\
\widehat\phi
\end{pmatrix}
=
\begin{pmatrix}
-|{\mathbf k}|^2&|{\mathbf k}|^2&|{\mathbf k}|^2\\
-1&-|{\mathbf k}|^2&0\\
-2&0&0
\end{pmatrix}
\begin{pmatrix}
\widehat d\\
\widehat\theta\\
\widehat\phi
\end{pmatrix}.
$$
The characteristic polynomial of the coefficient matrix is given by
$$
\lambda^3+2|{\mathbf k}|^2\lambda^2
+\bigl(|{\mathbf k}|^4+3|{\mathbf k}|^2\bigr)\lambda
+2|{\mathbf k}|^4=0.
$$
By the Routh-Hurwitz criterion, all roots of this polynomial have strictly negative real parts for every $\mathbf k\neq0$. Hence, each nonzero Fourier mode of the homogeneous system is spectrally stable. This stability results from the coupling among $d$, $\theta$, and $\phi$, through which the diffusive effects acting on $d$ and $\theta$ are transferred to the nondiffusive variable $\phi$. In this sense,  \eqref{linear-phi-divu-theta} shows that the dissipation inherent in $\partial_2u_2$ in equation \eqref{heat-mhd-linear} determines the stability of $d$, $\theta$, and $\phi$.

\subsection{Difficulties and scheme of the proof} The first main difficulty  we encounter here stems from the absence of dissipation or damping terms for $a,\mathbf h$ and the strong nonlinear interactions arising from  the velocity vector field,  the magnetic field and the temperature field.  In particular,  neither the density nor the magnetic field admits a direct diffusive mechanism. Their higher-order regularity must therefore be controlled indirectly via coupling to the velocity field, leveraging the viscous and thermal dissipation available in the system. A second main difficulty is that we aim to establish uniform-in-time estimates for the highest-order energy. This is a highly nontrivial task  since the highest-order energy may grow in the direction where the stabilizing effect of the background magnetic field degenerates. The central challenge lies in constructing an appropriate energy framework that captures the couplings among velocity, magnetic field and temperature, closes the nonlinear estimates at lower regularity, and achieves uniform-in-time control.

 To overcome the first difficulty, we need to  exploit the stabilizing mechanism for $a,\mathbf h$  induced by the background magnetic field \(\mathbf e_2\).  From a simple analysis of \eqref{linear-h-1} and \eqref{linear-rho}, we  observe that magnetic coupling yields a dissipative effect of $a,\mathbf h$  in the \(x_2\)-direction, which is confined to that direction. This motivates us to fully extract this  dissipation by using the energy method. %This motivates us to adopt separate treatments for the \(x_1\)- and \(x_2\)-directions in our energy estimates.
Concretely,  for any $m\geq2$,  we can  establish the following two key estimates:
 \begin{align}
{\frac2R}\frac{d}{dt}&\sum_{|\alpha|=m-1}\int\partial^\alpha\partial_2\mathbf u_2\partial^\alpha\partial^2_2adx+\|D^{m-1}\partial^2_2a\|_0^2\leq C_0\|D^{m+1}\partial_2\mathbf u,D^{m+1}\partial_2\theta\|^2_0\nonumber\\
&\quad+(1+\|a,\theta\|_{m+1})\|a,\mathbf u,\theta,\mathbf h\|_{m+1}\|D^{m-1}\partial_2^2a,D^{m+1}\partial_2\mathbf u,D^{m+1}\mathbf u,D^{m-1}\partial^2_2\mathbf h\|_0^2,\label{10-16-partial2-a}
\end{align}
  and
\begin{align}
2\frac{d}{dt}&\sum_{|\alpha|=m-1}\int\partial^\alpha\partial^2_2\mathbf u\cdot\partial^\alpha\partial_2\mathbf hdx+\|D^{m-1}\partial^2_2\mathbf h\|_0^2\leq C_1\|D^{m+1}\partial_2\mathbf u,D^{m+1}\partial_2\theta\|^2_0\nonumber\\
&\quad+(1+\|a\|_{m+1})\|a,\mathbf u,\mathbf h\|_{m+1}\|D^{m-1}\partial_2^2a,D^{m+1}\partial_2\mathbf u,D^{m+1}\mathbf u,D^{m+2}\theta,D^{m-1}\partial^2_2\mathbf h\|_0^2,\label{10-16-partial2-b}
\end{align}
(see Lemmas \ref{dissipation-rho}-\ref{dissipation-h} for details). However, another principal difficulty we face lies in how to use the resulting dissipative effect on \(a,\mathbf h\) in the \(x_2\)-direction to  control the energy associated with \(x_1\)-directional modes. Our idea here  mainly hinges on  new anisotropic commutator estimates (see Lemma   \ref{commutator111}) and  the {structure} of the equations, for example, $\mathrm{div}\mathbf h=0$, which enables shift of derivative from the space variable \(x_1\)-direction  to the space variable \(x_2\)-direction. %More specifically,
More precisely, by employing  analysis tools and  the {structure} of the
equations,  we can establish the following  higher-order energy estimates with $m\geq2$ for the solution  $(a,\mathbf u,\theta,{\mathbf h})$  according to $x_1$-direction and  $x_2$-direction, respectively
\begin{align}
&\frac{d}{dt}\|\sqrt{\frac{R(1+\theta)}{1+a}}D^m\partial_1a,\sqrt{\rho}D^m\partial_1\mathbf u,\sqrt{c_\nu\rho}D^m\partial_1\theta,D^m\partial_1\mathbf h\|^2_0+2\min\{\mu,\kappa\}\|D^{m+1}\partial_1\mathbf u,D^{m+1}\partial_1\theta\|_0^2\nonumber\\
&\lesssim(1+\|a,\mathbf u,\theta,\mathbf h\|_{m+1})^2\|D^2\partial_2^2a,D^3\mathbf u,D^{m}\mathrm{div}\mathbf u,D^m\partial_2\mathbf u,D^4\theta,D^2\partial_2^2\mathbf h\|_0\|D^{m+1}a,D^{m+1}\mathbf u,D^{m+1}\mathbf h\|_0^2\nonumber\\
&\quad+(1+\|a,\mathbf u,\theta,\mathbf h\|_{m+1})^2(\|Da,D\mathbf u,D\mathbf h\|_m+\|\theta\|_{m+1})\nonumber\\
&\qquad{\times}\|{D^{m-1}\partial_2^2a,}D^{m+2}\mathbf u,D^{m+2}\theta,D^{m-1}\partial_2^2\mathbf h,D^2\partial_2^2\mathbf h\|_0^2,\label{1-h-older-x}
\end{align}
and
\begin{align}
&\frac{d}{dt}\|\sqrt{\frac{R(1+\theta)}{1+a}}D^m\partial_2a,\sqrt{\rho}D^m\partial_2\mathbf u,\sqrt{c_\nu\rho}D^m\partial_2\theta,D^m\partial_2\mathbf h\|^2_0+2\min\{\mu,\kappa\}\|D^{m+1}\partial_2\mathbf u,D^{m+1}\partial_2\theta\|_0^2\nonumber\\
&\lesssim(1+\|a,\mathbf u,\theta,\mathbf h\|_{m+1})^3\|a,\mathbf u,\theta,\mathbf h\|_{m+1}\|D^{m-1}\partial_2^2a,D^{m+1}\partial_2\mathbf u,D^{m+1}\mathbf u, D^{m+1}\mathrm{div}\mathbf u,D^{m+2}\theta,D^{m-1}\partial_2^2\mathbf h\|_0^2\nonumber\\
&\quad+\Big|\int{\frac{1+\theta}{1+a}}\partial_2\mathbf u_1\partial_1^{m+1}a\partial_1^m\partial_2adx\Big|+\Big|\int\partial_2\mathbf u_1\partial_1^{m+1}\mathbf h\partial_1^m\partial_2\mathbf hdx\Big|.\label{2-h-older-y}
\end{align}
Obviously, the terms \(D^{m-1}{\partial_2^2} a\) and \(D^{m-1}{\partial_2^2}\mathbf h\) appearing on the right-hand sides of \eqref{1-h-older-x} and \eqref{2-h-older-y} can be absorbed into the left-hand sides of \eqref{10-16-partial2-a} and \eqref{10-16-partial2-b}. For the remaining parts,
we have to exploit the dissipation of {$(D^{m+1}\mathrm{div}\mathbf u,D^{m+2}\theta,D^{m+1}\mathbf u)$}.
To address the  difficulty,  we introduce a suitably defined total pressure $\Phi:=R\nabla a(1+\theta)+\nabla \mathbf h_2+\nabla\frac{|\mathbf h|^2}{2}.$  %Here we provide an explanation.
For the term $R\nabla a+R\nabla a\theta+\nabla\mathbf h_2+\nabla\frac{|\mathbf h|^2}{2}$ stated in $\eqref{heat-mhd}_2$, the four terms $R\nabla a,R\nabla a\theta,\nabla\mathbf h_2$ and $\nabla\frac{|\mathbf h|^2}{2}$ do not contribute any dissipation individually. However,
the choice of \(\Phi\) relying on all combinations of these nondissipative terms  shall present favorable dissipative effect  in all directions via coupling with the equation of $\mathrm{div}\mathbf u$
\begin{align*}
\left\{
\begin{aligned}
&\partial_{t}\Phi+(\mathbf u\cdot\nabla)\Phi+(R+1)\nabla\mathrm{div}u-\nabla\partial_2\mathbf u_2=\cdots,\\
&\partial_{t}\mathrm{div}\mathbf u-(2\mu+\lambda)\Delta\mathrm{div}\mathbf u+R\Delta\theta+\mathrm{div}\Phi=\cdots,
\end{aligned}
\right.
\end{align*}
(see \eqref{Phi-div-mhd} for details). Here,  the core idea of this process originates precisely from the stability structure of the linear system \eqref{linear-phi-divu-theta}.    Based  on this crucial observation, we  finally obtain the dissipative estimates of {$D^{m+1}\mathrm{div}\mathbf u,D^{m+2}\theta, D^{m+1}\mathbf u$ and $D^m\Phi$} (see Lemmas \ref{dissipation-divu}-\ref{dissipation-Du} for details).

For the second difficulty, combining with the estimates  of {$D^{m+1}\mathrm{div}\mathbf u,D^{m+2}\theta, D^{m+1}\mathbf u$ and $D^m\Phi$} developed above with  \eqref{10-16-partial2-a} \eqref{10-16-partial2-b} and \eqref{2-h-older-y}, we can  derive the higher-order ($m\geq2$) energy estimates as follows
\begin{align}\label{2-h-older-yy}
\frac{d}{dt}E_m+2\mathcal D_m&\lesssim\|D\partial_2\mathbf h\|^2_0\|D^{m+2}\mathbf u\|^2_0\nonumber\\
&\quad+\big|\int{\frac{1+\theta}{1+a}}\partial_2\mathbf u_1\partial_1^{m+1}a\partial_1^m\partial_2adx\big|+\big|\int\partial_2\mathbf u_1\partial_1^{m+1}\mathbf h\partial_1^m\partial_2\mathbf hdx\big|
\end{align}
with
\begin{align*}&E_m\approx\|D^m\partial_2a,D^m\partial_2\mathbf u,D^m\mathbf u,D^m\mathrm{div}\mathbf u, D^m\partial_2\mathbf h,D^m\partial_2\theta,D^{m+1}\theta,D^{m}\Phi\|_0^2,
\end{align*}
\begin{align*}&\mathcal D_m\approx\|D^{m-1}\partial^2_2 a,D^{m+1}\partial_2\mathbf u,D^{m+1}\mathrm{div}\mathbf u,D^{m+1}\mathbf u,D^{m+1}\partial_2\theta,D^{m+2}\theta,D^{m-1}\partial^2_2\mathbf h,D^{m}\Phi\|_0^2.
\end{align*}
Employing  the smallness assumption on $\|a,\mathbf u,\theta,\mathbf h\|_N$, we conclude that for any $2\leq m\leq N-2$,
\begin{equation*}
\frac{d}{dt}E_m+\mathcal D_m\leq0,
\end{equation*}
which together with the  Gagliardo-Nirenberg interpolation inequality, \eqref{dengjia-1} and Poincar\'{e}'s inequality yields that
\begin{equation}\label{2026-10-17-5-1}
E_m(t)\lesssim  \big( t+1\big)^{-(N-m-1)}.
\end{equation}
However, this decay estimate is not sufficient for us to establish  the uniform estimates for the highest-order derivatives because the time integrability   of $\mathcal D_{N-2}$  is required (see \eqref{2026-10-18-7-1} below). To resolve this issue,  using \eqref{2026-10-17-5-1},  we need  further exploit the following decay estimates
\begin{align}
&\sup_{t\in[0,T)}(1+t)^{s+l}E_{N-1-l}+\int_0^T(1+t)^{s+l}\mathcal D_{N-1-l}dt<\infty,\label{2026-10-19-3-1}
\end{align}
where $s\in(0,1)$, and $l=0,1,2\cdots,N-3$ (see Lemma \ref{decay} for details).
 %It then  follows from  \eqref{2026-10-19-3-1} that $\mathcal D_{N-2}$  is integrable in time, which guarantees the uniform boundedness of the highest-order derivatives of the solution.

Furthermore, taking $m=N-1$ in \eqref{1-h-older-x} and \eqref{2-h-older-yy}, and adding them together,  we can  derive the highest-order energy estimates as follows
\begin{align}\label{2026-10-18-7-1}
\frac{d}{dt}\mathcal E+c_0\|D(\mathbf u,\theta)\|^2_{N}\lesssim(\mathcal D_{N-2})^\frac12\mathcal E,
\end{align} where $\mathcal E \approx\|a,\mathbf u,\theta,\mathbf h\|^2_{N}.$
This along with the time integrability   of $\mathcal D_{N-2}$ in \eqref{2026-10-19-3-1} and Gronwall's inequality, guarantees the uniform boundedness of the highest-order energy estimates. Finally, applying  the
bootstrapping argument
%the continuity method
gives uniform \emph{a priori} estimates (see Proposition \ref{Uniform-boundedness}). The proof of Theorem \ref{th:main1} is then completed at the end of Section \ref{Theorem-th:main2}.

\section{Preliminaries}
This section reviews several analysis tools including the several calculus inequalities that were needed in this paper. The following standard Moser-type product and commutator estimates on $\mathbb{T}^2$ will be used throughout the paper (see, e.g., \cite{Kato-Ponce,Kenig-Ponce-Vega}).
\begin{lemma}\label{Hm-estimate-1}
Let $m\geq1$ be an integer, $|\alpha|=m,$ and $f,g$ be two smooth functions on $\mathbb T^2$. There exists a positive constant $K_1$  such that
\begin{align} \|\partial^\alpha(fg)\|_{L^p} &\leq K_1 \big(\|f\|_{L^{p_1}}\|D^m g\|_{L^{p_2}} + \|D^m f\|_{L^{p_3}}\|g\|_{L^{p_4}}\big), \label{product}\\ \|[\partial^\alpha,f]g\|_{L^p} &\leq K_1\big(\|\nabla f\|_{L^{p_1}}\|D^{m-1}g\|_{L^{p_2}} +\|D^mf\|_{L^{p_3}} \|g\|_{L^{p_4}}\big), \label{commutator} \end{align}
with $p_2,p_3\in(1,\infty)$ and $\frac{1}{p}=\frac{1}{p_1}+\frac{1}{p_2}=\frac{1}{p_3}+\frac{1}{p_4}$.
\end{lemma}

With the above results in hand, we have  the following  anisotropic
commutator estimates.
\begin{lemma}\label{commutator111}
Let $m\geq1$ be an integer and $|\alpha|=m.$ Let $f=(f_1, f_2)$ and $g$ be two smooth functions on $\mathbb T^2$.  There is a positive constant $K_2$ such that
\begin{align}
\|[\partial^\alpha\partial_1,f\cdot\nabla]g\|_0\leq &K_2(\|D^2\mathrm{div}f,D^2\partial_2f\|_0\|D^{m+1}g\|_{0}+\|D^{m}\mathrm{div}f,D^{m}\partial_2f\|_0\|D^3g\|_{0})\nonumber\\
&+K_2(\|D^3f\|_{0}\|D^{m}\partial_2g\|_{0}+\|D^{m+1}f\|_0\|\partial_2D^2g\|_{0}){,}\label{commutator-1}\\
\|[\partial^\alpha\partial_2,f\cdot\nabla]g\|_0\leq & K_2(\|D\partial_2\mathrm{div}f,D\partial_2^2f\|_{0}\|D^{m+1}g\|_{0}+\|D^{m-1}\partial_2\mathrm{div}f,D^{m-1}\partial^2_2f\|_0\|D^3g\|_{0})\nonumber\\
&+{K_2}(\|D^2\mathrm{div}f,D^2\partial_2f\|_{0}\|D^{m}\partial_2g\|_{0}+\|D^{m}\mathrm{div}f,D^{m}\partial_2f\|_{0}\|D^2\partial_2g\|_{0})\nonumber\\
&+{K_2}(\|D^3f\|_{0}\|D^{m-1}\partial^2_2g\|_{0}+\|D^{m+1}f\|_{0}\|D\partial^2_2g\|_{0}).\label{commutator-2}
\end{align}
\end{lemma}
\begin{proof}
Since $[\partial^\alpha\partial_i,f\cdot\nabla]g=\sum_{j=1}^2{\big(}\partial^\alpha(\partial_if_j\partial_jg)+[\partial^\alpha,f_j]\partial_j\partial_ig{\big)}$ for any $i=1,2$, we infer by Lemma \ref{Hm-estimate-1} that
\begin{align*}
\|[\partial^\alpha\partial_i,f\cdot\nabla]g\|_0\leq& K_1(\|\partial_if_1\|_{L^\infty}\|D^{m+1}g\|_{0}+\|D^{m}\partial_if_1\|_0\|Dg\|_{L^\infty})\\
&+K_1(\|\partial_if\|_{L^\infty}\|D^{m}\partial_2g\|_{0}+\|D^{m}\partial_if\|_0\|\partial_2g\|_{L^\infty})\\
&+K_1(\|Df_1\|_{L^\infty}\|D^{m}\partial_ig\|_{0}+\|D^{m}f_1\|_{L^4}\|D\partial_ig\|_{L^4})\\
&+K_1(\|Df_2\|_{L^\infty}\|D^{m-1}\partial_2\partial_ig\|_{0}+\|D^{m}f_2\|_{L^4}\|\partial_2\partial_ig\|_{L^4}).
\end{align*}
This together with Sobolev's embedding inequality, $\|D(\cdot)\|_{L^\infty}\lesssim\|D^3(\cdot)\|_0,\ \|D(\cdot)\|_{L^4}\lesssim\|D^2(\cdot)\|_0$ and the fact, $\|Df_1\|_0\leq\|\partial_2f\|_0+\|\mathrm{div}f\|_0$, yield the desired conclusion. This completes the proof of Lemma \ref{commutator111}.
\end{proof}
\begin{lemma}(\cite{Wu-Zhai} Lemma 2.6)\label{nolinear-estimate}
Let $m>0$ and $f\in H^m(\mathbb T^2)\cap L^\infty(\mathbb T^2).$ Assume that $F$ is a smooth function on $\mathbb R$ with $F(0)=0$. Then,
$$\|F(f)\|_{m}\leq K_3(1+\|f\|_{L^\infty})^{[m]+1}\|f\|_m,$$
where the positive constant $K_3$ depends on $|\mathbb T^2|,m$ and {$\sup_{k\leq[m]+2,\,|t|\leq\|f\|_{L^\infty}}|F^{(k)}(t)|.$}
\end{lemma}

\section{A priori estimates }\label{Theorem-th:main1}
In this section, we first derive \emph{a priori} estimates for smooth  solutions to the  system \eqref{heat-mhd}. For this purpose, we assume throughout that
$(a,\mathbf{u},\theta,\mathbf{h})$ is a smooth solution of the  system \eqref{heat-mhd} on $\mathbb{T}^2\times[0,T]$ for some $T>0$, with $(a,\mathbf{u},\theta,\mathbf{h})\in
C([0,T];H^{N}(\mathbb{T}^2))$.% At the end of this section, we complete the proof of Theorem \ref{th:main1}.
\subsection{Basic energy estimate} %First, we present the basic energy estimate.
\begin{lemma}\label{priori-basic} Let $(a,\mathbf u,\theta,\mathbf h)$ be smooth solutions to  the system $\eqref{heat-mhd}$
 on $\mathbb T^2\times[0,T)$ for some $T>0$,  then there holds,
\begin{equation}\label{basic-1}
\frac{d}{dt}\|\sqrt{R}a,\sqrt{\rho}\mathbf u,\sqrt{c_\nu\rho}\theta,\mathbf h\|_0^2+2\mu\|D\mathbf u\|_0^2+2\kappa\|D\theta\|_0^2\lesssim(1+\|a\|_{2})\|Du\|_{0}\|a,\theta,D\mathbf u\|^2_2.
\end{equation}
\end{lemma}
\begin{proof}
Performing a basic $L^2$ energy estimate for equation $\eqref{heat-mhd}_{{1,2,3,4}}$ with $(Ra,\mathbf u,\theta,\mathbf h)$, integration by parts, using
$\eqref{heat-mhd}_{5}$ and then using the following facts:
$$\int(\mathbf h\cdot\nabla)\mathbf h\cdot \mathbf u+(\mathbf h\cdot\nabla)\mathbf u\cdot \mathbf hdx=0,$$
$$\int\frac12(\mathbf u\cdot\nabla)(|\mathbf h|^2)+(\mathbf u\cdot\nabla)\mathbf h\cdot \mathbf h+\mathrm{div}\mathbf u|\mathbf h|^2dx=0,$$
$$\int\mu|\nabla\mathbf u|^2+(\mu+\lambda)(\mathrm{div}\mathbf u)^2{dx}=\int{\mathcal Q}( \mathbf u)dx,$$
we obtain by H\"{o}lder's inequality that
\begin{align*}
&\frac{1}{2}\frac{d}{dt}\|\sqrt{R}a,\sqrt{\rho}\mathbf u,\sqrt{c_\nu\rho}\theta,\mathbf h\|_0^2+\mu\|D\mathbf u\|_0^2+(\mu+\lambda)\|\mathrm{div}\mathbf u\|_0{^2}+\kappa\|D\theta\|_0^2\\
&=-\frac R2\int\mathrm{div}\mathbf u|a|^2dx-R\int (1+a)\mathrm{div}\mathbf u|\theta|^2dx+\int\mathcal Q(\mathbf u)\theta dx\\
&\lesssim(1+\|a\|_{L^\infty})\|Du\|_{L^2}\|a,\theta,D\mathbf u\|^2_{L^4}.
\end{align*}
This together with  Sobolev's embedding inequality yields \eqref{basic-1}. This completes the proof of Lemma \ref{priori-basic}.
\end{proof}
\subsection{Higher-order energy estimates} In this subsection, we  shall construct  some higher-order energy estimates for smooth solution  $(a,\mathbf u,\theta,\mathbf h)$ to  the system $\eqref{heat-mhd}$ according to $x_1$-direction and  $x_2$-direction, respectively.
\begin{lemma}\label{priori-vt-thetat} Let $(a,\mathbf u,\theta,\mathbf h)$ be smooth solutions to  the system $\eqref{heat-mhd}$  on $\mathbb T^2\times[0,T)$ for some $T>0$,  and suppose that \eqref{junzhi} holds. Then,
\begin{align}
&\|\mathbf h_1\|_1\lesssim\|\partial_2\mathbf h\|_0\label{h-junzhi}.
%&\|\partial_t\mathbf u\|_1\lesssim\|D\partial_t\mathbf u\|_0+\|\mathbf u\|_0\|D\mathbf u\|_0,\label{u-t}\\
%&\|\partial_t\mathbf \theta\|_1\lesssim\|D\partial_t\mathbf \theta\|_0+\|\mathbf u,\theta\|_1\|D\mathbf u,D\theta\|_0.\label{theta-t}
\end{align}
\end{lemma}
\begin{proof}
It follows from the equality $\eqref{heat-mhd}_{4}$ and \eqref{junzhi} that
 $$\int\mathbf h_1dx=\int\mathbf (\mathbf h_{0})_1dx=0,$$
which together with Poincar\'{e}'s inequality and $\mathrm{div}\mathbf h=0$ implies \eqref{h-junzhi}. %It remain to proof \eqref{u-t}.
%On the one hand, one has by Poincar\'{e} and H\"{o}lder inequality and \eqref{density-1}
%\begin{align*}
%\|\partial_t\mathbf u\|_{0}\leq C\|D\partial_t\mathbf u\|_0+\left|\int\partial_t\mathbf udx\right|\leq C\|D\partial_t\mathbf u\|_0+\left|\int(1+a)\partial_t\mathbf udx\right|+\|a\|_{0}\|\partial_t\mathbf u\|_0
%\end{align*}
%On the other hand, integration $\eqref{heat-mhd}_{1}$ with respect to $x$ and using \eqref{density-1} and H\"{o}lder inequality yield
%$$\left|\int(1+a)\partial_t\mathbf udx\right|=\left|\int(1+a)(\mathbf u\cdot\nabla)\mathbf udx\right|\lesssim\|\mathbf u\|_0\|D\mathbf u\|_0.$$
%Thus, using the above two inequalities and  \eqref{density-1}, yields
%$$\|\partial_t\mathbf u\|_{0}\leq C\|D\partial_t\mathbf u\|_0+\|\mathbf u\|_0\|D\mathbf u\|_0+\frac12\|\partial_t\mathbf u\|_0,$$
%and then this leads to \eqref{u-t}. Now we prove \eqref{theta-t}. Similar it follows Poincar\'{e} and H\"{o}lder inequality and \eqref{density-1} that
%\begin{align*}
%\|\partial_t\theta\|_{0}\leq C\|D\partial_t\theta\|_0+\left|\int(1+a)\partial_t\theta dx\right|+\frac12\|\partial_t\theta\|_0
%\end{align*}
%One has by $\eqref{heat-mhd}_{4}$, H\"{o}lder inequality and \eqref{density-1}
%$$\left|\int(1+a)\partial_t\theta dx\right|\lesssim\|\mathbf u,\theta\|_1\|\nabla(\mathbf u,\theta)\|_0.$$
%Thus, \eqref{theta-t} follows by the above two inequalities.
This completes the proof of Lemma \ref{priori-vt-thetat}.\end{proof}
\begin{lemma}\label{H-priori-v} Let $(a,\mathbf u,\theta,\mathbf h)$ be smooth solutions to the system  $\eqref{heat-mhd}$  on $\mathbb T^2\times[0,T]$ for some $T>0$, and suppose that
\begin{align}
&|a|<\frac12\label{density-1},\\
&|\theta|<\frac12,\label{theta-1}
\end{align}
and \eqref{junzhi} hold. Then,  for any $m\geq2$,
\begin{align}
&\frac{d}{dt}\|\sqrt{\frac{R(1+\theta)}{1+a}}D^m\partial_1a,\sqrt{\rho}D^m\partial_1\mathbf u,\sqrt{c_\nu\rho}D^m\partial_1\theta,D^m\partial_1\mathbf h\|^2_0+2\min\{\mu,\kappa\}\|D^{m+1}\partial_1\mathbf u,D^{m+1}\partial_1\theta\|_0^2\nonumber\\
&\lesssim(1+\|a,\mathbf u,\theta,\mathbf h\|_{m+1})^2\|D^2\partial_2^2a,D^3\mathbf u,D^{m}\mathrm{div}\mathbf u,D^m\partial_2\mathbf u,D^4\theta,D^2\partial_2^2\mathbf h\|_0\|D^{m+1}a,D^{m+1}\mathbf u,D^{m+1}\mathbf h\|_0^2\nonumber\\
&\quad+(1+\|a,\mathbf u,\theta,\mathbf h\|_{m+1})^2(\|Da,D\mathbf u,D\mathbf h\|_m+\|\theta\|_{m+1})\nonumber\\
&\qquad{\times}\|{D^{m-1}\partial_2^2a,}D^{m+2}\mathbf u,D^{m+2}\theta,D^{m-1}\partial_2^2\mathbf h,D^2\partial_2^2\mathbf h\|_0^2,\label{1-h-older}
\end{align}
and
\begin{align}
&\frac{d}{dt}\|\sqrt{\frac{R(1+\theta)}{1+a}}D^m\partial_2a,\sqrt{\rho}D^m\partial_2\mathbf u,\sqrt{c_\nu\rho}D^m\partial_2\theta,D^m\partial_2\mathbf h\|^2_0+2\min\{\mu,\kappa\}\|D^{m+1}\partial_2\mathbf u,D^{m+1}\partial_2\theta\|_0^2\nonumber\\
&\lesssim(1+\|a,\mathbf u,\theta,\mathbf h\|_{m+1})^3\|a,\mathbf u,\theta,\mathbf h\|_{m+1}\|D^{m-1}\partial_2^2a,D^{m+1}\partial_2\mathbf u,D^{m+1}\mathbf u,D^{m+1}\mathrm{div}\mathbf u,D^{m+2}\theta,D^{m-1}\partial_2^2\mathbf h\|_0^2\nonumber\\
&\quad+\Big|\int{\frac{1+\theta}{1+a}}\partial_2\mathbf u_1\partial_1^{m+1}a\partial_1^m\partial_2adx\Big|+\Big|\int\partial_2\mathbf u_1\partial_1^{m+1}\mathbf h\partial_1^m\partial_2\mathbf hdx\Big|.\label{2-h-older}
\end{align}
\end{lemma}
\begin{proof}
Given any $i=1,2$ and multi-indices $\alpha$ with $|\alpha|=m$,  applying $\partial^{\alpha}\partial_i$ to \eqref{heat-mhd}$_{1,2,3,4}$, testing against $(\frac{R(1+\theta)}{1+a}\partial^\alpha\partial_ia,\partial^{\alpha}\partial_i\mathbf u,\partial^{\alpha}\partial_i\theta,\partial^{\alpha}\partial_i\mathbf h),$ and then using integration by parts and  $\eqref{heat-mhd}_{1,5}$, we conclude that 
\begin{align}
&\frac{1}{2}\frac{d}{dt}\|\sqrt{\frac{R(1+\theta)}{1+a}}\partial^\alpha\partial_ia,\sqrt{\rho}\partial^\alpha\partial_i\mathbf u,\sqrt{c_\nu\rho}\partial^\alpha\partial_i\theta,\partial^\alpha\partial_i\mathbf h\|^2_0\nonumber\\
&\quad+\|\sqrt{\mu}D\partial^\alpha\partial_i\mathbf u,\sqrt{\mu+\lambda}\partial^\alpha\partial_i\mathrm{div}\mathbf u,\sqrt{\kappa}D\partial^\alpha\partial_i\theta\|_0^2\nonumber\\
&=\frac{R}{2}\int\Big(\partial_t\big(\frac{1+\theta}{1+a}\big)+\mathrm{div}\big(\frac{1+\theta}{1+a}\mathbf u\big)\Big)|\partial^\alpha\partial_ia|^2dx-R\int{\frac{1+\theta}{1+a}}[\partial^\alpha\partial_i,\mathbf u\cdot\nabla]a\partial^\alpha\partial_iadx\nonumber\\
&\quad-R\int\frac{1+\theta}{1+a}[\partial^\alpha\partial_i,a]\mathrm{div}\mathbf u\partial^\alpha\partial_iadx+R\int\Big([\partial^\alpha,a]\partial_i\theta+[\partial^\alpha,\theta]\partial_ia\Big)\partial^\alpha\partial_i\mathrm{div}\mathbf udx\nonumber\\
&\quad-R\int[\partial^\alpha\partial_i,a]\mathrm{div}\mathbf u\partial^\alpha\partial_i\theta dx-\int[\partial^\alpha\partial_i,a]\partial_t\mathbf u\cdot \partial^\alpha\partial_i\mathbf udx-\int[\partial^\alpha\partial_i,(1+a)\mathbf u\cdot\nabla]\mathbf u\cdot\partial^\alpha\partial_i\mathbf udx\nonumber\\
&\quad+\sum_{j=1}^2\int[\partial^\alpha,\mathbf h_j]\partial_i\mathbf h_j\partial^\alpha\partial_i\mathrm{div}\mathbf udx
+\int[\partial^\alpha\partial_i,\mathbf h\cdot\nabla]\mathbf h\cdot\partial^\alpha\partial_i\mathbf udx-c_\nu\int[\partial^\alpha\partial_i,a]\partial_t\theta\cdot \partial^\alpha\partial_i\theta dx\nonumber\\
&\quad-c_\nu\int[\partial^\alpha\partial_i,(1+a)\mathbf u\cdot\nabla]\theta\partial^\alpha\partial_i\theta dx-R\int\partial^\alpha\partial_i((1+a)\theta\mathrm{div}\mathbf u)\partial^\alpha\partial_i\theta dx\nonumber\\
&\quad+\int\partial^\alpha\partial_i\mathcal Q(\nabla\mathbf u)\partial^\alpha\partial_i\theta dx-\int\partial^\alpha\partial_i((\mathbf u\cdot\nabla)\mathbf h)\partial^\alpha\partial_i\mathbf hdx\nonumber\\
&\quad+\int[\partial^\alpha\partial_i,\mathbf h\cdot\nabla]\mathbf u\cdot\partial^\alpha\partial_i\mathbf hdx-\sum_{j=1}^2\int[\partial^\alpha\partial_i,\mathbf h_j]\mathrm{div}\mathbf u\partial^\alpha\partial_i\mathbf h_jdx\nonumber\\&\triangleq\sum_{j=1}^{16}I_j^{i},\label{zuigao-estiamte}
\end{align}
where we have used the following facts
\begin{align*}
R&\int\frac{1+\theta}{1+a}\partial^\alpha\partial_i((1+a)\mathrm{div}\mathbf u)\partial^\alpha\partial_ia+\nabla\partial^\alpha\partial_i((1+a)(1+\theta))\cdot\partial^\alpha\partial_i\mathbf u\\
&\quad+\partial^\alpha\partial_i((1+a)\mathrm{div}\mathbf u)\partial^\alpha\partial_i\theta dx\\
&=R\int\frac{1+\theta}{1+a}[\partial^\alpha\partial_i,a]\mathrm{div}\mathbf u\partial^\alpha\partial_iadx-R\int\Big([\partial^\alpha,a]\partial_i\theta+[\partial^\alpha,\theta]\partial_ia\Big)\partial^\alpha\partial_i\mathrm{div}\mathbf udx\\
&\quad+R\int[\partial^\alpha\partial_i,a]\mathrm{div}\mathbf u\partial^\alpha\partial_i\theta dx,
\end{align*}
\begin{align*}
\int&\partial^\alpha\partial_i(\mathbf h\cdot\nabla\mathbf h)\cdot\partial^\alpha\partial_i\mathbf u+\partial^\alpha\partial_i(\mathbf h\cdot\nabla\mathbf u)\cdot\partial^\alpha\partial_i\mathbf hdx\\
&=\int[\partial^\alpha\partial_i,\mathbf h\cdot\nabla]\mathbf h\cdot\partial^\alpha\partial_i\mathbf u+[\partial^\alpha\partial_i,\mathbf h\cdot\nabla]\mathbf u\cdot\partial^\alpha\partial_i\mathbf hdx,
\end{align*}and
\begin{align*}\frac12&\int\partial^\alpha\partial_i\mathbf u\nabla\partial^\alpha\partial_i(|\mathbf h|^2)+\partial^\alpha\partial_i(\mathbf h\mathrm{div}\mathbf u)\cdot\partial^\alpha\partial_i\mathbf hdx\\
&=\sum_{j=1}^2\int-[\partial^\alpha,\mathbf h_j]\partial_i\mathbf h_j\partial^\alpha\partial_i\mathrm{div}\mathbf u+[\partial^\alpha\partial_i,\mathbf h_j]\mathrm{div}\mathbf u\partial^\alpha\partial_i\mathbf h_jdx.
\end{align*}
In what follows, we shall  successively bound  all the above terms $I_j^i(j=1,2\cdots,16)$.\\
\textbf{Estimate of $I_1^i$.} It follows from $\eqref{heat-mhd}_{1,3}$ that
$$\partial_t\big(\frac{1+\theta}{1+a}\big)+\mathrm{div}\big(\frac{1+\theta}{1+a}\mathbf u\big)=\frac{(2c_\nu-R)(1+\theta)}{c_\nu(1+a)}\mathrm{div}\mathbf u+\frac{\mathcal Q(\mathbf u)+\kappa\Delta\mathbf\theta}{c_\nu(1+a)^2}.$$
Thus, using \eqref{density-1}, \eqref{theta-1} and  H\"{o}lder's and  Young's inequalities, Sobolev's embedding inequality  yields
$$I_1^i\lesssim(1+\|D\mathbf u\|_{L^\infty})\|D\mathbf u,D^2\theta\|_{L^\infty}\|\partial^\alpha\partial_ia\|_0^2\lesssim(1+\|D^3\mathbf u\|_{0})\|D^3\mathbf u,D^3D\theta\|_{0}\|D^{m+1}a\|_0\|D^{m-1}\partial_i^2a\|_0,$$
for any $i=1,2$. Here,  we have used the following fact
\begin{equation}\label{part-1}
\|D\partial_i(\cdot)\|_0^2\lesssim\|D^2(\cdot)\|_0\|\partial_i^2(\cdot)\|_0.
\end{equation}
\textbf{Estimates of $I_2^i$ and $I_{14}^i$.}
When $i=1$,   it then  follows from H\"{o}lder's inequality, Sobolev's embedding inequality{, \eqref{density-1}, \eqref{theta-1}} and \eqref{commutator-1} that
\begin{align*}
I_2^1\lesssim\|[\partial^\alpha\partial_1,\mathbf u\cdot\nabla]a\|_0\|\partial^\alpha\partial_1a\|_0\lesssim\|D^3\mathbf u,D^2\partial_2^2a,D^m\mathrm{div}\mathbf u,D^m\partial_2\mathbf u\|_0\|D^{m+1}a,D^{m+1}\mathbf u\|_0^2.
\end{align*}
When $i=2,$  to estimate $I_2^2$, we further divide it into three terms,
\begin{align*}
I_2^2&=-R\int{\frac{1+\theta}{1+a}}\left([\partial^\alpha,\mathbf u\cdot\nabla]\partial_2a+\partial^\alpha(\partial_2\mathbf u_2\partial_2a)\right)\partial^{\alpha}\partial_2adx-R\int{\frac{1+\theta}{1+a}}[\partial^\alpha,\partial_2\mathbf u_1]\partial_1a\partial^{\alpha}\partial_2adx\\
&\quad-R\int{\frac{1+\theta}{1+a}}\partial_2\mathbf u_1\cdot\partial^{\alpha}\partial_1a\partial^{\alpha}\partial_2adx\\&\triangleq I_{2,1}^2+I_{2,2}^2+I_{2,3}^2.
\end{align*}We bound term by term above in what follows.
First, for $I^2_{2,1},$   employing  H\"{o}lder's and  Young's inequalities, Sobolev's embedding inequality, Lemma \ref{Hm-estimate-1} and \eqref{part-1},   we conclude that
\begin{align*}
I^2_{2,1}&\lesssim(\|D\mathbf u\|_{L^\infty}\|D^m\partial_2a\|_0+\|D^m\mathbf u\|_{L^4}\|D\partial_2a\|_{L^4}+\|D^{m+1}\mathbf u\|_0\|\partial_2 a\|_{L^\infty})\|D^{m}\partial_2a\|_0\\
&\lesssim\|D^{m+1}\mathbf u\|_0\|D^m\partial_2a\|_0^2\lesssim\|D^{m+1}a\|_{{0}}\|D^{m+1}\mathbf u,D^{m-1}\partial^2_2a\|^2_0.
\end{align*}
For $I_{2,2}^2,$ due to $|\alpha|>0$, we may choose a component $\alpha_k (k\in\{1,2\})$ of the multi-index $\alpha$ such that $\alpha_k>0$. Then, choosing the corresponding unit multi-index $\mathbf e_k$,  and employing  integration by parts, H\"{o}lder's and  Young's inequalities, Sobolev's embedding inequality, Poincar\'{e}'s inequality, and Lemma \ref{Hm-estimate-1}, we deduce that
\begin{align*}
I_{2,2}^2&=R\int\partial^{\mathbf e_k}{\Big(\frac{1+\theta}{1+a}}[\partial^\alpha,\partial_2\mathbf u_1]\partial_1a{\Big)}\partial^{\alpha-\mathbf e_k}\partial_2adx\\
&=R\int{\frac{1+\theta}{1+a}}\left([\partial^{\alpha+\mathbf e_k},\partial_2\mathbf u_1]\partial_1a-\partial^{\mathbf e_k}{\partial_2}\mathbf u_1\partial^\alpha\partial_1a\right)\partial^{\alpha-\mathbf e_k}\partial_2adx\\
&\quad{+R\int\partial^{\mathbf e_k}\Big(\frac{1+\theta}{1+a}\Big)[\partial^\alpha,\partial_2\mathbf u_1]\partial_1a\,\partial^{\alpha-\mathbf e_k}\partial_2adx}\\
&\lesssim{(1+\|Da,D\theta\|_{L^\infty})}(\|D^{m+1}\partial_2\mathbf u\|_0\|Da\|_{L^\infty}+\|D{\partial_2}\mathbf u\|_{L^\infty}\|D^{m+1}a\|_{0})\|D^{m-1}\partial_2a\|_0\\
&\lesssim{(1+\|a,\theta\|_{m+1})}\|D^{m+1}a\|_0\|D^{m+1}\partial_2\mathbf u,D^{m+1}\mathbf u,D^{m-1}\partial^2_2a\|^2_0.
\end{align*}
It should be noted that,  $\mathbf e_k$ is used consistently whenever it appears  in the subsequent estimates.
If $\alpha_2\neq0$,  it  then follows  from  H\"{o}lder's and  Young's inequalities, Sobolev's embedding inequality, that
$$I_{2,3}^2\lesssim\|D\mathbf u\|_{L^\infty}\|D^{m+1}a\|_{0}\|D^{m-1}\partial_2^2a\|_0\lesssim\|D^{m+1}a\|_0\|D^{m+1}\mathbf u,D^{m-1}\partial^2_2a\|^2_0.$$
{If $\alpha_2=0$, i.e. $\alpha=(m,0)$, then $I_{2,3}^2$ is kept as it stands.} Thus, combining with  both cases, we get
$$I_{2,3}^2\lesssim\|D^{m+1}a\|_0\|D^{m+1}\mathbf u,D^{m-1}\partial^2_2a\|^2_0+\Big|\int{\frac{1+\theta}{1+a}}\partial_2\mathbf u_1\cdot\partial^{m+1}_1a\partial^m_1\partial_2adx\Big|.$$
Based on the above  estimates of $I_{2,1}^2,I_{2,2}^2$ and $I_{2,3}^2$, we deduce that 
$$I_{2}^2\lesssim{(1+\|a,\theta\|_{m+1})}\|D^{m+1}a\|_0\|D^{m+1}\partial_2\mathbf u,D^{m+1}\mathbf u,D^{m-1}\partial^2_2a\|^2_0+\Big|\int{\frac{1+\theta}{1+a}}\partial_2\mathbf u_1\cdot\partial^{m+1}_1a\partial^m_1\partial_2adx\Big|.$$
It remains to estimate $I_{14}^i$. By using  integration by parts, we have
\begin{align*}
I_{14}^i=\frac{1}{2}\int\mathrm{div}\mathbf u|\partial^\alpha\partial_i\mathbf h|^2dx-\int[\partial^\alpha\partial_i,\mathbf u\cdot\nabla]\mathbf h\cdot\partial^\alpha\partial_i\mathbf hdx.
\end{align*}
 Thus, by a similar derivation in $I_2^i$,  we obtain 
\begin{align*}
I_{14}^1&\leq \|D\mathbf u\|_{L^\infty}\|D^{m+1}\mathbf h\|_0^2-\int[\partial^\alpha{\partial_1},\mathbf u\cdot\nabla]\mathbf h\cdot\partial^\alpha{\partial_1}\mathbf hdx\\
&\lesssim\|D^3\mathbf u,\partial_2^2D^2\mathbf h,D^m\mathrm{div}\mathbf u,D^m\partial_2\mathbf u\|_0\|D^{m+1}\mathbf u,D^{m+1}\mathbf h\|_0^2,
\end{align*}and
\begin{align*}
I_{14}^2&\leq \|D\mathbf u\|_{L^\infty}\|D^{m+1}\mathbf h\|_0\|D^{m-1}\partial_2^2\mathbf h\|_0-\int[\partial^\alpha{\partial_2},\mathbf u\cdot\nabla]\mathbf h\cdot\partial^\alpha{\partial_2}\mathbf hdx\\
&\lesssim\|D^{m+1}\mathbf h\|_0\|D^{m+1}\partial_2\mathbf u,D^{m+1}\mathbf u,D^{m-1}\partial^2_2\mathbf h\|^2_0+\Big|\int\partial_2\mathbf u_1\cdot\partial^{m+1}_1\mathbf h\cdot\partial^m_1\partial_2\mathbf hdx\Big|.
\end{align*}
\textbf{Estimate of $I_3^i$.}
When $i=1$,  it follows from \eqref{density-1}, \eqref{theta-1},  H\"{o}lder's inequality,  Sobolev's embedding  inequality and  Lemma \ref{Hm-estimate-1},  that
\begin{align*}
I_3^1\lesssim(\|D\mathbf u\|_{L^\infty}\|D^{m+1}a\|+\|Da\|_{L^\infty}\|D^m\mathrm{div}\mathbf u\|_{0})\|D^{m+1}a\|_0\lesssim\|D^3\mathbf u,D^m\mathrm{div}\mathbf u\|_0\|D^{m+1}a\|_0^2.
\end{align*}
When  $i=2$,  choosing $\alpha_k>0$,  we further divide $I_3^2$ into three terms,
\begin{align*}
I_3^2&=-R\int\frac{1+\theta}{1+a}\partial^\alpha(\partial_2a\mathrm{div}\mathbf u)\partial^\alpha\partial_2adx+R\int\partial^{\mathbf e_k}(\frac{1+\theta}{1+a})[\partial^\alpha,a]\partial_2\mathrm{div}\mathbf u\partial^{\alpha-\mathbf e_k}\partial_2adx\\
&\quad+R\int\frac{1+\theta}{1+a}([\partial^{\alpha+\mathbf e_k},a]\partial_2\mathrm{div}\mathbf u-\partial^{\mathbf e_k}a\partial^{\alpha}\partial_2\mathrm{div}\mathbf u)\partial^{\alpha-\mathbf e_k}\partial_2adx
\\&\triangleq I_{3,1}^2+I_{3,2}^2+I_{3,3}^2.
\end{align*}
Let us now estimate term by term in the above  identity. Employing  H\"{o}lder's inequality,  Sobolev's embedding  inequality, Lemma \ref{Hm-estimate-1} and \eqref{part-1}, we have
\begin{align*}
I_{3,1}^2&\lesssim(\|\partial_2a\|_{L^\infty}\|D^{{m}}\mathrm{div}\mathbf u\|_0+\|D^{m}\partial_2a\|_{0}\|\mathrm{div}\mathbf u\|_{L^\infty})\|D^m\partial_2a\|_0\\
&\lesssim\|D^{m+1}a\|_0\|D^{m}\mathrm{div}\mathbf u\|_0\|D^{m-1}\partial^2_2a\|_0,
\end{align*}
\begin{align*}
I_{3,2}^2&\lesssim\|D\theta,Da\|_{{L^\infty}}(\|Da\|_{L^\infty}\|D^{m}\mathrm{div}\mathbf u\|_0+\|D^{m}a\|_{0}\|\partial_2\mathrm{div}\mathbf u\|_{L^\infty})\|D^{m-1}\partial_2a\|_0\\
&\lesssim\|D^{m+1}\theta,D^{m+1}a\|_0^2\|D^{m+1}\mathrm{div}\mathbf u\|_0\|D^{m-1}\partial_2a\|_0,
\end{align*}and
\begin{align*}
I_{3,3}^2&\lesssim(\|D^{m+1}a\|_{0}\|\partial_2\mathrm{div}\mathbf u\|_{L^\infty}+\|Da\|_{L^\infty}\|D^{m+1}\mathrm{div}\mathbf u\|_{0})\|D^{m-1}\partial_2a\|_0\\&\lesssim\|D^{m+1}a\|_{0}\|D^{m+1}\mathrm{div}\mathbf u\|_{0}\|D^{m-1}\partial_2a\|_0.
\end{align*}
This along with  Young's inequality ensures
$$I_3^2\lesssim{(1+\|a,\theta\|_{m+1})}\|D^{m+1}a\|_0\|D^{m+1}\mathrm{div}\mathbf u,D^{m-1}\partial^2_2a\|^2_0.$$
\textbf{Estimates of $I_4^i, I_5^i, I_7^i$ and $I_{11}^i$.}
For any $i=1,2$, we deduce  by H\"{o}lder's and Young's inequalities, Sobolev's embedding inequality, \eqref{density-1}, \eqref{theta-1} and Lemma \ref{Hm-estimate-1},  that
\begin{align*}
I^i_{4}&\lesssim(\|Da\|_{L^\infty}\|D^{m}\theta\|_0+\|D^{m}a\|_0\|D\theta\|_{L^\infty})\|D^{m+1}\partial_i\mathbf u\|_0\lesssim\|D^{m+1}a,D^{m+1}\theta\|_0\|D^{m+1}\partial_i\mathbf u,D^{m+1}\theta\|^2_0,
\end{align*}
\begin{align*}
I^i_{5}&\lesssim(\|Da\|_{L^\infty}\|D^{m}\mathrm{div}\mathbf u\|_0+\|D^{m+1}a\|_0\|\mathrm{div}\mathbf u\|_{L^\infty})\|D^{m+1}\theta\|_0\lesssim\|D^{m+1}a,D^{m+1}\theta\|_0\|D^{m}\mathrm{div}\mathbf u,D^{m+1}\theta\|^2_0,
\end{align*}
\begin{align*}
I^i_{7}&\lesssim((1+\|a,\mathbf u\|_{L^\infty})\|Da,D\mathbf u\|_{L^\infty}\|D^{m+1}\mathbf u\|_0+(1+\|a,\mathbf u\|_{L^\infty})\|D^{m+1}a,D^{m+1}\mathbf u\|_0\|D\mathbf u\|_{L^\infty})\|D^{m+1}\mathbf u\|_0\\
&\lesssim(1+\|a,\mathbf u\|_{m+1})\|Da,D\mathbf u\|_m\|D^{m+1}\mathbf u\|^2_0,
\end{align*}and
\begin{align*}
I^i_{11}&\lesssim((1+\|a,\mathbf u\|_{L^\infty})\|Da,D\mathbf u\|_{L^\infty}\|D^{m+1}\theta\|_0+(1+\|a,\mathbf u\|_{L^\infty})\|D^{m+1}a,D^{m+1}\mathbf u\|_0\|D\theta\|_{L^\infty})\|D^{m+1}\theta\|_0\\
&\lesssim(1+\|a,\mathbf u\|_{m+1})\|Da,D\mathbf u\|_m\|D^{m+2}\theta\|^2_0.
\end{align*}
\textbf{Estimate of $I_6^i$.}
When $i=1$,  by \eqref{heat-mhd}$_2$,  we can divide $I^1_{6}$ into four terms,
\begin{align*}
I^1_{6}&=\int[\partial^\alpha\partial_1,a]\big((\mathbf u\cdot\nabla)\mathbf u+R\nabla\theta-\frac{\mu\Delta\mathbf u+(\mu+\lambda)\nabla\mathrm{div}\mathbf u}{1+a}\big)\cdot \partial^\alpha\partial_1\mathbf udx\\
&\quad +\int[\partial^\alpha\partial_1,a]\big(\frac{{R}\partial_1a(1+\theta)}{1+a}+\frac{\partial_1\mathbf h_2}{1+a}+\frac{\partial_1(\mathbf h_2)^2}{2(1+a)}-\frac{\partial_2\mathbf h_1}{1+a}-\frac{\mathbf h_2\partial_2\mathbf h_1}{1+a}\big)\partial^\alpha\partial_1\mathbf u_1dx\\
&\quad+\int[\partial^\alpha\partial_1,a]\big(\frac{\partial_2(\mathbf h_1)^2}{2(1+a)}-\frac{\mathbf h_1\partial_1\mathbf h_2}{1+a}\big)\partial^\alpha\partial_1\mathbf u_2dx+R\int [\partial^\alpha\partial_1,a]\big(\frac{\partial_2a(1+\theta)}{1+a}\big)\partial^\alpha\partial_1\mathbf u_2dx
\\&\triangleq \sum_{j=1}^4I^1_{6,j}.
\end{align*} %Let us now bound term by term in the above equality as follows.
We bound term by term above in what follows. For $I^1_{6,1}$,
thanks to  H\"{o}lder's and Young's inequalities, Sobolev's embedding inequality, \eqref{density-1}, Lemmas \ref{Hm-estimate-1} and  \ref{nolinear-estimate},  we get
\begin{align*}
I^1_{6,1}&\lesssim\Big(\|Da\|_{L^\infty}\big(\|\mathbf u\|_{L^\infty}\|D^{m+1}\mathbf u\|_0+\|D^m\mathbf u\|_{0}\|D\mathbf u\|_{L^\infty}+\|D^{m+1}\theta\|_0+\|D^{m+2}\mathbf u\|_0\\
&\quad+\|D^2\mathbf u\|_{L^\infty}\|D^{m}\big((1+a)^{-1}-1\big)\|_0\big)+\|D^{m+1}a\|_0(\|\mathbf u\|^2_{L^\infty}\|D\mathbf u\|_{L^\infty}+\|\nabla\theta,D^2\mathbf u\|_{L^\infty})\Big)\|D^{m+1}\mathbf u\|_0\\
&\lesssim\|D^{m+1}a\|_0(1+\|a,u\|_m)\|D^{m+2}\mathbf u,D^{m+2}\theta\|_0^2.
\end{align*}
Similarly, for $I^1_{6,2}$, by virtue of  \eqref{theta-1},  one has
\begin{align*}
I^1_{6,2}&\lesssim\|D^{m+1}a\|_0\big((1+\|a,\theta\|_m)^2\|D^{m+1}a\|_0+(1+\|a\|_m)(1+\|\mathbf h\|_m)\|D^{m+1}\mathbf h\|_0\big)\|D^{m+1}\mathbf u_1\|_0\\
&\lesssim(1+\|a,\theta,\mathbf h\|_{m+1})^2\|D^{m}\partial_2\mathbf u,D^{m}\mathrm{div}\mathbf u\|_0\|D^{m+1}a,D^{m+1}{\mathbf h}\|_0^2.
\end{align*}
For $I^1_{6,3}$, it follows from
 H\"{o}lder's and Young's inequalities, Sobolev's embedding inequality,  Poincar\'{e} inequality,  \eqref{h-junzhi}, \eqref{density-1},  Lemmas \ref{Hm-estimate-1} and  \ref{nolinear-estimate}, that
\begin{align*}
I^1_{6,3}&\lesssim(1+\|a\|_m)\|D^{m+1}a\|_0\|D^{m+1}\mathbf h\|_0\|\mathbf h_1\|_m\|D^{m+1}\mathbf u_2\|_0\\
&\lesssim(1+\|a\|_m)\|D^{m+1}a\|_0\|D^{m+1}\mathbf h\|_0\|D^{m-1}\partial^2_2\mathbf h,D^{m+2}\mathbf u\|_0^2.
\end{align*}
To bound  $I^1_{6,4}$, by using  integration by parts, we further divide it into two terms,
\begin{align*}
I^1_{6,4}&=-R\int [\partial^\alpha,a]\big(\frac{\partial_2a(\theta+1)}{1+a}\big)\partial^\alpha\partial^2_1\mathbf u_2dx+R\int\partial_1a\partial^\alpha \big(\frac{\partial_2a(1+\theta)}{1+a}\big)\partial^\alpha\partial_1\mathbf u_2dx
\\&\triangleq I^1_{6,4,1}+I^1_{6,4,2}.
\end{align*}
For $I^1_{6,4,1}$, using H\"{o}lder's and Young's inequalities, Sobolev's embedding inequality, \eqref{density-1}, \eqref{theta-1}, Lemmas \ref{Hm-estimate-1} and  \ref{nolinear-estimate}, we infer that
\begin{align*}I^1_{6,4,1}&\lesssim\Big(\|Da\|_{L^\infty}\big(\|D^{m-1}\partial_2a\|_0+\|\partial_2a\|_{L^4}(\|D^m\theta\|_{0}+\|a\|_m)\big)\\
&\quad+\|D^{m-1}a\|_{L^4}\|\partial_2a\|_{L^4}\Big)\|D^{m+2}\mathbf u_2\|_0\\&\lesssim\|D^{m+1}a\|_0(1+\|a,\theta\|_m)^2\|\partial^2_2D^{m-1}a,D^{m+2}\mathbf u\|^2_0.
\end{align*}
For $I^1_{6,4,2}$,  choosing $\alpha_k>0$,   then using H\"{o}lder's and Young's inequalities, Sobolev's embedding inequality, Poincar\'{e}'s inequality,   Lemmas \ref{Hm-estimate-1} and  \ref{nolinear-estimate} and  integration by parts,  we deduce that
\begin{align*}I^1_{6,4,2}&=-\int\partial^{\alpha-\mathbf e_k}\big(\frac{\partial_2a(1+\theta)}{1+a}\big)\partial^{\mathbf e_k}(\partial_1a\partial^\alpha\partial_1\mathbf u_2)dx\\
&\lesssim\|D^{m+1}a\|_0(1+\|a,\theta\|_m)^2\|\partial^2_2D^{m-1}a,D^{m+2}\mathbf u\|^2_0.
\end{align*}
Hence,
\begin{align*}I^1_{6,4}\lesssim\|D^{m+1}a\|_0(1+\|a,\theta\|_m)^2\|\partial^2_2D^{m-1}a,D^{m+2}\mathbf u\|^2_0.
\end{align*}
Summing up the above estimates of $I^1_{6,1},I^1_{6,2},I^1_{6,3}$ and $I^1_{6,4}$, we have
\begin{align*}I^1_{6}&\lesssim\|D^{m+1}a\|_0(1+\|a,\mathbf u,{\theta,}\mathbf h\|_{m+1})^2\|{D^{m-1}\partial_2^2a,}D^{m+2}\mathbf u,D^{m+2}\theta,D^{m-1}\partial_2^2\mathbf h\|_0^2\\
&\quad +(1+\|a,\theta,\mathbf h\|_{m+1})^2\|D^{m}\partial_2\mathbf u,D^{m}\mathrm{div}\mathbf u\|_0\|D^{m+1}a,D^{m+1}{\mathbf h}\|_0^2.
\end{align*}
When $i=2$, for $\alpha_k>{0}$,  we can  divide $I^2_{6}$ into three terms,
\begin{align*}
I^2_{6}=&\int\partial^{\alpha-\mathbf e_k}(\partial_2a\partial_t\mathbf u)\cdot \partial^{\alpha+\mathbf e_k}\partial_2\mathbf udx+\int[\partial^{\alpha-\mathbf e_k},a]\partial_2\partial_t\mathbf u\cdot \partial^{\alpha+\mathbf e_k}\partial_2\mathbf udx\\
&\quad-\int\partial^{\mathbf e_k}a\partial^{\alpha-\mathbf e_k}\partial_2\partial_t\mathbf u\cdot \partial^{\alpha}\partial_2\mathbf udx
\\&\triangleq I^2_{6,1}+I^2_{6,2}+I^2_{6,3}.
\end{align*}In what follows,  we estimate each term in the above identity separately.  For  $I^2_{6,1}$,
 thanks to $\eqref{heat-mhd}_2$, H\"{o}lder's and Young's inequalities, Sobolev's embedding inequality, Poincar\'{e}'s inequality, Lemmas \ref{Hm-estimate-1} and  \ref{nolinear-estimate}, we conclude that
\begin{align*}
I^2_{6,1}&\lesssim\Big(\|\partial_2a\|_{L^4}\big(\|D^{m-1}((\mathbf u \cdot\nabla)\mathbf u)\|_{L^4}+\|D^{m} \theta\|_{L^4}+\|D^{m-1}\big(\frac{\nabla a(1+\theta)}{1+a}\big)\|_{L^4}+\|D^{m-1}\big(\frac{\nabla\mathbf h_2+\frac12{\nabla}|\mathbf h|^2}{1+a}\big)\|_{L^4}\\
&\quad +\|D^{m-1}\big(\frac{\partial_2\mathbf h+(\mathbf h\cdot\nabla)\mathbf h}{1+a}\big)\|_{L^4}\big)+\|\partial_2D^{m-1}a\|_{0}\big(\|\nabla(\theta,a,\mathbf h)\|_{L^\infty}+\|\mathbf u,\mathbf h\|_{L^\infty}\|D(\mathbf u,\mathbf h)\|_{L^\infty}\big)\\
&\quad +\|\partial_2a\|_{L^\infty}\|D^{m-1}\big(\frac{\mu\Delta\mathbf u+(\mu+\lambda)\nabla\mathrm{div}\mathbf u}{1+a}\big)\|_{0}+\|D^{m-1}\partial_2a\|_{L^4}\|D^2\mathbf u\|_{L^4}\Big)\|D^{m+1}\partial_2\mathbf u\|_0\\
&\lesssim\Big(\|D^{m-1}\partial_2a\|_0\|Da,D\mathbf u,D\theta,D\mathbf h\|_m(1+\|a,\mathbf u,\theta,\mathbf h\|_{m+1})^2\\
&\quad+\|Da\|_{m}(1+\|a\|_{m+1})\|D^{m+1}\mathbf u\|_0\Big)\|D^{m+1}\partial_2\mathbf u\|_{{0}}\\
&\lesssim\|Da,D\mathbf u,D\theta,D\mathbf h\|_m(1+\|a,\mathbf u,\theta,\mathbf h\|_{m+1})^2\|D^{m-1}\partial^2_2a,D^{m+1}\partial_2\mathbf u,D^{m+1}\mathbf u\|^2_0.
\end{align*}
In order to  bound $I^2_{6,2}$, noting that
\begin{align*}
\partial_2\partial_t\mathbf u&=-(\partial_2\mathbf u\cdot\nabla)\mathbf u-(\mathbf u\cdot\nabla)\partial_2\mathbf u-R\partial_2\nabla\theta-{R}\frac{\nabla \partial_2a(1+\theta)+\nabla a\partial_2\theta}{1+a}+R\frac{\nabla a\partial_2a(1+\theta)}{(1+a)^2}\\
&\quad +\frac{\mu\Delta\partial_2\mathbf u+(\mu+\lambda)\nabla\mathrm{div}\partial_2\mathbf u}{1+a}-\frac{\partial_2a\big(\mu\Delta\mathbf u+(\mu+\lambda)\nabla\mathrm{div}\mathbf u\big)}{(1+a)^2}-\frac{2\nabla\partial_2\mathbf h_2+\nabla\partial_2(|\mathbf h|^2)}{2(1+a)}\\
&\quad +\frac{(2\nabla\mathbf h_2+\nabla|\mathbf h|^2)\partial_2a}{2(1+a)^2}+\frac{\partial_2^2\mathbf h+(\partial_2\mathbf h\cdot\nabla)\mathbf h+(\mathbf h\cdot\nabla)\partial_2\mathbf h}{1+a}-\frac{\partial_2a(\partial_2\mathbf h+(\mathbf h\cdot\nabla)\mathbf h)}{(1+a)^2},
\end{align*}
we conclude by \eqref{product}, Lemma \ref{nolinear-estimate}, Sobolev's embedding inequality, Poincar\'{e}'s inequality, \eqref{density-1} and \eqref{theta-1}, that
\begin{align}\label{10-12-1}
\|D^{m-2}\partial_2\partial_t\mathbf u\|_0\lesssim(1+\|a,\mathbf u,\theta,\mathbf h\|_{m+1})^3\|\partial_2D^{m+1}\mathbf u,D^{m+1}\mathbf u,D^{m+2}\theta,D^{m-1}\partial_2^2a,D^{m-1}\partial_2^2\mathbf h\|_0.\end{align}
%In particular, when $m=2$, the above inequality follows  from  H\"{o}lder's inequality, Sobolev's embedding  inequality and Poincar\'{e}'s inequality, \eqref{density-1} and \eqref{theta-1}.
For $m=2$,  it  then follows from  H\"{o}lder's and Young's inequalities, Sobolev's embedding inequality,  Poincar\'{e}'s inequality and \eqref{10-12-1}, that
\begin{align*}I^2_{6,2}&\lesssim\|Da\|_{L^\infty}\|D^{m-2}\partial_2\partial_t\mathbf u\|_0\|D^{m+1}\partial_2\mathbf u\|_0\\
&\lesssim(1+\|a,\mathbf u,\theta,\mathbf h\|_{m+1})^3\|Da\|_{m}\|\partial_2D^{m+1}\mathbf u,D^{m+1}\mathbf u,D^{m+2}\theta,D^{m-1}\partial_2^2a,D^{m-1}\partial_2^2\mathbf h\|^2_0.
\end{align*}
For $m>2$,  employing  Lemma \ref{Hm-estimate-1} and \eqref{10-12-1}  leads to
\begin{align*}I^2_{6,2}&\lesssim(\|Da\|_{L^\infty}\|D^{m-2}\partial_2\partial_t\mathbf u\|_0+\|D^{m-1}a\|_{L^4}\|\partial_2\partial_t\mathbf u\|_{L^4})\|D^{m+1}\partial_2\mathbf u\|_0\\
&\lesssim(1+\|a,\mathbf u,\theta,\mathbf h\|_{m+1})^3\|Da\|_{m}\|D^{m+1}\partial_2\mathbf u,D^{m+1}\mathbf u,D^{m+2}\theta,D^{m-1}\partial_2^2a,D^{m-1}\partial_2^2\mathbf h\|^2_0.
\end{align*}
Next, we deal with  {$I^2_{6,3}$}. Due to $|\alpha-\mathbf e_k|>0$,  without loss of generality,  we may assume  that its $l$-th component is nonzero. Thus,  we deduce  by integration by parts, H\"{o}lder's and Young's inequalities, Sobolev's embedding inequality,  Poincar\'{e}'s inequality  and \eqref{10-12-1}, that
\begin{align*}
I^2_{6,3}&=\int\partial^{\alpha-\mathbf e_k-\mathbf e_l}\partial_2\partial_t\mathbf u\cdot \partial^{\mathbf e_l}(\partial^{\mathbf e_k}a\partial^{\alpha}\partial_2\mathbf u)dx\\
&\lesssim\|D^{m-2}\partial_2\partial_t\mathbf u\|_0(\|Da\|_{L^\infty}\|D^{m+1}\partial_2\mathbf u\|_0+\|D^2a\|_{L^4}\|\|D^{m}\partial_2\mathbf u\|_{L^4})\\
&\lesssim(1+\|a,\mathbf u,\theta,\mathbf h\|_{m+1})^3\|Da\|_{m}\|D^{m+1}\partial_2\mathbf u,D^{m+1}\mathbf u,D^{m+2}\theta,D^{m-1}\partial_2^2a,D^{m-1}\partial_2^2\mathbf h\|^2_0.
\end{align*}
Gathering the above estimates of $I^2_{6,1},I^2_{6,2}$ and $I^2_{6,3}$, we have
\begin{align*}I^2_{6}&\lesssim(1+\|a,\mathbf u,\theta,\mathbf h\|_{m+1})^3\|D(a,\mathbf u,\theta,\mathbf h)\|_{m}\|D^{m+1}\partial_2\mathbf u,D^{m+1}\mathbf u,D^{m+2}\theta,D^{m-1}\partial_2^2a,D^{m-1}\partial_2^2\mathbf h\|^2_0.
\end{align*}
 \textbf{Estimate of $I_8^i$.}
When $i=1$,  it then follows from  integration by parts, H\"{o}lder's and Young's inequalities, Sobolev's embedding inequality and Lemma \ref{Hm-estimate-1}, that
\begin{align*}
I_8^1&=-\sum_{j=1}^2\int[\partial^\alpha\partial_1,\mathbf h_j]\partial_1\mathbf h_j\partial^\alpha\mathrm{div}\mathbf udx+\sum_{j=1}^2\int\partial_1\mathbf h_j\partial^\alpha\partial_1\mathbf h_j\partial^\alpha\mathrm{div}\mathbf udx\\
&\lesssim\|D^{m}\mathrm{div}\mathbf u\|_0\|D^{m+1}\mathbf h\|_0^2.
\end{align*}
When $i=2$,  thanks to  H\"{o}lder's and Young's inequalities, Sobolev's embedding inequality and Lemma \ref{Hm-estimate-1}, we infer that 
\begin{align*}
I_8^2&\lesssim(\|D\mathbf h\|_{L^\infty}\|D^{m-1}\partial_2\mathbf h\|_0+\|D^m\mathbf h\|_{L^4}\|\partial_2\mathbf h\|_{L^4})\|D^{m+1}\partial_2\mathbf u\|_0\\
&\lesssim\|D\mathbf h\|_{m}\|D^{m-1}\partial_2^{{2}}\mathbf h,D^{m+1}\partial_2\mathbf u\|_0^2.
\end{align*}
 \textbf{Estimate of $I_9^i$.}
When $i=2$,  using H\"{o}lder's and Young's inequalities, \eqref{commutator-2}, $\mathrm{div}\mathbf h=0$ and \eqref{part-1}, we have
\begin{align*}
I_9^2&\lesssim\|[\partial^\alpha\partial_2,\mathbf h\cdot\nabla]\mathbf h\|_0\|\partial^\alpha\partial_2\mathbf u\|_0\\
&\lesssim(\|D\mathbf h\|_{m}\|D^{m-1}\partial_2^2\mathbf h\|+\|D^m\partial_2\mathbf h\|_0^2)\|D^{m}\partial_2\mathbf u\|_0\\
&\lesssim\|D\mathbf h\|_{m}\|D^{m-1}\partial_2^2\mathbf h,D^{m+1}\partial_2\mathbf u\|_0^2.
\end{align*}
When $i=1$, if $\alpha_2>0,$ then the estimate coincides with that in the case $i=2$ above.
If $\alpha_2=0$, then  using integration by parts, \eqref{commutator-1},  H\"{o}lder's and Young's inequalities, Sobolev's embedding inequality and  Poincar\'{e}'s inequality, we get 
\begin{align*}
I_9^1&=-\int[\partial_1^{m},\mathbf h\cdot\nabla]\mathbf h\cdot\partial^{m+2}_1\mathbf udx+\int(\partial_1\mathbf h\cdot\nabla)\partial_1^m\mathbf h\cdot\partial^{m+1}_1\mathbf udx\\
&=-\int[\partial_1^{m},\mathbf h\cdot\nabla]\mathbf h\cdot\partial^{m+2}_1\mathbf udx-\int(\partial_1\mathbf h\cdot\nabla)\partial_1^{m-1}\mathbf h\cdot\partial^{m+2}_1\mathbf udx-\int(\partial^2_1\mathbf h\cdot\nabla)\partial_1^{m-1}\mathbf h\cdot\partial^{m+1}_1\mathbf udx\\
&\lesssim\|D^m\mathbf h\|_0\|D^2\partial^2_2\mathbf h,D^{m-1}\partial^2_2\mathbf h\|_0\|D^{m+2}\mathbf u\|_0+\|\partial_2\mathbf h_2\|_{L^4}\|D^m\mathbf h\|_{L^4}\|D^{m+2}\mathbf u\|_0\\
&\quad+\|\partial_1\mathbf h_2\|_{L^\infty}\|D^{m-1}\partial_2\mathbf h\|_0\|D^{m+2}\mathbf u\|_0+\|\partial_1\partial_2{\mathbf h_2}\|_{0}\|D^m\mathbf h\|_{L^4}\|D^{m+1}\mathbf u\|_{L^4}\\
&\quad+\|D^2{\mathbf h_2}\|_{L^4}\|D^{m-1}\partial_2\mathbf h\|_{0}\|D^{m+1}\mathbf u\|_{L^4}\\
&\lesssim\|D\mathbf h\|_m\|D^{m+2}\mathbf u,D^{m-1}\partial_2^2\mathbf h,D^2\partial_2^2\mathbf h\|_0^2.
\end{align*}
 \textbf{Estimate of $I_{10}^i$.}
It follows from $\eqref{heat-mhd}_{{3}}$, H\"{o}lder's and Young's inequalities, Sobolev's embedding inequality,  Poincar\'{e}'s inequality,
\eqref{density-1} and \eqref{product} that
\begin{align*}
\|D^{m}\partial_t\theta\|_{0}&\lesssim\|D^m\big((\mathbf u\cdot\nabla)\theta\big)\|_0+\|D^m(\theta\mathrm{div}\mathbf u)\|_0+\|D^{m}\big(\frac{\Delta\theta}{1+a}\big)\|_0+\|D^m\big(\frac{\mathcal Q(\mathbf u)}{1+a}\big)\|_0\\
&\lesssim(1+\|a,\mathbf u,\theta\|_{1+m})^2\|D^{m+1}\mathbf u,D^{m+2}\theta\|_0,
\end{align*}and
\begin{align*}
\|\partial_t\theta\|_{L^\infty}&\lesssim\|(\mathbf u\cdot\nabla)\theta\|_{L^\infty}+\|\theta\mathrm{div}\mathbf u\|_{L^\infty}+\|\frac{\Delta\theta}{1+a}\|_{L^\infty}+\|\frac{\mathcal Q(\mathbf u)}{1+a}\|_{L^\infty}\\
&\lesssim(1+\|a,\mathbf u,\theta\|_{1+m})\|D^{m+1}\mathbf u,D^{m+2}\theta\|_0.
\end{align*}
From which, we have
\begin{align*}I_{10}^i&\lesssim(\|Da\|_{L^\infty}\|D^m\partial_t\theta\|_{0}+\|D^{{m+1}}a\|_0\|\partial_t\theta\|_{L^\infty})\|D^{m+1}\theta\|_0\\
&\lesssim\|Da\|_m(1+\|a,\mathbf u,\theta\|_{1+m})\|D^{m+1}\mathbf u,D^{m+2}\theta\|^2_0.\end{align*}
 \textbf{Estimates of $I_{12}^i$ and $I_{13}^i$.} For any $i=1,2$, we deduce  by integration by parts, H\"{o}lder's and Young's inequalities, Sobolev's embedding inequality, \eqref{density-1}, \eqref{theta-1} and \eqref{product}, that
\begin{align*}
I_{12}^i&=R\int\partial^\alpha((1+a)\theta\mathrm{div}\mathbf u)\partial^\alpha\partial^2_i\theta dx\\
&\lesssim(\|\theta\|_{L^\infty}\|D^{m+1}\mathbf u\|_{0}+\|D^ma,D^m\theta\|_0\|D\mathbf u\|_{L^\infty})\|D^{m+2}\theta\|_0\\
&\lesssim(\|\theta\|_{2}+\|D(a,\theta)\|_{m})\|D^{m+1}\mathbf u,D^{m+2}\theta\|_0^2,
\end{align*}and
\begin{align*}I_{13}^i&=-\int\partial^\alpha\mathcal Q(\mathbf u)\partial^\alpha\partial^2_i\theta dx\lesssim\|D\mathbf u\|_{m}\|D^{m+1}\mathbf u,D^{m+2}\theta\|_0^2.
\end{align*}
 \textbf{Estimates of $I_{15}^i$ and $I_{16}^i$.}
For  $I_{15}^1$, it follows from H\"{o}lder's and Young's inequalities, \eqref{commutator-1} and \eqref{commutator-2} that
\begin{align*}
I_{15}^1&\lesssim(\|D^2\partial_2\mathbf h\|_0\|D^{m+1}\mathbf u\|_0+\|D^m\partial_2\mathbf h\|_0\|D^{3}\mathbf u\|_0+\|D^3\mathbf h\|_0\|D^{m}\partial_2\mathbf u\|_0\\
&\quad+\|D^{m+1}\mathbf h\|_0\|D^{2}\partial_2\mathbf u\|_0)\|D^{m+1}\mathbf h\|_0\\
&\lesssim\|D^2\partial_2\mathbf h,D^3\mathbf u,D^m\partial_2\mathbf u\|_0\|D^{m+1}\mathbf u,D^{m+1}\mathbf h\|_0^2.
\end{align*}
For $I_{15}^2$, choosing $\alpha_k>0$ and by integration by parts, we have
\begin{align*}
I_{15}^2&=-\int[\partial^\alpha,\mathbf h\cdot\nabla]\mathbf u\cdot\partial^\alpha\partial^2_2\mathbf hdx+\int(\partial_2\mathbf h\cdot\nabla)\partial^\alpha\mathbf u\cdot\partial^\alpha\partial_2\mathbf hdx\\
&=\int[\partial^{\alpha+\mathbf e_k},\mathbf h\cdot\nabla]\mathbf u\cdot\partial^{\alpha-\mathbf e_k}\partial^2_2\mathbf hdx-\int(\partial^{\mathbf e_k}\mathbf h\cdot\nabla)\partial^{\alpha}\mathbf u\cdot\partial^{\alpha-\mathbf e_k}\partial^2_2\mathbf hdx
\\&\quad+\int(\partial_2\mathbf h\cdot\nabla)\partial^\alpha\mathbf u\cdot\partial^\alpha\partial_2\mathbf hdx\\
&\triangleq I_{15,1}^2+I_{15,2}^2+I_{15,3}^2.
\end{align*} We bound term by term above in what follows.  For $I_{15,1}^2$,
by virtue of H\"{o}lder's and Young's inequalities,   Poincar\'{e}'s inequality
 and  \eqref{commutator}, we infer that
$$I_{15,1}^2\lesssim\|[\partial^{\alpha+\mathbf e_k},\mathbf h\cdot\nabla]\mathbf u\|_0\|D^{m-1}\partial_2^2\mathbf h\|\lesssim\|D\mathbf h\|_m\|D^{m+1}\mathbf u,D^{m-1}\partial_2^2\mathbf h\|_0^2.$$
 For $I_{15,2}^2$, thanks to H\"{o}lder's and Young's inequalities, Sobolev's embedding inequality and Poincar\'{e}'s inequality, we have
\begin{align*}
&I_{15,2}^2\lesssim\|D\mathbf h\|_{L^\infty}\|D^{m+1}\mathbf u\|_0\|D^{m-1}\partial_2^2\mathbf h\|_0\lesssim\|D\mathbf h\|_m\|D^{m+1}\mathbf u,D^{m-1}\partial_2^2\mathbf h\|_0^2,
\end{align*}
In the same
manner, we get
\begin{align*}
&I_{15,3}^2\lesssim\|D^{m+1}\mathbf u\|_0\|D^{m}\partial_2\mathbf h\|^2_0\lesssim\|D\mathbf h\|_m\|D^{m+1}\mathbf u,D^{m-1}\partial_2^2\mathbf h\|_0^2.
\end{align*}
Collecting these estimates of $I_{15,1}^2,I_{15,2}^2$ and $I_{15,3}^2$, we arrive at
$$I_{15}^2\lesssim\|D\mathbf h\|_m\|D^{m+1}\mathbf u,D^{m-1}\partial_2^2\mathbf h\|_0^2.$$
It remains to estimate $I_{16}^i$. Employing H\"{o}lder's and Young's inequalities, Sobolev's embedding inequality,  Poincar\'{e}'s inequality
 and \eqref{product}, we obtain
$$I_{16}^1\lesssim(\|D{\mathbf h}\|_{L^\infty}\|D^{m}\mathrm{div}\mathbf u\|_0+\|D^{m+1}\mathbf h\|_0\|\mathrm{div}\mathbf u\|_{L^\infty})\|D^{m+1}\mathbf h\|_0\lesssim\|D^m\mathrm{div}\mathbf u\|_0\|D^{m+1}\mathbf h\|_0^2.$$
 By a similar derivation in $I_{15}^2$, we get
\begin{align*}
I_{16}^2&=-\sum_{j=1}^2\int[\partial^{\alpha+\mathbf e_k},\mathbf h_j]\mathrm{div}\mathbf u\partial^{\alpha-\mathbf e_k}\partial^2_2\mathbf h_jdx+\sum_{j=1}^2\int\partial^{\mathbf e_k}\mathbf h_j\partial^\alpha\mathrm{div}\mathbf u\partial^{\alpha-\mathbf e_k}\partial^2_2\mathbf h_jdx\\
&\quad-\sum_{j=1}^2\int\partial_2\mathbf h_j\partial^\alpha\mathrm{div}\mathbf u\partial^{\alpha}\partial_2\mathbf h_jdx\\
&\lesssim\|D\mathbf h\|_m\|D^{m+1}\mathbf u,D^{m-1}\partial_2^2\mathbf h\|_0^2.
\end{align*}
Finally, substituting the above estimates of $I_j^i(j=1,2,\cdots,16)$  into \eqref{zuigao-estiamte}, summing over all multi-indices
 $\alpha$ with $|\alpha|=m$, and taking $i=1$ and $i=2$, respectively, we obtain the corresponding conclusions \eqref{1-h-older} and \eqref{2-h-older}. This completes the proof of Lemma \ref{H-priori-v}.
\end{proof}
\subsection{The dissipation for $\partial_2^2a$ and $\partial_2^2\mathbf h$} In this subsection, we exploit the coupling between the velocity and magnetic field equations to recover the dissipation missing for \(\partial_2^2a\) and \(\partial_2^2\mathbf h\).
\begin{lemma}\label{dissipation-rho} Let $(a,\mathbf u,\theta,\mathbf h)$ be smooth solutions to  the system $\eqref{heat-mhd}$
  on $\mathbb T^2\times[0,T)$ for some $T>0$, and suppose that \eqref{density-1}, \eqref{theta-1} and \eqref{junzhi} hold.  Then, for any $m\geq2$,
 \begin{align}
{\frac2R}\frac{d}{dt}&\sum_{|\alpha|=m-1}\int\partial^\alpha\partial_2\mathbf u_2\partial^\alpha\partial^2_2adx+\|D^{m-1}\partial^2_2a\|_0^2\leq C_0\|D^{m+1}\partial_2\mathbf u,D^{m+1}\partial_2\theta\|^2_0\nonumber\\
&\quad+(1+\|a,\theta\|_{m+1})\|a,\mathbf u,\theta,\mathbf h\|_{m+1}\|D^{m-1}\partial_2^2a,D^{m+1}\partial_2\mathbf u,D^{m+1}\mathbf u,D^{m-1}\partial^2_2\mathbf h\|_0^2,\label{10-16-partial2a}
\end{align}
where  $C_0$ is  a positive constant depending only on $|\mathbb T|,R,\mu,\lambda,m$.
\end{lemma}
\begin{proof}
It follows from \eqref{heat-mhd}$_2$ that
\begin{equation*}
R\partial_2a=\frac{R(a-\theta)\partial_2a}{1+a}-\partial_t{\mathbf u_2}-(\mathbf u\cdot\nabla)\mathbf u_2-R\partial_2\theta-\frac{\mathbf h_1\partial_2\mathbf h_1}{1+a}+\frac{\mathbf h_1\partial_1\mathbf h_2}{1+a}+\frac{\mu\Delta\mathbf u_2+(\mu+\lambda)\partial_2\mathrm{div}\mathbf u}{1+a}.
\end{equation*}
Given any multi-indices $\alpha$ with $|\alpha|=m-1$, applying $\partial^{\alpha}\partial_2$ to the above identity, and  then performing $L^2$ inner product
with $\partial^{\alpha}\partial_2^2a$ yield
\begin{equation}R\|\partial^{\alpha}\partial^2_2a\|_0^2=\sum_{j=1}^{6}II_j\label{rho-10-13-1}
\end{equation}
with
\begin{align*}
&II_1:=R\int\partial^{\alpha}\partial_2\big(\frac{(a-\theta)\partial_2a}{1+a}\big)\partial^{\alpha}\partial^2_2adx,\quad II_2:=-\int\partial^{\alpha}\partial_2\partial_t{\mathbf u_2}\partial^{\alpha}\partial^2_2adx,\\
&II_3:=-\int\partial^{\alpha}\partial_2((\mathbf u\cdot\nabla){\mathbf u_2})\partial^{\alpha}{\partial_2^2}adx,\quad II_4:=-R\int\partial^\alpha\partial^2_2\theta\partial^{\alpha}\partial^2_2adx,\\
&II_5:=-\int\partial^\alpha\partial_2\big(\frac{\mathbf h_1\partial_2\mathbf h_1}{1+a}-\frac{\mathbf h_1\partial_1\mathbf h_2}{1+a}\big)\partial^\alpha\partial^2_2adx,\quad II_6:=\int\partial^\alpha\partial_2\big(\frac{\mu\Delta{\mathbf u_2}+(\mu+\lambda){\partial_2}\mathrm{div}\mathbf u}{1+a}\big)\partial^\alpha\partial_2^2adx.
\end{align*}
We bound term by term above in what follows.  For $II_1$, it follows from \eqref{product},  Sobolev's embedding inequality, Poincar\'{e}'s inequality,  \eqref{density-1}, \eqref{theta-1} and Lemma \ref{nolinear-estimate}, that
\begin{align*}
&\|\partial^{\alpha}\partial_2\big(\frac{(a-\theta)\partial_2a}{1+a}\big)\|_{0}
\\&\lesssim\|D^{m-1}\big(\frac{(a-\theta)\partial^2_2a}{1+a}\big)\|_0+\|D^{m-1}\big(\frac{(\partial_2a-\partial_2\theta)\partial_2a}{1+a}\big)\|_0+\|D^{m-1}\big(\frac{(a-\theta)(\partial_2a)^2}{(1+a)^2}\big)\|_0\\
&\lesssim\|a,\theta\|_{L^\infty}\|D^{m-1}\partial_2^2a\|_0+\|\partial^2_2a\|_{L^4}\|a,\theta\|_{m}+\|D(a,\theta)\|_{L^\infty}\|D^{m-1}\partial_2a\|+\|\partial_2a\|_{L^4}(\|D(a,\theta)\|_{m}\\
&\quad +\|D(a,\theta)\|_{L^\infty}\|a\|_m)+\|D^{m-1}\partial_2a\|_0\|\partial_2a\|_{L^\infty}+\|\partial_2a\|_{L^4}(\|Da\|_{m}+\|Da\|_{L^\infty}\|a,\theta\|_{m+1})\\
&\lesssim (1+\|a,\theta\|_{m+1})\|a,\theta\|_{m+1}\|D^{m-1}\partial_2^2a\|_0,
\end{align*}
which together with  H\"{o}lder's and Young's inequalities,  yields
\begin{align}
II_1\lesssim\|\partial^{\alpha}\partial_2\big(\frac{(a-\theta)\partial_2a}{1+a}\big)\|_{0}\|\partial^{\alpha}\partial^2_2a\|_0\lesssim (1+\|a,\theta\|_{m+1})\|a,\theta\|_{m+1}\|D^{m-1}\partial_2^2a\|^2_0.\label{rho-10-13-2}
\end{align}
For $II_2$, by using  integration by parts and $\eqref{heat-mhd}_1$,  H\"{o}lder's and Young's inequalities, Sobolev's embedding inequality,  Poincar\'{e}'s inequality, \eqref{density-1} and Lemma \ref{Hm-estimate-1},   we get
\begin{align*}
II_2&=-\frac{d}{dt}\int\partial^\alpha\partial_2\mathbf u_2\partial^\alpha\partial^2_2adx-\int(\partial^\alpha\partial_2((1+a)\mathbf u)\cdot\nabla)\partial^\alpha\partial^2_2\mathbf u_2dx\\
&\leq-\frac{d}{dt}\int\partial^\alpha\partial_2\mathbf u_2\partial^\alpha\partial^2_2adx+\big(2\|D^{m-1}\partial_2\mathbf u\|_0+C\|D^{m-1}a\|_{0}\|\partial_2\mathbf u\|_{L^\infty}+C\|Da\|_{L^\infty}\|D^{m-2}\partial_2\mathbf u\|_0\\
&\quad+\|\partial_2a\|_{L^4}\|D^{m+1}\mathbf u\|_{0}+\|D^{m-1}\partial_2a\|_0\|\mathbf u\|_{L^\infty}\big)\|D^{m+1}\partial_2\mathbf u\|_0\\
&\leq-\frac{d}{dt}\int\partial^\alpha\partial_2\mathbf u_2\partial^\alpha\partial^2_2adx+C_{0,1}\|D^{m+1}\partial_2\mathbf u\|^2_0+C\|a,\mathbf u\|_{3}\|D^{m+1}\partial_2\mathbf u,D^{m-1}\partial^2_2a\|^2_0,
\end{align*}
%One has by H\"{o}lder, Young, Sobolev embedding and Poincar\'{e} inequality, Lemma \ref{Hm-estimate-1} and \eqref{density-1}
%\begin{align*}
%II_{2,1}\leq&\big(2\|D^{m-1}\partial_2\mathbf u\|_0+C\|D^{m-1}a\|_{0}\|\partial_2\mathbf u\|_{L^\infty}+C\|Da\|_{L^\infty}\|D^{m-2}\partial_2\mathbf u\|_0\\
%&+\|\partial_2a\|_{L^4}\|D^{m+1}\mathbf u\|_{0}+\|D^{m-1}\partial_2a\|_0\|\mathbf u\|_{L^\infty}\big)\|D^{m+1}\partial_2\mathbf u\|_0\\
%&\leq C_{0,1}\|D^{m+1}\partial_2\mathbf u\|^2_0+C\|a,\mathbf u\|_{3}\|D^{m+1}\partial_2\mathbf u,D^{m-1}\partial^2_2a\|^2_0,
%\end{align*}
where $C_{0,1}$ is a positive constant depending only on $|\mathbb T|$.
%Combining the above equality with \eqref{rho-10-13-1}, yields
%$$II_2\leq -R\frac{d}{dt}\int\partial^\alpha\partial_2\mathbf u_2\partial^\alpha\partial^2_2ad+C_{0,1}\|D^{m+1}\partial_2\mathbf u\|^2_0+C\|a,\mathbf u\|_{3}\|D^{m+1}\partial_2\mathbf u,D^{m-1}\partial^2_2a\|^2_0.$$
For $II_3$, it follows from H\"{o}lder's and Young's inequalities, Sobolev's embedding inequality,  Poincar\'{e}'s inequality and \eqref{product}, that
\begin{align*}
II_{3}&\lesssim\big(\|\mathbf u\|_{L^\infty}\|D^{m+1}{\mathbf u}\|_0+{\|D^m\mathbf u\|_0\|D\mathbf u\|_{L^\infty}}\big)\|D^{m-1}\partial^2_2a\|_{{0}}\\
&\lesssim\|\mathbf u\|_{m+1}\|D^{m-1}\partial^2_2a,D^{m+1}\mathbf u\|^2_0.
\end{align*}
Employing H\"{o}lder's and Young's inequalities,   Poincar\'{e}'s inequality, leads to
$$II_4\leq C_{0,2}\|D^{m+1}\partial_2\theta\|_0^2+\frac{R}{4}\|\partial^{\alpha}\partial^2_2a\|^2_0,$$
where $C_{0,2}$ is a positive constant depending only on $|\mathbb T|,R$.
 By a similar derivation of $II_{3}$, we conclude from \eqref{density-1} and \eqref{h-junzhi}, that
\begin{align*}II_{5}&\lesssim(1+\|a\|_m)\|\mathbf h_1\|_m\|D\mathbf h\|_m\|D^{m-1}\partial^2_2a\|_0\\&\lesssim(1+\|a\|_m)\|D^{m-1}\partial_2\mathbf h\|_0\|D\mathbf h\|_m\|D^{m-1}\partial^2_2a\|_0\\&\lesssim(1+\|a\|_m)\|D\mathbf h\|_{m}\|D^{m-1}\partial^2_2a,D^{m-1}\partial^2_2h\|^2_0.\end{align*}
For $II_{6}$, using H\"{o}lder's and Young's inequalities, Sobolev's embedding inequality,  Poincar\'{e}'s inequality,  \eqref{density-1} and Lemma \ref{Hm-estimate-1},  we obtain
\begin{align*}
II_{6}&\leq \Big(\|\frac{\mu\Delta\partial^\alpha\partial_2\mathbf u_2+(\mu+\lambda)\partial^\alpha\partial^2_2\mathrm{div}\mathbf u}{1+a}\|_0+\|[\partial^\alpha\partial_2,\frac{1}{1+a}](\mu\Delta\mathbf u_2+(\mu+\lambda)\partial_2\mathrm{div}\mathbf u)\|_0\Big)\|D^\alpha\partial^2_2a\|_0\\
&\leq \Big(2(2\mu+\lambda)\|D^{m+1}\partial_2\mathbf u\|_0+C\|Da\|_{L^\infty}\|D^{m+1}\mathbf u\|+C\|a\|_{m+1}\|D^2\mathbf u\|_{L^4}\Big)\|D^{m-1}\partial^2_2a\|_0\\
&\leq C_{0,3}\|D^{m+1}\partial_2\mathbf u\|^2_0+\|a\|_{m+1}\|D^{m+1}\mathbf u,D^{m-1}\partial^2_2a\|_0^2+\frac R4\|D^\alpha\partial^2_2a\|^2_0,
\end{align*}
where $C_{0,3}$ is a positive constant depending only on $\mu,\lambda,R$. Eventually, plugging the above  estimates of $II_{j}(j=1,2,\cdots,6)$ into \eqref{rho-10-13-1}, summing over all multi-indices $\alpha$ with $|\alpha|=m-1$, {and multiplying the resulting inequality by $\frac2R$,} we deduce the desired conclusion. This completes the proof of Lemma \ref{dissipation-rho}.
\end{proof}
\begin{lemma}\label{dissipation-h} Let $(a,\mathbf u,\theta,\mathbf h)$ be smooth solutions to  the system $\eqref{heat-mhd}$
  on $\mathbb T^2\times[0,T)$ for some $T>0$, and suppose that \eqref{density-1}, \eqref{theta-1} and \eqref{junzhi} hold.  Then,  for any $m\geq2$,
 \begin{align}
&2\frac{d}{dt}\sum_{|\alpha|=m-1}\int\partial^\alpha\partial^2_2\mathbf u\cdot\partial^\alpha\partial_2\mathbf hdx+\|D^{m-1}\partial^2_2\mathbf h\|_0^2\leq C_1\|D^{m+1}\partial_2\mathbf u,D^{m+1}\partial_2\theta\|^2_0\nonumber\\
&\quad+(1+\|a\|_{m+1})\|a,\mathbf u,\mathbf h\|_{m+1}\|D^{m-1}\partial_2^2a,D^{m+1}\partial_2\mathbf u,D^{m+1}\mathbf u,D^{m+2}\theta,D^{m-1}\partial^2_2\mathbf h\|_0^2,\label{10-16-partial2h}
\end{align}
where {$C_{1}$} is a positive constant depending only on $|\mathbb T|,R,\mu,\lambda,m.$
\end{lemma}
\begin{proof}
By \eqref{heat-mhd}$_2$, we have
\begin{align*}
\partial_2\mathbf h= &{\frac{a\partial_2\mathbf h}{1+a}}+\partial_t\mathbf u+(\mathbf u\cdot\nabla)\mathbf u+{R\nabla\theta}+R\nabla \ln(1+a)(1+\theta)\\
&+\frac{\nabla\mathbf h_2}{1+a}+\frac{\nabla(|\mathbf h|^2)}{2(1+a)}-\frac{(\mathbf h\cdot\nabla)\mathbf h}{1+a}-\frac{\mu\Delta\mathbf u+(\mu+\lambda)\nabla\mathrm{div}\mathbf u}{1+a}.
\end{align*}
Given any multi-indices $\alpha$ with $|\alpha|=m-1$,  applying $\partial^{\alpha}\partial_2$ to the above identity, and taking $L^2$ inner product
with $\partial^{\alpha}\partial_2^2\mathbf h$ yield
\begin{equation}\|\partial^{\alpha}\partial^2_2\mathbf h\|_0^2=\sum_{j=1}^{9}III_j\label{hh-10-14-1}
\end{equation}
with
\begin{align*}
&{III_1:=\int}\partial^\alpha\partial_2\big(\frac{a\partial_2\mathbf h}{1+a}\big)\cdot\partial^\alpha\partial_2^2\mathbf hdx,\quad III_2:=\int\partial^\alpha\partial_2\partial_t\mathbf u\cdot\partial^\alpha\partial_2^2\mathbf hdx,\\
&III_3:=\int\partial^\alpha\partial_2\big((\mathbf u\cdot\nabla)\mathbf u\big)\cdot\partial^\alpha\partial_2^2\mathbf hdx,\quad III_4:=R\int\partial^\alpha\partial_2{\nabla\theta}\cdot\partial^\alpha\partial_2^2\mathbf hdx,\\
&III_5:=R\int\partial^\alpha\partial_2\big(\nabla \ln(1+a)(1+\theta)\big)\cdot\partial^\alpha\partial_2^2\mathbf hdx,\quad
III_6:=\int\partial^\alpha\partial_2\big(\frac{\nabla\mathbf h_2}{1+a}\big)\cdot\partial^\alpha\partial_2^2\mathbf hdx,\\
&III_7:=\int\partial^\alpha\partial_2\big(\frac{\nabla(|\mathbf h|^2)}{2(1+a)}\big)\cdot\partial^\alpha\partial_2^2\mathbf hdx,\quad
III_8:=-\int\partial^\alpha\partial_2\big(\frac{(\mathbf h\cdot\nabla)\mathbf h}{1+a}\big)\cdot\partial^\alpha\partial_2^2\mathbf hdx,\\
&III_9:=-\int\partial^\alpha\partial_2\big(\frac{\mu\Delta\mathbf u+(\mu+\lambda)\nabla\mathrm{div}\mathbf u}{1+a}\big)\cdot\partial^\alpha\partial_2^2\mathbf hdx.
\end{align*}
In what follows, we bound  $III_j(j=1,2,\cdots,9)$  successively. For $III_1$, by  H\"{o}lder's and Young's inequalities, Sobolev's embedding inequality,  Poincar\'{e}'s inequality, \eqref{density-1}  and Lemma \ref{Hm-estimate-1},   we have
\begin{align*}
III_1&\lesssim(\|a\|_{L^\infty}\|D^{m-1}\partial_2^2\mathbf h\|_0+\|D^m\big(\frac{a}{1+a}\big))\|_{0}\|\partial^2_2h\|_{L^4}+\|Da\|_{L^\infty}\|D^{m-1}\partial_2\mathbf h\|_0\\
&\quad+\|D^{m}\big(\frac{\partial_2a}{(1+a)^2}\big)\|_0\|\partial_2\mathbf h\|_{L^4}\\&\lesssim(1+\|a\|_{m+1})\|a\|_{m+1}\|D^{m-1}\partial^2_2\mathbf h\|_0^2.\\
\end{align*}
For $III_2$, by integration by parts and \eqref{heat-mhd}$_4$, we easily get
\begin{align}
III_2&=-\frac{d}{dt}\int\partial^\alpha\partial^2_2\mathbf u\cdot\partial^\alpha\partial_2\mathbf hdx+\int\partial^\alpha\partial^2_2\mathbf u\cdot\partial^\alpha\partial_2\partial_t\mathbf hdx\nonumber\\
&=-\frac{d}{dt}\int\partial^\alpha\partial^2_2\mathbf u\cdot\partial^\alpha\partial_2\mathbf hdx+\sum_{l=1}^{{4}}III_{2,l}\label{hh-10-14-2}
\end{align}
with
\begin{align*}
III_{2,1}:=-\int\partial^\alpha\partial^2_2\mathbf u\cdot\partial^\alpha\partial_2((\mathbf u\cdot\nabla)\mathbf h)dx,\quad III_{2,2}:=-\int\partial^\alpha\partial^2_2\mathbf u\cdot\partial^\alpha\partial_2(\mathbf h\mathrm{div}\mathbf u)dx,\\
III_{2,3}:=\int\partial^\alpha\partial^2_2\mathbf u\cdot\partial^\alpha\partial_2((\mathbf h\cdot\nabla)\mathbf u)dx,\quad III_{2,4}:=\int\partial^\alpha\partial^2_2\mathbf u\cdot\partial^\alpha\partial_2(\partial_2\mathbf u-\mathrm{div}\mathbf u\mathbf e_2)dx.
\end{align*}
For $III_{2,1}$, choosing $\alpha_k>0$,  and then  using integration by parts and  H\"{o}lder's and Young's inequalities, Sobolev's embedding inequality,  Poincar\'{e}'s inequality and Lemma \ref{Hm-estimate-1}, we infer that
\begin{align*}
III_{2,1}&=\int\partial^{\alpha+\mathbf e_k}\partial^2_2\mathbf u\cdot\partial^{\alpha-\mathbf e_k}\partial_2((\mathbf u\cdot\nabla)\mathbf h)dx\\
&\lesssim\|D^{m+1}\partial_2\mathbf u\|_0\big(\|(\mathbf u\cdot\nabla)\partial^{\alpha-\mathbf e_k}\partial_2\mathbf h\|_0+\|[\partial^{\alpha-\mathbf e_k}\partial_2,\mathbf u\cdot\nabla]\mathbf h\|_0\big)\\
&\lesssim\|D^{m+1}\partial_2\mathbf u\|_0\big(\|\mathbf u\|_{L^\infty}\|D^{m-1}\partial_2\mathbf h\|_0+\|D\mathbf u\|_{L^\infty}\|D^{m-1}\mathbf h\|_0+\|D^{m-1}\mathbf u\|_{0}\|D\mathbf h\|_{L^\infty}\big)\\
&\lesssim\|\mathbf u,\mathbf h\|_3\|D^{m+1}\partial_2\mathbf u,D^{m+1}\mathbf u,D^{m-1}\partial^2_2\mathbf h\|_0^2.
\end{align*}
For $III_{2,2}$, employing H\"{o}lder's and Young's inequalities, Sobolev's embedding inequality,  Poincar\'{e}'s inequality and Lemma \ref{Hm-estimate-1}, leads to
\begin{align*}
III_{2,2}\lesssim\|D^{m+1}\mathbf u\|_0(\|\mathbf h\|_{L^\infty}\|D^m\mathrm{div}\mathbf u\|_0+\|D^{m}\mathbf h\|_{0}\|\mathrm{div}\mathbf u\|_{L^\infty})\lesssim\|\mathbf h\|_3\|D^{m+1}\mathbf u\|_0^2.
\end{align*}
And a similar argument  gives
$$III_{2,3}\lesssim\|\mathbf h\|_3\|D^{m+1}\mathbf u\|_0^2.$$
According to Poincar\'{e}'s inequality,  we consequently arrive at
$$III_{2,4}\leq C_{1,1}\|D^{m+1}\partial_2\mathbf u\|^2_0,$$
where $C_{1,1}$ is a positive constant depending only on $|\mathbb T|.$ Then, combining the above estimates of $III_{2,1}-III_{2,4}$ with \eqref{hh-10-14-2}, yields
$$III_2\leq-\frac{d}{dt}\int\partial^\alpha\partial^2_2\mathbf u\cdot\partial^\alpha\partial_2\mathbf hdx+C_{1,1}\|D^{m+1}\partial_2\mathbf u\|^2_0+\|\mathbf u,\mathbf h\|_3\|D^{m+1}\partial_2\mathbf u,D^{m+1}\mathbf u,D^{m-1}\partial^2_2\mathbf h\|_0^2.$$
For $III_3$,  by a similar derivation in $III_{2,2}$, we also  have
$$III_{3}\lesssim\|\mathbf u\|_3\|D^{m+1}\mathbf u\|_0\|D^{m-1}\partial_2^2\mathbf h\|_0\lesssim\|\mathbf u\|_3\|D^{m+1}\mathbf u,D^{m-1}\partial_2^2\mathbf h\|^2_0.$$
{For $III_4$,  using integrating by parts and  $\mathrm{div}\mathbf h=0$, we find that this term vanishes identically
$$III_{4}=-R\int\partial^\alpha\partial_2\theta\,\partial^\alpha\partial_2^2\mathrm{div}\mathbf h\,dx=0.$$}
For $III_5$, by integration by parts and $\mathrm{div}\mathbf h=0$ {again (the gradient part $R\nabla\big(\ln(1+a)(1+\theta)\big)$ does not contribute)}, we obtain
\begin{align*}
III_{5}=-R\int\partial^\alpha\partial_2\big(\nabla\theta\ln(1+a)\big)\cdot\partial^\alpha\partial_2^2\mathbf hdx\lesssim\|a\|_{{m+1}}\|D^{m+1}\theta,D^{m-1}\partial_2^2\mathbf h\|^2_0.
\end{align*}
For $III_6$, employing  integration by parts, $\mathrm{div}\mathbf h=0$, H\"{o}lder's and Young's inequalities, Sobolev's embedding inequality,  Poincar\'{e}'s inequality, \eqref{density-1} and  Lemma \ref{Hm-estimate-1},  we infer that
\begin{align*}
III_{6}&=-\int\partial^\alpha\big(\frac{\partial_2a\nabla\mathbf h_2}{(1+a)^2}\big)\cdot\partial^\alpha\partial_2^2\mathbf hdx+\int\partial^\alpha\big(\frac{\nabla a\partial_2\mathbf h_2}{(1+a)^2}\big)\cdot\partial^\alpha\partial_2^2\mathbf hdx\\
&\lesssim\Big(\|D^{m-1}\big(\frac{\partial_2a}{(1+a)^2}\big)\|_0\|D\mathbf h\|_{L^\infty}+\|\partial_2a\|_{L^4}\|D^{m}\mathbf h\|_{L^4}+\|D^{m-1}\big(\frac{\nabla a}{(1+a)^2}\big)\|_{L^4}\|\partial_2\mathbf h\|_{L^4}\\
&\quad+\|Da\|_{L^\infty}\|D^{m-1}\partial_2\mathbf h\|_0\Big)\|D^{m-1}\partial_2^2\mathbf h\|_0\\
&\lesssim\big((1+\|a\|_{m+1})\|\mathbf h\|_{m+1}\|D^{m-1}\partial_2a\|_0+(1+\|a\|_{m+1})\|a\|_{m+1}\|D^{m-1}\partial_2\mathbf h\|_0\big)\|D^{m-1}\partial_2^2\mathbf h\|_0\\
&\lesssim(1+\|a\|_{m+1})\|a,\mathbf h\|_{m+1}\|D^{m-1}\partial_2^2a,D^{m-1}\partial_2^2\mathbf h\|^2_0.
\end{align*}
In the same way, we can get
\begin{align*}
III_{7}&\lesssim\Big(\|D^{m-1}\big(\frac{\partial_2a}{(1+a)^2}\big)\|_0\|D(|\mathbf h|^2)\|_{L^\infty}+\|\partial_2a\|_{L^4}\|D^{m}(|\mathbf h|^2)\|_{L^4}\\
&\quad+\|D^{m-1}\big(\frac{\nabla a}{(1+a)^2}\big)\|_{L^4}\|\partial_2(|\mathbf h|^2)\|_{L^4}+\|Da\|_{L^\infty}\|D^{m-1}\partial_2(|\mathbf h|^2)\|_0\Big)\|D^{m-1}\partial_2^2\mathbf h\|_0\\
&\lesssim\big((1+\|a\|_{m+1})\|\mathbf h\|^2_{m+1}\|D^{m-1}\partial_2a\|_0+(1+\|a\|_{m+1})\|a\|_{m+1}\|h\|_{m+1}\|D^{m-1}\partial_2\mathbf h\|_0\big)\|D^{m-1}\partial_2^2\mathbf h\|_0\\
&\lesssim(1+\|a\|_{m+1})\|a,\mathbf h\|^2_{m+1}\|D^{m-1}\partial_2^2a,D^{m-1}\partial_2^2\mathbf h\|^2_0.
\end{align*}
For $III_{8}$, it then  follows from H\"{o}lder's and Young's inequalities, Sobolev's embedding inequality,  Poincar\'{e}'s inequality, {\eqref{h-junzhi}}, \eqref{density-1}  and Lemma \ref{Hm-estimate-1},  that
\begin{align*}
III_{8}&\lesssim\Big(\|\partial^\alpha\partial_2\big(\frac{\mathbf h_1\partial_1\mathbf h}{1+a}\big)\|_0+\|\frac{\mathbf h_2\partial^\alpha\partial^2_2\mathbf h}{1+a}\|_0+\|[\partial^\alpha\partial_2,\frac{\mathbf h_2}{1+a}]\partial_2\mathbf h\|_0\Big)\|D^{m-1}\partial_2^2\mathbf h\|_0\\
&\lesssim\Big(\|D^m\mathbf h_1\|_0\|D\mathbf h\|_{L^\infty}+(1+\|a\|_{m+1})\|\mathbf h\|_{m+1}\|\mathbf h_1\|_{L^\infty}+\|\mathbf h\|_{L^\infty}\|D^{m-1}\partial_2^2\mathbf h\|_0\\
&\quad+\|D^m\big(\frac{\mathbf h_2}{1+a}\big)\|_{L^4}\|\partial_2\mathbf h\|_{L^4}+\|D\mathbf h_2\|_{L^\infty}\|D^{m-1}\partial_2\mathbf h\|_0\Big)\|D^{m-1}\partial_2^2\mathbf h\|_0\\
&\lesssim(1+\|a\|_{m+1})\|\mathbf h\|_{m+1}(\|\mathbf h_1\|_m+\|D^{m-1}\partial_2^2\mathbf h\|_0)\|D^{m-1}\partial_2^2\mathbf h\|_0\\
&\lesssim(1+\|a\|_{m+1})\|\mathbf h\|_{m+1}\|D^{m-1}\partial_2^2\mathbf h\|^2_0.
\end{align*}
For $III_{9}$, by H\"{o}lder's and Young's inequalities, Sobolev's embedding inequality,  Poincar\'{e}'s inequality, \eqref{density-1}  and Lemma \ref{Hm-estimate-1},  we get
\begin{align*}
III_{9}&\leq \Big((4\mu+2\lambda)\|D^{m+1}\partial_2\mathbf u\|_0+\|[\partial^\alpha\partial_2,\frac{1}{1+a}](\mu\Delta\mathbf u+(\mu+\lambda)\nabla\mathrm{div}\mathbf u)\|_0\Big)\|\partial^\alpha\partial_2^2\mathbf h\|_0\\
&\leq C_{1,3}\|D^{m+1}\partial_2\mathbf u\|^2_0+\frac{1}{4}\|\partial^\alpha\partial_2^2\mathbf h\|^2_0+\|a\|_{m+1}\|D^{m+1}\mathbf u,D^{m-1}\partial_2^2\mathbf h\|^2_0,
\end{align*}
where $C_{1,3}$ is a positive constant depending only on $\mu,\lambda.$ In the end, plugging the above estimates of $III_{j}(j=1,2,\cdots,9)$ into {\eqref{hh-10-14-1}}, and summing over all multi-indices $\alpha$ with $|\alpha|=m-1$, we conclude the desired conclusion. This completes the proof of Lemma \ref{dissipation-h}.
\end{proof}
\subsection{The dissipation for $\mathrm{div}\mathbf u$ and $\Phi$}In this subsection, we introduce the total pressure $\Phi=R\nabla a(1+\theta)+\nabla({\frac{|\mathbf h|^2}{2}+\mathbf h_2})$. By virtue of the coupling  between $\Phi$ and $\mathrm{div}\mathbf u,$ we further establish the dissipation estimates for $\Phi$ and $\mathrm{div}\mathbf u.$ The equations governing $\Phi$ and $\mathrm{div}\mathbf u$ are presented below
\begin{align}\label{Phi-div-mhd}
\left\{
\begin{aligned}
&\partial_{t}\Phi+(\mathbf u\cdot\nabla)\Phi+(R+1)\nabla\mathrm{div}u-\nabla\partial_2\mathbf u_2=-R\nabla\mathrm{div}\mathbf u\theta+\mathcal V_1+R\nabla a\mathcal W_1+\mathcal V_2+\mathcal V_3,\\
&\partial_{t}\mathrm{div}\mathbf u-(2\mu+\lambda)\Delta\mathrm{div}\mathbf u+R\Delta\theta+\mathrm{div}\Phi=\mathcal W_2,
\end{aligned}
\right.
\end{align}
where $\mathcal{V}_p=(\mathcal V_{p,1},\mathcal V_{p,2}) (p=1,2,3)$ are vectors, and $\mathcal{W}_q (q=1,2)$ are scalars, with the expressions given by:
\begin{align*}
\mathcal V_{1,i}:=&-R(1+\theta)\big((\partial_i\mathbf u\cdot\nabla)a+\partial_i(a\mathrm{div}\mathbf u)\big),\\
\mathcal V_{2,i}:=&-(\partial_i\mathbf u\cdot\nabla)\mathbf h_2+\partial_i((\mathbf h\cdot\nabla)\mathbf u_2)-\partial_i(\mathbf h_2\mathrm{div}\mathbf u),\\
\mathcal V_{3,i}:=&-\frac12(\partial_i\mathbf u\cdot\nabla)(|\mathbf h|^2)+\partial_i\big((\mathbf h\cdot\nabla)\mathbf u\cdot\mathbf h+\partial_2\mathbf u\cdot\mathbf h-|\mathbf h|^2\mathrm{div}\mathbf u-\mathbf h_2\mathrm{div}\mathbf u\big),
\end{align*}
with $i=1,2$, and
\begin{align*}
\mathcal W_{1}:=&-\frac{R(1+\theta)\mathrm{div}\mathbf u}{c_\nu}+\frac{\mathcal Q( u)}{c_\nu(1+a)}+\frac{\kappa\Delta\theta}{c_{\nu}(1+a)},\\
\mathcal W_{2}:=&-\mathrm{div}((\mathbf u\cdot\nabla)\mathbf u)+\frac{a}{1+a}\mathrm{div}\Phi+\frac{(\Phi\cdot\nabla)a}{(1+a)^2}-\frac{(2\mu+\lambda)a}{1+a}\Delta\mathrm{div}\mathbf u-\frac{\big((\mu\Delta\mathbf u+(\mu+\lambda)\nabla\mathrm{div}\mathbf u)\cdot\nabla\big)a}{(1+a)^2}\\
&+\sum_{k,i=1}^2\frac{\partial_k\mathbf h_i\partial_i\mathbf h_k}{1+a}-\frac{\big((\mathbf h\cdot\nabla)\mathbf h\cdot\nabla\big)a}{(1+a)^2}-\frac{(\partial_2\mathbf h\cdot\nabla)a}{(1+a)^2}.
\end{align*}
\begin{lemma}\label{dissipation-vvvww} Let $(a,\mathbf u,\theta,\mathbf h)$ be smooth solutions to  the system $\eqref{heat-mhd}$
  on $\mathbb T^2\times[0,T)$ for some $T>0$, and suppose that \eqref{density-1}, \eqref{theta-1} and \eqref{junzhi} hold. Then,  for any $m\geq2$,
 \begin{align}
&\|D^{m}\mathcal{V}_1\|_0\lesssim(1+\|\theta\|_{m+1})\|a\|_{m+1}\|D^{m+1}\mathbf u,D^{m+1}\mathrm{div}\mathbf u\|_0,\label{V-1}
\end{align}
 \begin{align}
&\|D^{m}\mathcal{V}_2\|_0\lesssim\|\mathbf h\|_{m+1}\|D^{m+1}\mathbf u,D^{m+1}\partial_2\mathbf u,D^{m+1}\mathrm{div}\mathbf u\|_0+\|D\partial_2\mathbf h\|_{0}\|D^{m+2}\mathbf u\|_0,\label{V-2}
\end{align}
 \begin{align}
&\|D^{m}\mathcal{V}_3\|_0\lesssim(1+\|\mathbf h\|_{m+1})\big(\|\mathbf h\|_{m+1}\|D^{m+1}\mathbf u,D^{m+1}\partial_2\mathbf u,D^{m+1}\mathrm{div}\mathbf u\|_0+\|D\partial_2\mathbf h\|_{0}\|D^{m+2}\mathbf u\|_0\big),\label{V-3}
\end{align}
 \begin{align}
&\|\mathcal{W}_1\|_m\lesssim(1+\|a,\mathbf u,\theta\|_{m+1})^2\|D^{m+1}\mathbf u,D^{m+2}\theta\|_0,\label{W-1}
\end{align}
 \begin{align}
&\|\mathcal{W}_2\|_{m-1}\lesssim(1+\|a\|_{m+1})\|a,\mathbf u,\mathbf h\|_{m+1}\|D^{m+1}\mathbf u,D^{m+1}\mathrm{div}\mathbf u,D^{m-1}\partial_2^2\mathbf h,D^{m+2}\theta,D^{m}\Phi\|_0.\label{W-2}
\end{align}
\end{lemma}
\begin{proof}
By the definition of $\mathcal{V}_1$,  H\"{o}lder's  inequality, Sobolev's embedding inequality,  Poincar\'{e}'s inequality and Lemma \ref{Hm-estimate-1}, we have
\begin{align*}
\|D^{m}\mathcal{V}_1\|_0&\lesssim(1+\|\theta\|_{m+1})\|(\nabla\mathbf u\cdot\nabla)a+\nabla(a\mathrm{div}\mathbf u)\|_m\\
&\lesssim(1+\|\theta\|_{m+1})\big(\|D\mathbf u\|_{L^\infty}\|Da\|_{m}+\|D\mathbf u\|_m\|Da\|_{L^\infty}+\|a\|_{L^\infty}\|\mathrm{div}\mathbf u\|_{m+1}\\
&\quad+\|a\|_{m+1}\|\mathrm{div}\mathbf u\|_{L^\infty}\big)\\&\lesssim(1+\|\theta\|_{m+1})\|a\|_{m+1}\|D^{m+1}\mathbf u,D^{m+1}\mathrm{div}\mathbf u\|_0.
\end{align*}
It follows from the definition of $\mathcal{V}_2$,  H\"{o}lder's  inequality, Sobolev's embedding inequality,  Poincar\'{e}'s inequality, \eqref{h-junzhi} and Lemma \ref{Hm-estimate-1},   that
\begin{align*}
\|D^{m}\mathcal{V}_2\|_0&\lesssim\|D\mathbf u\|_{L^\infty}\|D^{m+1}\mathbf h\|_0+\|D\mathbf h\|_{L^\infty}\|D^{m+1}\mathbf u\|_0+\|\mathbf h_1\|_{L^\infty}\|D^{m+2}\mathbf u\|_0\\
&\quad+\|D^{m+1}\mathbf h_1\|_{0}\|D\mathbf u\|_{L^\infty}+\|\mathbf h_2\|_{L^\infty}\|D^{m+1}\partial_2\mathbf u\|_0+\|D^{m+1}\mathbf h_2\|_{0}\|\partial_2\mathbf u\|_{L^\infty}\\
&\quad+\|\mathbf h_2\|_{L^\infty}\|D^{m+1}\mathrm{div}\mathbf u\|_0+\|D^{m+1}\mathbf h_2\|_{0}\|\mathrm{div}\mathbf u\|_{L^\infty}\\
&\lesssim\|\mathbf h\|_{m+1}\|D^{m+1}\mathbf u,D^{m+1}\partial_2\mathbf u,D^{m+1}\mathrm{div}\mathbf u\|_0+\|\mathbf h_1\|_{2}\|D^{m+2}\mathbf u\|_0\\
&\lesssim\|\mathbf h\|_{m+1}\|D^{m+1}\mathbf u,D^{m+1}\partial_2\mathbf u,D^{m+1}\mathrm{div}\mathbf u\|_0+\|D\partial_2\mathbf h\|_{0}\|D^{m+2}\mathbf u\|_0.
\end{align*}
Similarly, we  also obtain
\begin{align*}
\|D^{m}\mathcal{V}_3\|_0&\lesssim\|D\mathbf u\|_{L^\infty}\|D^{m+1}(|\mathbf h|^2)\|_0+\|D(|\mathbf h|^2)\|_{L^\infty}\|D^{m+1}\mathbf u\|_0+\|\mathbf h_1\|_{L^\infty}\|D^{m+1}(\partial_1\mathbf u\cdot\mathbf h)\|_0\\
&\quad+\|D^{m+1}\mathbf h_1\|_{0}\|\partial_1\mathbf u\cdot\mathbf h\|_{L^\infty}+\|\mathbf h_2\|_{L^\infty}\|D^{m+1}(\partial_2\mathbf u\cdot\mathbf h)\|_0+\|D^{m+1}\mathbf h_2\|_{0}\|\partial_2\mathbf u\cdot\mathbf h\|_{L^\infty}\\
&\quad+\|\partial_2\mathbf u\|_{L^\infty}\|\mathbf h\|_{m+1}+\|D^{m+1}\partial_2\mathbf u\|_0\|\mathbf h\|_{L^\infty}+\||\mathbf h|^2\|_{L^\infty}\|D^{m+1}\mathrm{div}\mathbf u\|_0\\
&\quad+\|D^{m+1}(|\mathbf h|^2)\|_{0}\|\mathrm{div}\mathbf u\|_{L^\infty}+\|\mathbf h_2\|_{L^\infty}\|D^{m+1}\mathrm{div}\mathbf u\|_0+\|D^{m+1}\mathbf h_2\|_{0}\|\mathrm{div}\mathbf u\|_{L^\infty}\\
&\lesssim(1+\|\mathbf h\|_{m+1})\big(\|\mathbf h\|_{m+1}\|D^{m+1}\mathbf u,D^{m+1}\partial_2\mathbf u,D^{m+1}\mathrm{div}\mathbf u\|_0+\|D\partial_2\mathbf h\|_{0}\|D^{m+2}\mathbf u\|_0\big).
\end{align*}
Employing H\"{o}lder's inequality, Sobolev's embedding inequality,  Poincar\'{e}'s inequality,  \eqref{density-1}, \eqref{theta-1},   Lemmas \ref{Hm-estimate-1} and \ref{nolinear-estimate}, we conclude that
\begin{align*}
\|\mathcal{W}_1\|_m&\lesssim(1+\|\theta\|_{m+1})\|\mathrm{div}\mathbf u\|_{m}+(1+\|a\|_{m+1})(\|D\mathbf u\|^2_{m}+\|D^2\theta\|_{{m}})\\
&\lesssim(1+\|a,\mathbf u,\theta\|_{m+1})^2\|D^{m+1}\mathbf u,D^{m+2}\theta\|_0,
\end{align*}
and
\begin{align*}
&\|\mathcal{W}_2\|_{m-1}\\&\lesssim\|\mathbf u\|_{m+1}\|D\mathbf u\|_{m}+\|a\|_{L^\infty}\|D\Phi\|_{m-1}+\|\frac{a}{1+a}\|_{W^{m-1,4}}\|D\Phi\|_{L^4}+(1+\|a\|_{m+1})\|a\|_{m+1}\|\Phi\|_{m}\\
&\quad+\|a\|_{L^\infty}\|D^2\mathrm{div}\mathbf u\|_{m-1}+\|\frac{a}{1+a}\|_{W^{m-1,4}}\|D^2\mathrm{div}\mathbf u\|_{L^4}+\|Da\|_{L^\infty}\|D^2\mathbf u\|_{m-1}+\|\partial_2\mathbf h\|_{m-1}\|\nabla a\|_{L^\infty}\\
&\quad+\|D^2\mathbf u\|_{L^4}\|\frac{a}{1+a}\|_{W^{m,4}}+\sum_{k,i=1}^2\big(\|\partial_k\mathbf h_i\partial_i\mathbf h_k\|_{m-1}+\|\frac{1}{1+a}\|_{W^{m-1,4}}\|\partial_k\mathbf h_i\partial_i\mathbf h_k\|_{L^4}\big)\\
&\quad+\|(\mathbf h\cdot\nabla)\mathbf h\|_{m-1}\|\nabla a\|_{L^\infty}+\|(\mathbf h\cdot\nabla)\mathbf h\|_{L^4}\|\frac{a}{1+a}\|_{W^{m,4}}+\|\partial_2\mathbf h\|_{L^4}\|\frac{a}{1+a}\|_{W^{m,4}}\\
&\lesssim(1+\|a\|_{m+1})\Big(\|a,\mathbf u\|_{m+1}(\|D\mathbf u\|_{m}+\|\Phi\|_m+\|D^2\mathrm{div}\mathbf u,D^2\mathbf u\|_{m-1})+\sum_{k,i=1}^2\|\partial_k\mathbf h_i\partial_i\mathbf h_k\|_{m-1}\Big)\\
&\quad +\|a\|_{m+1}(\|(\mathbf h\cdot\nabla)\mathbf h\|_{m-1}+\|\partial_2 \mathbf h\|_{m-1})\\
&\lesssim(1+\|a\|_{m+1})\|a,\mathbf u,\mathbf h\|_{m+1}\|D^{m+1}\mathbf u,D^{m+1}\mathrm{div}\mathbf u,D^{m-1}\partial_2^2\mathbf h,D^{m+2}\theta,D^{m}\Phi\|_0.
\end{align*}
Here, we have used  \eqref{h-junzhi} and  the following facts
\begin{align}
\sum_{k,i=1}^2\|\partial_k\mathbf h_i\partial_i\mathbf h_k\|_{m-1}&\lesssim\|D\mathbf h_1\|_{m-1}\|\partial_1\mathbf h\|_{L^\infty}+\|D\mathbf h_1\|_{L^4}\|\|\partial_1\mathbf h\|_{L^\infty}+\|D\mathbf h_2\|_{W^{m-1,4}}\|\partial_2\mathbf h\|_{L^4}\nonumber\\
&\quad+\|D\mathbf h_2\|_{L^\infty}\|\|\partial_2\mathbf h\|_{m-1}\nonumber\\
&\lesssim\|\mathbf h\|_{m+1}(\|D\mathbf h_1\|_{m-1}+\|\partial_2\mathbf h\|_{m-1})\lesssim\|\mathbf h\|_{m+1}\|D^{m-1}\partial^2_2\mathbf h\|_0,\label{2026-10-16-1}
\end{align}
\begin{align}
\|(\mathbf h\cdot\nabla)\mathbf h\|_{m-1}&\lesssim\|\mathbf h_1\|_{L^\infty}\|\partial_1\mathbf h\|_{m-1}+\|\mathbf h_1\|_{m-1}\|\partial_1\mathbf h\|_{L^\infty}\nonumber\\&\quad+\|\mathbf h_2\|_{L^\infty}\|\partial_2\mathbf h\|_{m-1}+\|\mathbf h_2\|_{W^{m-1,4}}\|\partial_2\mathbf h\|_{L^4}\nonumber\\
&\lesssim\|\mathbf h\|_{m+1}(\|\mathbf h_1\|_{m}+\|\partial_2\mathbf h\|_{m-1})\nonumber\\&\lesssim\|\mathbf h\|_{m+1}\|D^{m-1}\partial^2_2\mathbf h\|_0,\label{2026-10-16-2}
\end{align}
and
\begin{align}
\|\Phi\|_m&\lesssim\|D^m\Phi\|_0+\big|\int\Phi dx\big|\nonumber\\&\lesssim\|D^m\Phi\|_0+{R}\big|\int \nabla\theta a dx\big|\nonumber\\
&\lesssim\|D^m\Phi\|_0+\|a\|_0\|D\theta\|_0.\label{2026-10-18-13}
\end{align}
This completes the proof of Lemma \ref{dissipation-vvvww}.\end{proof}
\begin{lemma}\label{dissipation-divu} Let $(a,\mathbf u,\theta,\mathbf h)$ be smooth solutions to  the system $\eqref{heat-mhd}$
 on $\mathbb T^2\times[0,T)$ for some $T>0$, and suppose that \eqref{density-1}, \eqref{theta-1} and \eqref{junzhi} hold.  Then, for any $m\geq2$,
 \begin{align}
&\frac{d}{dt}\|(R+1)^{-\frac12}D^m\Phi,D^m\mathrm{div}\mathbf u,{\sqrt{c_\nu}}D^{m+1}\theta\|_0^2+2\min\{(2\mu+\lambda),{\kappa}\}\|D^{m+1}\mathrm{div}\mathbf u,D^{m+2}\theta\|_0^2\nonumber\\
&\leq C(1+\|a,\mathbf u,\theta,\mathbf h\|_{m+1})^2\|a,\mathbf u,\theta,\mathbf h\|_{m+1}\|D^{m+1}\mathbf u,D^{m+1}\partial_2\mathbf u,D^{m+1}\mathrm{div}\mathbf u,D^{m+2}\theta,D^{m-1}\partial_2^2\mathbf h,D^m\Phi\|_0^2\nonumber\\
&\quad+2\|D^{m}\Phi\|_0\|D^{m+1}\partial_2\mathbf u\|_0+C(1+\|\mathbf h\|_{m+1})\|D\partial_2\mathbf h\|_0\|D^{m+2}\mathbf u\|_0\|D^m\Phi\|_0,\label{10-16-1-divu}
\end{align}
where   {$C$}  is a positive constant depending only on $|\mathbb T|,R,\mu,\lambda{,c_\nu,\kappa,m}$.
\end{lemma}
\begin{proof}
Given any multi-indices $\alpha$ with $|\alpha|=m$,  applying $\partial^{\alpha}$ to \eqref{Phi-div-mhd}, taking $L^2$ inner product
with $(\frac{\partial^\alpha\Phi}{R+1},\partial^{\alpha}\mathrm{div}\mathbf u)$ and  using integration by parts, yield
\begin{align}
&\frac{d}{2dt}\|\partial^{\alpha}((R+1)^{-\frac12}\Phi,\mathrm{div}\mathbf u)\|_0^2+(2\mu+\lambda)\|D\partial^\alpha\mathrm{div}\mathbf u\|_0^2\\&=-\frac{1}{R+1}\int\partial^\alpha((\mathbf u\cdot\nabla)\Phi)\cdot\partial^\alpha\Phi dx+\frac{1}{R+1}\int(\partial^\alpha\Phi\cdot\nabla)\partial^\alpha\partial_2\mathbf u_2dx\nonumber\\
&\quad-\frac{R}{R+1}\int\partial^\alpha(\nabla\mathrm{div}\mathbf u\theta)\cdot\partial^\alpha\Phi dx+\frac{1}{R+1}\int{\partial^\alpha}\big(\mathcal V_1+R\nabla a\mathcal W_1+\mathcal V_2+\mathcal V_3\big)\cdot\partial^\alpha\Phi dx\nonumber\\
&\quad+\int\partial^\alpha\mathcal W_2\cdot\partial^\alpha\mathrm{div}\mathbf udx-R\int\partial^\alpha\Delta\theta\partial^\alpha\mathrm{div}\mathbf udx\\&\triangleq\sum_{j=1}^{6}IV_{j}.\label{10-16-divu}
\end{align}
In what follows, we shall deal with  $IV_j(j=1,2,\cdots,6)$  successively.  First, for $IV_{1}$, $IV_{2}$ and $IV_{3}$, by integration by parts, H\"{o}lder's and Young's inequalities, Sobolev's embedding inequality,  Poincar\'{e}'s inequality and Lemma \ref{Hm-estimate-1}, we obtain
\begin{align*}
IV_{1}&=\frac{1}{2(R+1)}\int\mathrm{div}\mathbf u|\partial^\alpha\Phi|^2dx{-\frac{1}{R+1}}\int[\partial^\alpha,\mathbf u\cdot\nabla]\Phi\cdot\partial^\alpha\Phi dx\\
&\lesssim\|D\mathbf u\|_{L^\infty}\|D^{m}\Phi\|_0^2+(\|D\mathbf u\|_{L^\infty}\|D^{m}\Phi\|_0+\|D^m\mathbf u\|_{L^4}\|D\Phi\|_{L^4})\|D^{m}\Phi\|_0\\
&\lesssim\|D\mathbf u\|_m\|D^m\Phi\|_0^2,
\end{align*}
\begin{align*}
IV_{2}&\leq\|D^{m}\Phi\|_0\|D^{m+1}\partial_2\mathbf u\|_0,
\end{align*}and
\begin{align*}
IV_{3}&\lesssim(\|D^{m+1}\mathrm{div}\mathbf u\|_0\|\theta\|_{L^\infty}+\|D\mathrm{div}\mathbf u\|_{L^\infty}\|\theta\|_{m})\|D^m\Phi\|_0\lesssim\|\theta\|_{m+1}\|D^{m+1}\mathrm{div}\mathbf u,D^m\Phi\|_0^2.
\end{align*}
For $IV_{4}$,  according to \eqref{V-1}-\eqref{W-1}, and then using  H\"{o}lder's and Young's inequalities, Sobolev's embedding inequality, and Lemma \ref{Hm-estimate-1}, we have
\begin{align*}
IV_{4}&\lesssim(\|D^m(\mathcal V_1,\mathcal V_2,\mathcal V_3)\|+\|a\|_{m+1}\|\mathcal W_1\|_m)\|D^m\Phi\|_0\\
&\lesssim(1+\|a,\mathbf u,\theta,\mathbf h\|_{m+1})^2\|a,\mathbf h\|_{m+1}\|D^{m+1}\mathbf u,D^{m+1}\partial_2\mathbf u,D^{m+1}\mathrm{div}\mathbf u,D^{m+2}\theta,D^m\Phi\|_0^2\\
&\quad+(1+\|\mathbf h\|_{m+1})\|D\partial_2\mathbf h\|_0\|D^{m+2}\mathbf u\|_0\|D^m\Phi\|_0.
\end{align*}
For $IV_{5}$, using integration by parts and \eqref{W-2}, leads to
\begin{align*}
IV_{5}&\lesssim\|D^{m-1}\mathcal W_2\|_0\|D^{m+1}\mathrm{div}\mathbf u\|_0\\
&\lesssim(1+\|a\|_{m+1})\|a,\mathbf u,\mathbf h\|_{m+1}\|D^{m+1}\mathbf u,D^{m+1}\mathrm{div}\mathbf u,D^{m-1}\partial_2^2\mathbf h,D^{m+2}\theta,D^{m}\Phi\|^2_0.
\end{align*}
It remains to bound $IV_6$. Employing  $\eqref{heat-mhd}_3$ {in the form
$$R\,\mathrm{div}\mathbf u=-c_\nu(\partial_{t}\theta+(\mathbf u\cdot\nabla)\theta)-R\theta\mathrm{div}\mathbf u+\frac{\kappa\Delta\theta}{1+a}+\frac{\mathcal Q(\mathbf u)}{1+a}$$}
and integration by parts, we get
\begin{align*}
IV_{6}&={\int\partial^\alpha\Delta\theta\partial^\alpha\big(c_\nu(\partial_{t}\theta+(\mathbf u\cdot\nabla)\theta)+R\theta\mathrm{div}\mathbf u-\frac{\kappa\Delta\theta}{1+a}-\frac{\mathcal Q(\mathbf u)}{1+a}\big)dx}\\
&={-\frac{c_\nu}{2}\frac{d}{dt}\|\partial^\alpha D\theta\|^2_0-\kappa\|D^2\partial^\alpha\theta\|_0^2}+IV_{6,1}+IV_{6,2},
\end{align*}
with
\begin{align*}
IV_{6,1}:={\int}\partial^\alpha\Delta\theta\partial^\alpha\big(c_\nu(\mathbf u\cdot\nabla)\theta+R\mathrm{div}\mathbf u\theta-\frac{\mathcal Q(\mathbf u)}{1+a}\big)dx,\quad IV_{6,2}:={\kappa\int}\partial^\alpha\Delta\theta\partial^\alpha\big(\frac{a\Delta\theta}{1+a}\big)dx.
\end{align*}
For $IV_{6,1}$, it follows from H\"{o}lder's and Young's inequalities, Sobolev's embedding inequality,  Poincar\'{e}'s inequality, Lemmas \ref{Hm-estimate-1} and  \ref{nolinear-estimate}, that
\begin{align*}
IV_{6,1}&\lesssim(1+\|a\|_{m+1})\|\mathbf u,\theta\|_{m+1}\|D(\mathbf u,\theta)\|_m\|D^{m+2}\theta\|\\
&\lesssim(1+\|a\|_{m+1})\|\mathbf u,\theta\|_{m+1}\|D^{m+1}\mathbf u,D^{m+2}\theta\|^2_0.
\end{align*}
For $IV_{6,2}$, thanks to  \eqref{density-1}, H\"{o}lder's and Young's inequalities, Sobolev's embedding inequality,  Poincar\'{e}'s inequality, Lemmas \ref{Hm-estimate-1} and  \ref{nolinear-estimate},  we deduce
\begin{align*}
IV_{6,2}&\lesssim\|\frac{a}{1+a}\|_{L^\infty}\|D^2\partial^\alpha\theta\|_0^2+\|[\partial^\alpha,\frac{a}{1+a}]\Delta\theta\|_0\|D^2\partial^\alpha\theta\|_0\\
&\lesssim\|a\|_{m+1}\|D^2\partial^\alpha\theta\|_0^2.
\end{align*}
Therefore,
\begin{align*}
IV_{6}&\lesssim{-\frac{c_\nu}{2}\frac{d}{dt}\|\partial^\alpha D\theta\|^2_0-\kappa\|D^2\partial^\alpha\theta\|_0^2}+\|a\|_{m+1}\|D^2\partial^\alpha\theta\|_0^2
\\&\quad+(1+\|a\|_{m+1})\|\mathbf u,\theta\|_{m+1}\|D^{m+1}\mathbf u,D^{m+2}\theta\|^2_0.
\end{align*}
Eventually, putting the above  estimates of $IV_{j}(j=1,2,\cdots,6)$ into \eqref{10-16-divu}, and summing over all multi-indices $\alpha$ with $|\alpha|=m$, we conclude  the desired conclusion.  This completes the proof of Lemma \ref{dissipation-divu}.
\end{proof}
\begin{lemma}\label{dissipation-Phi}Let $(a,\mathbf u,\theta,\mathbf h)$ be smooth solutions to  the system $\eqref{heat-mhd}$
 on $\mathbb T^2\times[0,T)$ for some $T>0$, and suppose that \eqref{density-1}, \eqref{theta-1} and \eqref{junzhi} hold.  Then, for any $m\geq2$,
 \begin{align}
&2\frac{d}{dt}\sum_{|\alpha|=m}\int\partial^\alpha\mathbf u\cdot\partial^\alpha\Phi dx+\|D^m\Phi\|_0^2\leq C_2\|D^{m+1}\partial_2\mathbf u,D^{m+1}\mathrm{div}\mathbf u,D^{m+2}\theta\|_0^2\nonumber\\
&\quad+C(1+\|a,\mathbf u,\theta,\mathbf h\|_{m+1})^2\|a,\mathbf u,\theta,\mathbf h\|_{m+1}\|D^{m+1}\mathbf u,D^{m+1}\partial_2\mathbf u,D^{m+1}\mathrm{div}\mathbf u,D^{m+2}\theta,D^{m-1}\partial_2^2\mathbf h,D^m\Phi\|_0^2\nonumber\\
&\quad+C(1+\|\mathbf h\|_{m+1})\|D\partial_2\mathbf h\|_{0}\|D^{m+2}\mathbf u\|_0\|D^{m+1}\mathbf u\|_0,\label{10-16-1-Phi}
\end{align}
where $C_2$ is a positive constant depending only on $|\mathbb T|,R,\mu,\lambda,m$.
\end{lemma}
\begin{proof}
It follows from $\eqref{heat-mhd}_2$ that
 \begin{align*}
\Phi=\frac{a}{1+a}\Phi-\partial_t\mathbf u-(\mathbf u\cdot\nabla)\mathbf u-R\nabla\theta+\frac{\mu\Delta\mathbf u+(\mu+\lambda)\nabla\mathrm{div}\mathbf u}{1+a}+\frac{(\mathbf h\cdot\nabla)\mathbf h}{1+a}+\frac{\partial_2\mathbf h}{1+a}.
\end{align*}
Given any multi-indices $\alpha$ with $|\alpha|=m$, applying $\partial^{\alpha}$ to the above equality and performing $L^2$ inner product
with $\partial^\alpha\Phi$, yield
 \begin{align}
\|\partial^\alpha\Phi\|_0^2&=\int\partial^\alpha\big(\frac{a}{1+a}\Phi-(\mathbf u\cdot\nabla)\mathbf u-R\nabla\theta\big)\cdot\partial^\alpha\Phi dx-\int\partial^\alpha\partial_t\mathbf u\cdot\partial^\alpha\Phi dx\nonumber\\
&\quad+\int\partial^\alpha\big(\frac{\mu\Delta\mathbf u+(\mu+\lambda)\nabla\mathrm{div}\mathbf u}{1+a}\big)\cdot\partial^\alpha\Phi dx+\int\partial^\alpha\big(\frac{(\mathbf h\cdot\nabla)\mathbf h}{1+a}\big)\partial^\alpha\Phi dx\nonumber\\
&\quad+\int\partial^\alpha\big(\frac{\partial_2\mathbf h}{1+a}\big)\cdot\partial^\alpha\Phi dx
\\&\triangleq\sum_{j=1}^{5}V_j.\label{10-16-Phi}
\end{align}
Subsequently,  we shall deal with  $V_j(j=1,2,3,4,5)$ term by term. For $V_1$, employing  H\"{o}lder's and Young's inequalities, Sobolev's embedding inequality,  Poincar\'{e}'s inequality, \eqref{density-1}, Lemmas \ref{Hm-estimate-1} and  \ref{nolinear-estimate},  we conclude  that
\begin{align*}
V_1&\leq C(\|a\|_{m+1}\|D^m\Phi\|_0+\|\mathbf u\|_m\|D^{m+1}\mathbf u\|_0)\|D^m{\Phi}\|_0+R\|D^{m+1}\theta\|_0\|\partial^\alpha\Phi\|_0\\
&\leq C\|a,\mathbf u\|_{m+1}\|D^{m+1}\mathbf u,D^m\Phi\|_0^2+C_{2,1}\|D^{m+2}\theta\|^2_0+\frac{1}{8}\|\partial^\alpha\Phi\|^2_0,
\end{align*}
where $C_{2,1}$ is a positive constant  depending only on $|\mathbb T|,R$.
We now  bound $V_2$. By integration by parts and using \eqref{Phi-div-mhd}$_1$, we arrive at
\begin{align}
V_2=-\frac{d}{dt}\int\partial^\alpha\mathbf u\cdot\partial^\alpha\Phi dx+V_{2,1}+V_{2,2}+V_{2,3} \label{10-16-Phi-1}
\end{align}
with
\begin{align*}
&V_{2,1}:=-\int\partial^\alpha\mathbf u\cdot\partial^\alpha\big((\mathbf u\cdot\nabla)\Phi\big)dx,\quad V_{2,2}:=\int\partial^\alpha\mathbf u\cdot\partial^\alpha\big(\nabla\partial_2\mathbf u_2-(R+1)\nabla\mathrm{div}u\big)dx,\end{align*}
and
\begin{align*}&{V_{2,3}:=\int}\partial^\alpha\mathbf u\cdot\partial^\alpha\big(-R\nabla\mathrm{div}\mathbf u\theta+\mathcal V_1+R\nabla a\mathcal W_1+\mathcal V_2+\mathcal V_3\big)dx.
\end{align*}
For $V_{2,1}$, using integration by parts, H\"{o}lder's and Young's inequalities, Sobolev's embedding inequality,  Poincar\'{e}'s inequality, and  Lemma \ref{Hm-estimate-1}, we have
\begin{align*}
V_{2,1}&\lesssim\|D^{m+1}\mathbf u\|_0\|D^{m-1}(\mathbf u\cdot\nabla)\Phi\|_0\\&\lesssim\|D^{m+1}\mathbf u\|_0(\|\mathbf u\|_{L^\infty}\|D^{m}\Phi\|_0+\|\mathbf u\|_{W^{m-1,4}}\|D\Phi\|_{L^4})\\
&\lesssim\|\mathbf u\|_{m+1}\|D^{m+1}\mathbf u,D^m\Phi\|_0^2.
\end{align*}
For $V_{2,2}$, by integration by parts,  H\"{o}lder's and  Young's inequalities,  we deduce that
\begin{align*}
V_{2,2}=\int\partial^\alpha\mathrm{div}\mathbf u\cdot\partial^\alpha\big((R+1)\mathrm{div}u-\partial_2\mathbf u_2\big)dx\leq C_{2,2}\|D^{m+1}\mathrm{div}u,D^{m+1}\partial_2\mathbf u\|_0^2,
\end{align*}
where $C_{2,2}$ is a positive constant  depending only on $|\mathbb T|, R$. We conclude by H\"{o}lder's and Young's inequalities, Sobolev's embedding inequality,  Poincar\'{e}'s inequality,  \eqref{V-1}-\eqref{W-1},  and  Lemma \ref{Hm-estimate-1}, that
\begin{align*}
V_{2,3}&\lesssim\|D^{m}\mathbf u\|_0(\|\theta\|_{m+1}\|D^{m+1}\mathrm{div}\mathbf u\|_0+\|D^{m}(\mathcal V_1,\mathcal V_2,\mathcal V_3)\|_0+\|a\|_{m+1}\|D^m\mathcal W_1\|_0)\\
&\lesssim(1+\|a,\mathbf u,\theta,\mathbf h\|_{m+1})^2\|a,\theta,\mathbf h\|_{m+1}\|D^{m+1}\mathbf u,D^{m+1}\partial_2\mathbf u,D^{m+1}\mathrm{div}\mathbf u,D^{m+2}\theta\|^2_0\\
&\quad+(1+\|\mathbf h\|_{m+1})\|D\partial_2\mathbf h\|_{0}\|D^{m+2}\mathbf u\|_0\|D^{m+1}\mathbf u\|_0.
\end{align*}
Combining the above  estimates of $V_{2,1},V_{2,2}$ and $V_{2,3}$ with \eqref{10-16-Phi-1}, yields
\begin{align*}
V_2&\leq-\frac{d}{dt}\int\partial^\alpha\mathbf u\cdot\partial^\alpha\Phi dx+C_{2,2}\|D^{m+1}(\mathrm{div}u,\partial_2\mathbf u)\|_0^2\\
&\quad+C(1+\|a,\mathbf u,\theta,\mathbf h\|_{m+1})^2\|a,\mathbf u,\theta,\mathbf h\|_{m+1}\|D^{m+1}\mathbf u,D^{m+1}\partial_2\mathbf u,D^{m+1}\mathrm{div}\mathbf u,D^{m+2}\theta,D^m\Phi\|^2_0\\
&\quad+C(1+\|\mathbf h\|_{m+1})\|D\partial_2\mathbf h\|_{0}\|D^{m+2}\mathbf u\|_0\|D^{m+1}\mathbf u\|_0.
\end{align*}
It remains to bound $V_3$. By using  the definition of $\Phi$ and integration by parts,  and choosing $\alpha_k>0$,  we get
\begin{align}
V_3&=-\mu\int\partial^\alpha\big(\frac{\Delta\mathrm{div}\mathbf u}{1+a}+(\Delta\mathbf u\cdot\nabla)((1+a)^{-1})\big)\partial^\alpha(Ra(1+\theta)+\frac{|\mathbf h|^2}{2}+\mathbf h_2)dx\nonumber\\
&\quad-\mu R\int\partial^\alpha\big(\frac{\Delta\mathbf u}{1+a}\big)\partial^\alpha(\nabla\theta a)dx+(\mu+\lambda)\int\partial^\alpha\big(\frac{\nabla\mathrm{div}\mathbf u}{1+a}\big)\cdot\partial^\alpha\Phi dx\nonumber\\
&=\mu\int\partial^{\alpha-\mathbf e_k}\big(\frac{\Delta\mathrm{div}\mathbf u}{1+a}\big)\partial^{\alpha+\mathbf e_k}(Ra(1+\theta)+\frac{|\mathbf h|^2}{2}+\mathbf h_2)dx+V_{3,2}+V_{3,3}\nonumber\\&\triangleq V_{3,1}+V_{3,2}+V_{3,3},\label{10-16-Phi-2}
\end{align}
with
\begin{align*}
&V_{3,1}:=\mu\int\partial^{\alpha-\mathbf e_k}\big(\frac{\Delta\mathrm{div}\mathbf u}{1+a}\big)\partial^{\alpha+\mathbf e_k}(Ra(1+\theta)+\frac{|\mathbf h|^2}{2}+\mathbf h_2)dx,\\
&V_{3,2}:=-\mu\int\partial^\alpha\big((\Delta\mathbf u\cdot\nabla)((1+a)^{-1})\big)\partial^\alpha(Ra(1+\theta)+\frac{|\mathbf h|^2}{2}+\mathbf h_2)dx,\\
&V_{3,3}:=-\mu R\int\partial^\alpha\big(\frac{\Delta\mathbf u}{1+a}\big)\partial^\alpha(\nabla\theta a)dx+(\mu+\lambda)\int\partial^\alpha\big(\frac{\nabla\mathrm{div}\mathbf u}{1+a}\big)\cdot\partial^\alpha\Phi dx.
\end{align*}
For $V_{3,1}$, thanks to  H\"{o}lder's and Young's inequalities, Sobolev's embedding inequality,  Poincar\'{e}'s inequality, \eqref{density-1},  Lemmas \ref{Hm-estimate-1} and \ref{nolinear-estimate}, we deduce  that
\begin{align*}
V_{3,1}&=\mu\int\partial^{\alpha-\mathbf e_k}\big(\frac{\Delta\mathrm{div}\mathbf u}{1+a}\big)\partial^{\alpha}(\Phi_k+Ra\partial_k\theta)dx\\
&\leq \big(2\mu \|D^{m+1}\mathrm{div}\mathbf u\|_0+C\|[\partial^{\alpha-\mathbf e_k},(1+a)^{-1}]\Delta\mathrm{div}\mathbf u\|_{0}\big)\|\partial^{\alpha}\Phi_k,R\partial^{\alpha}(a\partial_k\theta)\|_{0}\\
&\leq2\mu \|D^{m+1}\mathrm{div}\mathbf u\|_0\|\partial^{\alpha}\Phi\|_{0}+C\|a\|_{m+1}\|D^{m+1}\theta\|_0\|D^{m+1}\mathrm{div}\mathbf u\|_0\\
&\quad+C\big(\|Da\|_{L^\infty}\|D^{m+1}\mathrm{div}\mathbf u\|_0+\|\frac{a}{1+a}\|_{W^{m-1,4}}\|\Delta\mathrm{div}\mathbf u\|_{L^4}\big)(\|D^m\Phi\|_0+\|a\|_{m+1}\|D^{m+1}\theta\|_0)\\
&\leq C_{2,3}\|D^{m+1}\mathrm{div}\mathbf u\|^2_0+\frac{1}{8}\|\partial^{\alpha}\Phi\|^2_{0}+C(1+\|a\|_{m+1})\|a\|_{m+1}\|D^{m+1}\mathrm{div}\mathbf u,D^{m+2}\theta,D^m\Phi\|_0^2,
\end{align*}
where $C_{2,3}$ is  a positive constant  depending only on $\mu$. A similar argument gives
\begin{align*}
V_{3,2}&=\mu\int\partial^{\alpha-\mathbf e_k}\big((\Delta\mathbf u\cdot\nabla)((1+a)^{-1})\big)\partial^\alpha(\Phi_k+Ra\partial_k\theta)dx\\
&\lesssim\big(\|D^{m+1}\mathbf u\|_{0}\|Da\|_{L^\infty}+\|\Delta\mathbf u\|_{L^4}\|\frac{a}{1+a}\|_{W^{m,4}}\big)(\|D^m\Phi\|_0+\|a\|_{m+1}\|D^{m+1}\theta\|_0)\\
&\lesssim(1+\|a\|_{m+1})\|a\|_{m+1}\|D^{m+1}\mathbf u,D^{m+2}\theta,D^m\Phi\|_0^2,
\end{align*}
and
\begin{align*}
V_{3,3}&=\mu R\int\partial^{\alpha-\mathbf e_k}\big(\frac{\Delta\mathbf u}{1+a}\big)\partial^{\alpha+\mathbf e_k}(\nabla\theta a)dx+(\mu+\lambda)\int\partial^\alpha\big(\frac{\nabla\mathrm{div}\mathbf u}{1+a}\big)\cdot\partial^\alpha\Phi dx\\
&\leq C\|D^{m-1}\big(\frac{\Delta\mathbf u}{1+a}\big)\|_0\|a\|_{m+1}\|D^{m+2}\theta\|_0+2(\mu+\lambda)\|D^{m+1}\mathrm{div}\mathbf u\|_0\|\partial^\alpha\Phi\|_0\\
&\quad +\|[\partial^\alpha,(1+a)^{-1}]\nabla\mathrm{div}\mathbf u\|_0\|D^{m}\Phi\|_0\\
&\leq C(1+\|a\|_{m+1})\|a\|_m\|D^{m+1}\mathbf u,D^{m+2}\theta,D^{m}\Phi\|_0^2+ C_{2,4}\|D^{m+1}\mathrm{div}\mathbf u\|_0^2+\frac{1}{8}\|\partial^\alpha\Phi\|_0^2,
\end{align*}
where $C_{2,4}$ is a positive constant depending only on $\mu,\lambda$. Substituting the above  estimates of $V_{3,1}-V_{3,3}$ into \eqref{10-16-Phi-2}, yields
\begin{align*}V_3&\leq (C_{2,3}+C_{2,4})\|D^{m+1}\mathrm{div}\mathbf u\|_0^2+\frac{1}{8}\|\partial^\alpha\Phi\|_0^2\\& \quad+C(1+\|a\|_{m+1})\|a\|_{m+1}\|D^{m+1}\mathbf u,D^{m+1}\mathrm{div}\mathbf u,D^{m+2}\theta,D^m\Phi\|_0^2.
\end{align*}
Next, we deal with  $V_4$. By using the definition of $\Phi$, integration by parts, and $\mathrm{div}\mathbf h=0$,   and then  choosing $\alpha_k>0$,  we infer that
\begin{align}
V_4&=-\int\partial^{\alpha-\mathbf e_k}\big(\frac{(\mathbf h\cdot\nabla)\mathbf h}{1+a}\big)\partial^{\alpha+\mathbf e_k}\Phi dx\nonumber\\
&=\sum_{j=1}^2\int\partial^{\alpha-\mathbf e_k}\big(\frac{(\partial_j\mathbf h\cdot\nabla)\mathbf h_j}{1+a}+(\mathbf h\cdot\nabla)\mathbf h_j\partial_j((1+a)^{-1})\big)\partial^\alpha(\Phi_k+{R}a\partial_k\theta) dx\nonumber\\
&\quad+{R}\int\partial^{\alpha-\mathbf e_k}\big(\frac{(\mathbf h\cdot\nabla)\mathbf h}{1+a}\big){\cdot}\partial^{\alpha+\mathbf e_k}(\nabla\theta a)dx\nonumber\\&\triangleq V_{4,1}+V_{4,2}.\label{10-16-Phi-3}
\end{align}
For $V_{4,1}$, recalling \eqref{2026-10-16-1}, \eqref{2026-10-16-2} and using H\"{o}lder's and Young's inequalities, Sobolev's embedding inequality,  Poincar\'{e}'s inequality, \eqref{density-1}, \eqref{h-junzhi},  Lemmas \ref{Hm-estimate-1} and \ref{nolinear-estimate},  we get
\begin{align*}
V_{4,1}&\lesssim\Big(\sum_{j=1}^2\big(\|(\partial_j\mathbf h\cdot\nabla)\mathbf h_j\|_{m-1}+\|(1+a)^{-1}\|_{W^{m-1,4}}\|(\partial_j\mathbf h\cdot\nabla)\mathbf h_j\|_{L^4}\big)+\|(\mathbf h\cdot\nabla)\mathbf h\|_{m-1}\|Da\|_{L^\infty}\\
&\quad+\|(\mathbf h\cdot\nabla)\mathbf h\|_{L^4}\|D\big(\frac{a}{1+a}\big)\|_{W^{m-1,4}}\Big)(\|D^m\Phi\|_0+\|a\|_{m+1}\|D^{m+1}\theta\|_0)\\
&\lesssim(1+\|a\|_{m+1})^2\|\mathbf h\|_{m+1}\|D^{m+2}\theta,D^{m-1}\partial^2_2\mathbf h,D^m\Phi\|^2_0.
\end{align*}
In the same manner, one may deduce that
\begin{align*}
V_{4,2}\lesssim(1+\|a\|_{m+1})\|a\|_{m+1}\|\mathbf h\|_{m+1}\|D^{m+2}\theta,D^{m-1}\partial^2_2\mathbf h\|^2_0.
\end{align*}
Putting the above  estimates of $V_{4,1}$ and $V_{4,2}$ into  \eqref{10-16-Phi-3}, yields
$$V_4\lesssim(1+\|a\|_{m+1})^2\|\mathbf h\|_{m+1}\|D^{m+2}\theta,D^{m-1}\partial^2_2\mathbf h,D^m\Phi\|^2_0.$$
By the same argument as in the estimate for $V_4$, we  also have
\begin{align*}
V_5&=\int\partial^{\alpha-\mathbf e_k}\big((\partial_2\mathbf h\cdot\nabla)((1+a)^{-1})\big)\cdot\partial^\alpha(\Phi_k+{R}a\partial_k\theta) dx\\&\quad+{R}\int\partial^{\alpha-\mathbf e_k}\big(\frac{\partial_2\mathbf h}{1+a}\big)\cdot\partial^{\alpha+\mathbf e_k}(\nabla\theta a) dx\\
&\lesssim\big(\|D^{m-1}\partial_2\mathbf h\|_0\|Da\|_{L^\infty}+\|\partial_2\mathbf h\|_{L^4}\|\frac{a}{1+a}\|_{W^{m,4}}\big)(\|D^m\Phi\|_0+\|a\|_{m+1}\|D^{m+1}\theta\|_0)\\
&\quad+\big(\|D^{m-1}\partial_2\mathbf h\|+\|(1+a)^{-1}\|_{W^{m-1,4}}\|\partial_2\mathbf h\|_{L^4}\big)\|a\|_{m+1}\|D^{m+2}\theta\|_0\\
&\lesssim(1+\|a\|_{m+1})\|a\|_{m+1}\|D^{m+2}\theta,D^{m-1}\partial_2^2\mathbf h,D^m\Phi\|_0^2.
\end{align*}
At last, plugging the above  estimates of $V_j(j=1,2,3,4,5)$  into \eqref{10-16-Phi}, and summing over all multi-indices $\alpha$ with $|\alpha|=m$, we deduce the desired conclusion. This completes the proof of Lemma \ref{dissipation-Phi}.
\end{proof}
\subsection{The dissipation for $D\mathbf u$}
In this subsection, we shall construct the estimate of $\|D^{m+1}\mathbf u\|_{0}$ in order to close all the preceding estimates.
\begin{lemma}\label{dissipation-Du} Let $(a,\mathbf u,\theta,\mathbf h)$ be smooth solutions to  the system $\eqref{heat-mhd}$
  on $\mathbb T^2\times[0,T)$ for some $T>0$, and suppose that \eqref{density-1}, \eqref{theta-1} and \eqref{junzhi} hold.  Then,   for any $m\geq2$,
 \begin{align}
\frac{d}{dt}\|D^m\mathbf u\|_0^2+\mu\|D^{m+1}\mathbf u\|^2_0&\leq  C_{3}\|D^{m+2}\theta,D^{m-1}\partial^2_2\mathbf h,D^m\Phi\|_0^2\nonumber\\&\quad+C{(1+\|a\|_{m+1})\|a,\mathbf u,\mathbf h\|_{m+1}}\|D^{m+1}\mathbf u,D^{m-1}\partial_2^2\mathbf h,D^m\Phi\|^2_0,\label{10-16-1-Du}
\end{align}
where $C_3$ is  a positive constant  depending only on $|\mathbb T|,R,\mu,m$.
\end{lemma}
\begin{proof}
It follows from $\eqref{heat-mhd}_2$ that
\begin{equation*}\begin{split}
\partial_t\mathbf u-\mu\Delta\mathbf u-(\mu+\lambda)\nabla\mathrm{div}\mathbf u&=-\frac{\mu a\Delta\mathbf u+(\mu+\lambda)a\nabla\mathrm{div}\mathbf u}{1+a}-(\mathbf u\cdot\nabla)\mathbf u-R\nabla\theta\\&\quad-\frac{\Phi}{1+a}+\frac{(\mathbf h\cdot\nabla)\mathbf h}{1+a}+\frac{\partial_2\mathbf h}{1+a}.
\end{split}\end{equation*}
Given any multi-indices $\alpha$ with $|\alpha|=m$,  applying $\partial^{\alpha}$ to the above identity, and  then performing $L^2$ inner product
with $\partial^{\alpha}\mathbf u$ yield
\begin{align}
\frac{d}{2dt}\|\partial^\alpha\mathbf u\|_0^2+\mu\|D\partial^\alpha\mathbf u\|^2_0+(\mu+\lambda)\|\partial^\alpha\mathrm{div}\mathbf u\|_0^2\triangleq VI_1+VI_2+VI_3\label{10-16-Du}
\end{align}
with
\begin{align*}
&VI_1:=-\int\partial^\alpha\big(\frac{\mu a\Delta\mathbf u+(\mu+\lambda)a\nabla\mathrm{div}\mathbf u}{1+a}+(\mathbf u\cdot\nabla)\mathbf u\big)\cdot\partial^\alpha\mathbf udx,\end{align*}
\begin{align*}
&VI_2:=-\int\partial^\alpha\big(R\nabla\theta+\frac{\Phi}{1+a}-\frac{\partial_2\mathbf h}{1+a}\big)\cdot\partial^\alpha\mathbf udx,\end{align*}
\begin{align*}&VI_3:=\int\partial^\alpha\big(\frac{(\mathbf h\cdot\nabla)\mathbf h}{1+a}\big)\cdot\partial^\alpha\mathbf udx.
\end{align*}
In what follows, we shall estimate $VI_l(l=1,2,3)$  successively.  By integration by parts, H\"{o}lder's and Young's inequalities, Sobolev's embedding inequality,  Poincar\'{e}'s inequality, \eqref{density-1}, Lemmas \ref{Hm-estimate-1} and \ref{nolinear-estimate}, we conclude  that
\begin{align*}
VI_1&\lesssim\big(\|D^{m-1}\big(\frac{a\Delta\mathbf u}{1+a}\big)\|_0+\|D^{m-1}\big(\frac{a\nabla\mathrm{div}\mathbf u}{1+a}\big)\|_0+\|D^{m-1}\big((\mathbf u\cdot\nabla)\mathbf u\big)\|_0\big)\|D^{m+1}\mathbf u\|_0\\
&\lesssim(\|a\|_{L^\infty}\|D^{m+1}\mathbf u\|_0+\|\frac{a}{1+a}\|_{W^{m-1,4}}\|D^2\mathbf u\|_{L^4}+\|\mathbf u\|_{L^\infty}\|D^m\mathbf u\|_0+\|\mathbf u\|_{m-1}\|D\mathbf u\|_{L^\infty})\|D^{m+1}\mathbf u\|_0\\
&\lesssim\|a,\mathbf u\|_{m+1}\|D^{m+1}\mathbf u\|_0^2,
\end{align*}
and
\begin{align*}
VI_2&\lesssim (R\|D^{m}\theta\|_0+2\|D^{m-1}\Phi\|_0+2\|D^{m-1}\partial_2\mathbf h\|_0+C\|a\|_{m+1}\|D^{m-1}\partial_2^2\mathbf h,D^m\Phi\|_0)\|D\partial^\alpha\mathbf u\|_0\\
&\leq C_{3,1}\|D^{m+2}\theta,D^{m-1}\partial_2^{{2}}\mathbf h,D^m\Phi\|_0^2+\frac\mu2\|D\partial^\alpha\mathbf u\|_0^2+C\|a\|_{m+1}\|D^{m+1}\mathbf u,D^{m-1}\partial_2^2\mathbf h,D^m\Phi\|^2_0,
\end{align*}
where $C_{3,1}$ is a positive constant depending only on $|\mathbb T|,\mu,R$. At last,  let us turn to the term $VI_3$. It follows from integration by parts, H\"{o}lder, Young, Sobolev embedding, and Poincar\'{e} inequality, \eqref{density-1},  \eqref{2026-10-16-2}, Lemmas \ref{Hm-estimate-1} and  \ref{nolinear-estimate}, that
\begin{align*}
VI_3&\lesssim\big(\|(\mathbf h\cdot\nabla)\mathbf h\|_{m-1}+\|(1+a)^{-1}\|_{W^{m-1,4}}\|(\mathbf h\cdot\nabla)\mathbf h\|_{L^4}\big)\|D^{m+1}\mathbf u\|_0\\
&\lesssim(1+\|a\|_{m+1})\|\mathbf h\|_{m+1}\|D^{m+1}\mathbf u,D^{m-1}\partial^2_2\mathbf h\|^2_0.
\end{align*}
Consequently, substituting the above  estimates of  $VI_l(l=1,2,3)$ into \eqref{10-16-Du}, and summing over all multi-indices $\alpha$ with $|\alpha|=m$, we deduce the desired conclusion. This completes the proof of Lemma \ref{dissipation-Du}.
\end{proof}
\section{Proof of Theorem \ref{th:main1}}\label{Theorem-th:main2}

In this section, we shall  complete the proof of Theorem \ref{th:main1}. The main tool of the proof is the
bootstrapping argument. It starts with the ansatz that the smooth solution   $(a,\mathbf u,\theta,\mathbf h)$ to the system \eqref{heat-mhd}  satisfies
\begin{equation}\label{10-16-assume}
\sup_{t\in[0,T)}\|a,\mathbf u,\theta,\mathbf h\|^2_{N}+\int_0^T\|D\mathbf u,D\theta\|^2_{N}dt\leq\epsilon_0,
\end{equation}
for some $0<\epsilon_0\leq1$, where the small constant $\epsilon_0$ is to be determined in the next two lemmas and  depends only on the constants
$N,\mu,\lambda,R,c_\nu,\kappa,|\mathbb T|$ and $s$.

\subsection{Piecing together some  estimates} We now piece together several estimates from the preceding subsection to construct a key energy inequality.
\begin{lemma}\label{m-weifenbudengshi} %Let $M\geq2.$
Let $(a,\mathbf u,\theta,\mathbf h)$ be smooth solutions to  the system $\eqref{heat-mhd}$
on $\mathbb T^2\times[0,T)$ for some $T>0$, and suppose that \eqref{10-16-assume} holds. Then, for any $2\leq m\leq N-1$,
\begin{align}
\frac{d}{dt}E_m+2\mathcal D_m&\lesssim\|D\partial_2\mathbf h\|^2_0\|D^{m+2}\mathbf u\|^2_0+\big|\int{\frac{1+\theta}{1+a}}\partial_2\mathbf u_1\partial_1^{m+1}a\partial_1^m\partial_2adx\big|\nonumber\\
&\quad+\big|\int\partial_2\mathbf u_1\partial_1^{m+1}\mathbf h\partial_1^m\partial_2\mathbf hdx\big|\label{2026-10-17-2}
\end{align}
with
\begin{align*}
\mathcal D_m:=\|D^{m-1}\partial^2_2 a,D^{m+1}\partial_2\mathbf u,D^{m+1}\mathrm{div}\mathbf u,D^{m+1}\mathbf u,D^{m+1}\partial_2\theta,D^{m+2}\theta,D^{m-1}\partial^2_2\mathbf h,D^{m}\Phi\|_0^2,
\end{align*}
and
\begin{align}
E_m(t)\geq \frac{2}{\mu}\|D^m\partial_2a,D^m\partial_2\mathbf u,D^m\mathbf u,D^m\mathrm{div}\mathbf u, D^m\partial_2\mathbf h,D^m\partial_2\mathbf \theta,D^{m+1}\theta,D^{m}\Phi\|_0^2,\label{dengjia-2}
\end{align}
\begin{align}E_m(t)\leq Q_{N-1}\|D^m\partial_2a,D^m\partial_2\mathbf u,D^m\mathbf u,D^m\mathrm{div}\mathbf u, D^m\partial_2\mathbf h,D^m\partial_2\mathbf \theta,D^{m+1}\theta,D^{m}\Phi\|_0^2,\label{dengjia-1}
\end{align}
where $Q_{N-1}$ is a constant depending only on $|\mathbb T|, N ,\mu, \lambda, R, c_{\nu}, \kappa$. The explicit expression for $E_m$ is given {right after \eqref{2026-10-17-1}} below.
\end{lemma}
\begin{proof}
By Sobolev's embedding inequality, using \eqref{10-16-assume}, and choosing ${\epsilon_{0}}$ small enough, we have
\begin{align}\|a,\theta\|_{L^\infty}\lesssim\|a,\theta\|_{N}\leq\sqrt{{\epsilon_{0}}}\leq\frac12.\label{2026-10-17-3}\end{align}
%Thus, we can directly apply the conclusions of Lemma \ref{H-priori-v}, Lemma \ref{dissipation-rho}, Lemma \ref{dissipation-h}, Lemma \ref{dissipation-divu}, \ref{dissipation-Phi} and Lemma \ref{dissipation-Du}.
Let \(A>0\) be a sufficiently large constant, which will be fixed later. Multiplying \eqref{2-h-older} by $A^7$, \eqref{10-16-partial2a} by $4$, \eqref{10-16-partial2h} by $A^4$, \eqref{10-16-1-divu} by $A^3$, \eqref{10-16-1-Phi} by $A$, \eqref{10-16-1-Du} by $\frac{4}{\mu}$, and then adding the resulting equations together, we obtain
\begin{align}
&\frac{d}{dt}E_m+(2A^7\min\{\mu,\kappa\}-4C_0-A^4C_1-C_2A)\|D^{m+1}\partial_2\mathbf u,D^{m+1}\partial_2\theta\|_0^2+(A^4-\frac{4C_3}{\mu})\|D^{m-1}\partial^2_2\mathbf h\|_0^2\nonumber\\
&\quad+(2A^3\min\{(2\mu+\lambda),{\kappa}\}-C_2A-\frac{4C_3}{\mu})\|D^{m+1}\mathrm{div}\mathbf u,D^{m+2}\theta\|_0^2+(A-\frac{4C_3}{\mu})\|D^m\Phi\|_0^2\nonumber\\
&\quad+4\|D^{m-1}\partial^2_2a\|_0^2+4\|D^{m+1}\mathbf u\|^2_0-2A^3\|D^{m}\Phi\|_0\|D^{m+1}\partial_2\mathbf u\|_0\nonumber\\
&\lesssim(1+\|a,\mathbf u,\theta,\mathbf h\|_{m+1})^3\|a,\mathbf u,\theta,\mathbf h\|_{m+1}\|D^{m-1}(\partial_2^2a,D^2\partial_2\mathbf u,D^{2}\mathbf u,D^{2}\mathrm{div}\mathbf u,D^{3}\theta,\partial_2^2\mathbf h,D\Phi)\|_0^2\nonumber\\
&\quad+(1+\|\mathbf h\|_{m+1})\|D\partial_2\mathbf h\|_0\|D^{m+2}\mathbf u\|_0\|D^{m+1}\mathbf u,D^m\Phi\|_0\nonumber\\
&\quad+\big|\int{\frac{1+\theta}{1+a}}\partial_2\mathbf u_1\partial_1^{m+1}a\partial_1^m\partial_2adx\big|+\big|\int\partial_2\mathbf u_1\partial_1^{m+1}\mathbf h\partial_1^m\partial_2\mathbf hdx\big|,\label{2026-10-17-1}
\end{align}
where $C_0, C_1, C_2, C_3$ have already been defined in Lemmas \ref{dissipation-rho},  \ref{dissipation-h}, \ref{dissipation-Phi} and \ref{dissipation-Du}, and
\begin{align*}
E_m(t):&=A^7\|\sqrt{\frac{R(1+\theta)}{1+a}}D^m\partial_2a,\sqrt{\rho}D^m\partial_2\mathbf u,\sqrt{c_\nu\rho}D^m\partial_2\theta,D^m\partial_2\mathbf h\|^2_0+\frac{4}{\mu}\|D^m\mathbf u\|_0^2\\
&\quad+A^3\|(R+1)^{-\frac12}D^m\Phi,D^m\mathrm{div}\mathbf u,{\sqrt{c_\nu}}D^{m+1}\theta\|_0^2+2A^4\sum_{|\alpha|=m-1}\int\partial^\alpha\partial^2_2\mathbf u\cdot\partial^\alpha\partial_2\mathbf hdx\\
&\quad+2A\sum_{|\alpha|=m}\int\partial^\alpha\mathbf u\cdot\partial^\alpha\Phi dx+{\frac8R}\sum_{|\alpha|=m-1}\int\partial^\alpha\partial_2\mathbf u_2\partial^\alpha\partial^2_2adx.
\end{align*}
Due to $$2A^3\|D^{m}\Phi\|_0\|D^{m+1}\partial_2\mathbf u\|_0\leq\|D^{m}\Phi\|^2_0+A^6\|D^{m+1}\partial_2\mathbf u\|^2_0,$$
we take A sufficiently large such that
\begin{equation}\label{Adejiashe-1}\left\{\begin{split}&2A^7\min\{\mu,\kappa\}-4C_0-A^4C_1-C_2A-A^6>4, \quad A-\frac{4C_3}{\mu}-1>4, \\
&{A^4-\frac{4C_3}{\mu}>4,}\quad 2A^3\min\{(2\mu+\lambda),{\kappa}\}-C_2A-\frac{4C_3}{\mu}>4.
\end{split}\right.
\end{equation}
From which, \eqref{10-16-assume} and \eqref{2026-10-17-1} lead to
\begin{align}
&\frac{d}{dt}E_m+4\|D^{m-1}\partial^2_2 a,D^{m+1}\partial_2\mathbf u,D^{m+1}\mathrm{div}\mathbf u,D^{m+1}\partial_2\theta,D^{m+2}\theta,D^{m-1}\partial^2_2\mathbf h,D^{m+1}\mathbf u,D^{m}\Phi\|_0^2\nonumber\\
&\lesssim\sqrt{{\epsilon_{0}}}\|D^{m-1}(\partial_2^2a,D^2\partial_2\mathbf u,D^{2}\mathbf u,D^{2}\mathrm{div}\mathbf u,D^{3}\theta,\partial_2^2\mathbf h,D\Phi)\|_0^2+\|D\partial_2\mathbf h\|_0\|D^{m+2}\mathbf u\|_0\|D^{m+1}\mathbf u,D^m\Phi\|_0\nonumber\\
&\quad+\big|\int{\frac{1+\theta}{1+a}}\partial_2\mathbf u_1\partial_1^{m+1}a\partial_1^m\partial_2adx\big|+\big|\int\partial_2\mathbf u_1\partial_1^{m+1}\mathbf h\partial_1^m\partial_2\mathbf hdx\big|.\label{Em-biaoda}
\end{align}
Then employing  Young's inequality and taking ${\epsilon_{0}}$ sufficiently small, we conclude that  \eqref{2026-10-17-2} holds.
It follows from H\"{o}lder's and  Young's inequalities,    Poincar\'{e}'s inequality, that
\begin{align*}
&2A^4\Big|\sum_{|\alpha|=m-1}\int\partial^\alpha\partial^2_2\mathbf u\cdot\partial^\alpha\partial_2\mathbf hdx\Big|\leq C_4A^4\|D^{m}\partial_2\mathbf u,D^{m}\partial_2\mathbf h\|_0^2,\end{align*}
\begin{align*}&2A\Big|\sum_{|\alpha|=m}\int\partial^\alpha\mathbf u\cdot\partial^\alpha\Phi dx\Big|\leq \frac{2}{\mu}\|D^m\mathbf u\|^2+C_4A^2\|D^m\Phi\|_0^2,\end{align*}
\begin{align*}&{\frac8R}\Big|\sum_{|\alpha|=m-1}\int\partial^\alpha\partial_2\mathbf u_2\partial^\alpha\partial^2_2adx\Big|\leq C_4\|D^m\partial_2a,D^{m}\partial_2\mathbf u\|_0^2,
\end{align*}
where $C_4$ is a positive constant  depending only on $m,|\mathbb T|{,\mu,R}.$ This together with {the definition of $E_m$} and \eqref{2026-10-17-3} yields
 \begin{align*}
E_m&\geq(A^7-C_4A^4-C_4)\min\{\frac R3,\frac12,\frac{c_\nu}{2}\}\|D^m\partial_2a,D^m\partial_2\mathbf u,D^m\partial_2\theta,D^m\partial_2\mathbf h\|^2_0\\
&\quad+\frac{2}{\mu}\|D^m\mathbf u\|_0^2+(A^3-C_4A^2)\min\{\frac{1}{R+1},{c_\nu}\}\|D^m\Phi,D^m\mathrm{div}\mathbf u,D^{m+1}\theta\|_0^2,
\end{align*}
\begin{align*}E_m&\leq(A^7+C_4A^4+C_4)\max\{3R,\frac32,\frac{3c_\nu}{2}\}\|D^m\partial_2a,D^m\partial_2\mathbf u,D^m\partial_2\theta,D^m\partial_2\mathbf h\|^2_0\\
&\quad+\frac{6}{\mu}\|D^m\mathbf u\|_0^2+(A^3+C_4A^2)\max\{1,{c_\nu}\}\|D^m\Phi,D^m\mathrm{div}\mathbf u,D^{m+1}\theta\|_0^2.
\end{align*}
Thus, we may again choose $A$ large enough satisfying \eqref{Adejiashe-1} and
\begin{equation*}(A^7-C_4A^4-C_4)\min\{\frac R3,\frac12,\frac{c_\nu}{2}\}>\frac{2}{\mu},\ (A^3-C_4A^2)\min\{\frac{1}{R+1},{c_\nu}\}>\frac{2}{\mu},
\end{equation*} 
which ensures that
\eqref{dengjia-2} and
$$E_m\leq \tilde{Q}_m\|D^m\partial_2a,D^m\partial_2\mathbf u,D^m\mathbf u,D^m\mathrm{div}\mathbf u, D^m\partial_2\mathbf h,D^m\partial_2\mathbf \theta,D^{m+1}\theta,D^{m}\Phi\|_0^2$$ hold,  
where $\tilde{Q}_m$ is a constant depending only on $|\mathbb T|,m,\mu,\lambda,R,c_{\nu},\kappa$. Thus, \eqref{dengjia-1} follows by denoting
$Q_{N-1}:=\max\{\tilde{Q}_2,\tilde{Q}_3,\cdots,\tilde{Q}_{N-1}\}.$ This completes the proof of Lemma \ref{m-weifenbudengshi}.
\end{proof}

\subsection{Decay estimates of the higher-order energy} Building upon the dissipation estimates from the preceding subsection,  we now  establish decay estimates for $\partial_2a,D\mathbf u,\mathrm{div}\mathbf u, D\theta,\partial_2\mathbf h$ and $\Phi$.

\begin{lemma}\label{decay} %Let $M\geq4$ in Lemma \ref{m-weifenbudengshi}.
{Let $N\geq5$ and  the assumptions of Lemma \ref{m-weifenbudengshi} hold.} Given an arbitrarily fixed constant $s\in(0,1)$, then, for any $l=0,1,2\cdots,N-3,$
\begin{align}
&\sup_{t\in[0,T)}(1+t)^{s+l}E_{N-1-l}+\int_0^T(1+t)^{s+l}\mathcal D_{N-1-l}dt\lesssim E_{N-1}(0)+\sup_{t\in[0,T)}\|a,\mathbf u,\theta,\mathbf h\|^2_{N}.\label{2026-10-18-9}
\end{align}
\end{lemma}
\begin{proof} The proof is divided into two steps.

\textbf{Step 1.} We first derive the decay estimate for $E_2(t)$, which serves as a key ingredient for the second step of the proof.
For any $2\leq m\leq N-2$, it follows from integration by parts and    H\"{o}lder's, Young's and Poincar\'{e}'s inequalities that
\begin{align*}
&\big|\int{\frac{1+\theta}{1+a}}\partial_2\mathbf u_1\partial_1^{m+1}a\partial_1^m\partial_2adx\big|\\
&{\leq}\big|\int{\frac{1+\theta}{1+a}}\partial_2\mathbf u_1\partial_1^{m+2}a\partial_1^{m-1}\partial_2adx\big|+\big|\int{\partial_1\Big(\frac{1+\theta}{1+a}\partial_2\mathbf u_1\Big)}\partial_1^{m+1}a\partial_1^{m-1}\partial_2adx\big|\\
&\lesssim{(1+\|D(a,\theta)\|_{L^\infty})}\|\partial_2\mathbf u,D\partial_2\mathbf u\|_{L^\infty}\|D^{m+1}a\|_{1}\|D^{m-1}\partial_2a\|_{{0}}\\&\lesssim\|a\|_{N}\|D^{m-1}\partial^2_2a,D^{m+1}\partial_2\mathbf u\|_0^2{,}
\end{align*}
{where we have used \eqref{10-16-assume},  \eqref{2026-10-17-3} and $m+2\leq N$.}
Similarly, we also have
$$\big|\int\partial_2\mathbf u_1\partial_1^{m+1}\mathbf h\partial_1^m\partial_2\mathbf hdx\big|\lesssim\|\mathbf h\|_{N}\|D^{m+1}\partial_2\mathbf u,D^{m-1}\partial^2_2\mathbf h\|_0^2.$$
Combining the above two estimates with \eqref{2026-10-17-2}, and using Poincar\'{e}'s inequality, the condition \eqref{10-16-assume} and the definition of $\mathcal D_m$, we have
\begin{align*}
\frac{d}{dt}E_m+2\mathcal D_m&\lesssim{\epsilon_{0}}\|D\partial_2\mathbf h\|^2_0+\sqrt{\epsilon_{0}}\|D^{m-1}\partial^2_2a,D^{m+1}\partial_2\mathbf u,D^{m-1}\partial^2_2\mathbf h\|_0^2\lesssim{\epsilon_{0}}\mathcal D_m+\sqrt{\epsilon_{0}}\mathcal D_m,
\end{align*}
where we have used the following fact
$$\|D\partial_2\mathbf h\|_0\lesssim\|D\partial_{22}\mathbf h\|_0,$$ due to $\int_{\mathbb T} D\partial_2 \mathbf h dx_2=0.$
Thus, choosing ${\epsilon_{0}}$ sufficiently small again, we conclude that for any $2\leq m\leq N-2,$
\begin{equation}\label{M-1weifen}
\frac{d}{dt}E_m+\mathcal D_m\leq0.
\end{equation}
It follows from Gagliardo-Nirenberg interpolation inequality, \eqref{dengjia-1} and Poincar\'{e}'s inequality that,  for any $2\leq m\leq N-2,$
\begin{align}
E_{m}&\lesssim\big(\|D^{m-1}(\partial_2a,\partial_2\mathbf u,\mathbf u,\mathrm{div}\mathbf u, \partial_2\mathbf h,\partial_2\mathbf \theta,D\theta,D\Phi)\|_0^2\big)^{\frac{N-m-1}{N-m}}\nonumber\\
&\quad\times\big(\|D^{{N-1}}(\partial_2a,\partial_2\mathbf u,\mathbf u,\mathrm{div}\mathbf u, \partial_2\mathbf h,\partial_2\mathbf \theta,D\theta,{\Phi})\|_0^2\big)^{\frac{1}{N-m}}\nonumber\\
&\lesssim(\mathcal D_m)^{\frac{N-m-1}{N-m}}\|a,\mathbf u,\theta,\mathbf h\|_{N}^{\frac{2}{N-m}}\nonumber\\&\leq C{\epsilon_{0}}^{\frac{1}{N-m}}(\mathcal D_{m})^{\frac{N-m-1}{N-m}},\label{2026-10-17-4}
\end{align}
where %by the definition of $\Phi$, Lemma \ref{Hm-estimate-1}, Sobolev embedding inequality and \eqref{10-16-assume},
we have used the following fact
\begin{align}
\|\Phi\|_{N-1}&\lesssim\|R\nabla a(1+\theta)\|_{N-1}+\||\mathbf h|^2\|_{N}+\|\mathbf h_2\|_{N}\nonumber\\
&\lesssim(1+\|\mathbf h,\theta\|_{N})\|a,\theta,\mathbf h\|_{N}\nonumber\\&\lesssim\|a,\theta,\mathbf h\|_{N}.\label{2026-10-17-6}
\end{align}
Combining \eqref{2026-10-17-4} with \eqref{M-1weifen}, and choosing ${\epsilon_{0}}$ small enough such that $C{\epsilon_{0}}^{{\frac{1}{N-m}}}\leq1$ in \eqref{2026-10-17-4}, yield that for any $2\leq m\leq N-2,$
\begin{equation*}
\frac{d}{dt}E_m+E_m^{\frac{N-m}{N-m-1}}\leq0.
\end{equation*}
This leads to
\begin{equation}\label{2026-10-17-5}
E_m(t)\leq E_m(0)\big(\alpha t+1\big)^{-(N-m-1)}
\end{equation}
with $\alpha:=\frac{E_m^{\frac{1}{N-m-1}}(0)}{N-m-1}.$
Given any $s_0\in(0,1)$, we conclude  from \eqref{M-1weifen} and \eqref{2026-10-17-5} that
\begin{align}
\sup_{t\in[0,T)}&(1+t)^{1+s_0}E_2(t)+\int_0^T(1+t)^{1+s_0}\mathcal D_2dt\nonumber\\&\leq(1+s_0)\int_0^T(1+t)^{s_0}E_2(t)dt+E_2(0)\nonumber\\
&\leq (1+s_0)E_2(0)\sup_{t\in[0,\infty)}\big(\frac{1+t}{\alpha t+1}\big)^{s_0}\int_0^\infty(\alpha t+1)^{-(N-3-s_0)}dt+E_2(0)\nonumber\\
&\lesssim E_2(0)(1+\alpha^{-s_0})\alpha^{-1}+E_2(0)\nonumber\\
&\lesssim E_2^{\frac{N-4-s_0}{N-3}}(0)+E_2^{\frac{N-4}{N-3}}(0)+E_2(0)\lesssim1,\label{2026-10-17-8}
\end{align}
where we have used the facts $\int_0^\infty(1+y)^{-(N-3-s_0)}dy<\infty$,  due to $N\geq5${, the definition of $\alpha$ with $m=2$ (we may assume $E_2(0)>0$, since otherwise $E_2\equiv0$ by \eqref{M-1weifen}), and $E_2(0)\lesssim\|a_0,\mathbf u_0,\theta_0,\mathbf h_0\|_N^2\leq\epsilon_0\leq1$, which follows from \eqref{dengjia-1}, \eqref{2026-10-17-6} and \eqref{10-16-assume}. In particular, the estimate of \eqref{2026-10-17-8} is independent of $T$ and of the solution; this uniformity is used in Step 2 below}. %\eqref{10-16-assume}, and \eqref{2026-10-17-6},
%and $E_2(0)\lesssim{\epsilon_{0}},$ in the above inequality.

\textbf{Step 2. }It follows from \eqref{2026-10-17-2}, H\"{o}lder's inequality, Sobolev's embedding and Poincar\'{e}'s inequalities, and the definition of $\mathcal D_{2}$ that
\begin{align}
\frac{d}{dt}E_{N-1}+2\mathcal D_{N-1}&\lesssim\|D\partial_2\mathbf h\|^2_0\|D^{N+1}\mathbf u\|^2_0+\big|\int{\frac{1+\theta}{1+a}}\partial_2\mathbf u_1\partial_1^{N}a\partial_1^{N-1}\partial_2adx\big|\nonumber\\
&\quad+\big|\int\partial_2\mathbf u_1\partial_1^{N}\mathbf h\partial_1^{N-1}\partial_2\mathbf hdx\big|\nonumber\\
&\lesssim\|D\partial_2\mathbf h\|^2_0\|D^{N+1}\mathbf u\|^2_0+\|\partial_2\mathbf u_1\|_{L^{\infty}}\|a,\mathbf h\|_{N}\|D^{N-1}\partial_2a,D^{N-1}\partial_2\mathbf h\|_0\nonumber\\
&\lesssim\|D\partial_2\mathbf h\|^2_0\|D^{N+1}\mathbf u\|^2_0+\|\partial_2\mathbf u_1\|_{2}\|a,\mathbf h\|_{N}\|D^{N-1}\partial_2a,D^{N-1}\partial_2\mathbf h\|_0\nonumber\\&\lesssim\|D\partial_2\mathbf h\|^2_0\|D^{N+1}\mathbf u\|^2_0+\|\partial_2D^2\mathbf u\|_0\|a,\mathbf h\|_{N}\|D^{N-1}\partial_2a,D^{N-1}\partial_2\mathbf h\|_0\nonumber\\
&\lesssim E_{N-1}\|D\mathbf u\|^2_{N}+(\mathcal D_2)^{\frac12}\|a,\mathbf h\|_{N}(E_{N-1})^\frac{1}{2}.\label{2026-10-17-7}
\end{align}
By the fact $\|D^{N-1}\partial_2\mathbf h\|^2_0\leq\|D^{N-2}\partial_2^2\mathbf h\|_0\|D^{N}\mathbf h\|_0$, Poincar\'{e}'s inequality, and the estimates \eqref{dengjia-2}, \eqref{dengjia-1} and \eqref{2026-10-17-6}, we conclude
$$E_{N-1}\lesssim\|a,\mathbf u,\theta,\mathbf h\|_{N}(\mathcal D_{N-1})^\frac12.$$
This along with \eqref{2026-10-17-8}, \eqref{2026-10-17-7},  \eqref{10-16-assume}, and H\"{o}lder's and Young's inequalities,  yields for any $s\in(0,1)$.
\begin{align*}
&\sup_{t\in[0,T)}(1+t)^sE_{N-1}+2\int_0^T(1+t)^s\mathcal D_{N-1}dt\\&\leq E_{N-1}(0)+s\int_0^T(1+t)^{s-1}E_{N-1}dt+C\sup_{t\in[0,{T})}(1+t)^s E_{N-1}(t)\int_0^T\|D\mathbf u\|^2_{N}dt\\
&\quad+C\sup_{t\in[0,T)}\|a,\mathbf h\|_{N}\int_0^T\big((1+t)^s\mathcal D_2\big)^{\frac12}\big((1+t)^{s}E_{N-1}\big)^\frac{1}{2}dt\\
&\leq E_{N-1}(0)+C\sup_{t\in[0,T)}\|a,\mathbf u,\theta,\mathbf h\|_{N}\int_0^T(1+t)^{\frac{s-2}{2}}\big((1+t)^s\mathcal D_{N-1}\big)^\frac{1}{2}dt+C{\epsilon_{0}}\sup_{t\in[0,T)}(1+t)^s E_{N-1}(t)\\
&\quad+C_\delta\sup_{t\in[0,T)}\|a,\mathbf h\|^2_{N}\Big(\int_0^T\big((1+t)^s\mathcal D_2\big)^{\frac12}dt\Big)^2+\delta\sup_{t\in[0,T)}(1+t)^{s}E_{N-1}(t)\\
&\leq E_{N-1}(0)+C\int_0^T(1+t)^{s-2}dt\sup_{t\in[0,T)}\|a,\mathbf u,\theta,\mathbf h\|^2_{N}+\int_0^T(1+t)^s\mathcal D_{N-1}dt\\
&\quad+C_\delta\sup_{t\in[0,T)}\|a,\mathbf h\|^2_{N}+(C{\epsilon_{0}}+\delta)\sup_{t\in[0,T)}(1+t)^s E_{N-1}(t),
\end{align*}
where we have used the following inequality
\begin{align*}
\int_0^T\big((1+t)^s\mathcal D_2\big)^{\frac12}dt\leq\Big(\int_0^T(1+t)^{\frac{s-3}{2}}dt\Big)^\frac12\Big(\int_0^T(1+t)^{{1+s_0}}\mathcal D_2dt\Big)^\frac{1}{2}{\lesssim1}
\end{align*}
with $s_0:=\frac{s+1}2$ in \eqref{2026-10-17-8}.
By choosing ${\epsilon_{0}}$ and $\delta$ small enough, we have
\begin{align}\label{3.63AA}
&\sup_{t\in[0,T)}(1+t)^sE_{N-1}+2\int_0^T(1+t)^s\mathcal D_{N-1}dt\leq 2E_{N-1}(0)+C\sup_{t\in[0,T)}\|a,\mathbf u,\theta,\mathbf h\|^2_{N}.
\end{align}
Based on the above inequality and the fact $E_{N-2}\lesssim\mathcal D_{N-1}$ and $E_{N-2}\lesssim E_{N-1}$, following a similar derivation in \eqref{3.63AA},  we deduce from \eqref{M-1weifen} that
\begin{align*} 
\sup_{t\in[0,T)}&(1+t)^{s+1}E_{N-2}+\int_0^T(1+t)^{s+1}\mathcal D_{N-2}dt\\&\leq E_{N-2}(0)+(s+1)\int_0^T(1+t)^{s}E_{N-2}dt\\
&\leq E_{N-2}(0)+C\int_0^T(1+t)^{s}\mathcal D_{N-1}dt\\&\lesssim E_{N-1}(0)+\sup_{t\in[0,T)}\|a,\mathbf u,\theta,\mathbf h\|^2_{N}.
\end{align*}
Similarly, we also have
\begin{align*}
\sup_{t\in[0,T)}&(1+t)^{s+2}E_{N-3}+\int_0^T(1+t)^{s+2}\mathcal D_{N-3}dt\\&\leq E_{N-3}(0)+(s+2)\int_0^T(1+t)^{s+1}E_{N-3}dt\\
&\leq E_{N-3}(0)+C\int_0^T(1+t)^{s+1}\mathcal D_{N-2}dt\\&\lesssim E_{N-1}(0)+\sup_{t\in[0,T)}\|a,\mathbf u,\theta,\mathbf h\|^2_{N}.
\end{align*}
Thus, in general, it holds that for any $l=0,1,2\cdots, N-3$,
\begin{align*}
&\sup_{t\in[0,T)}(1+t)^{s+l}E_{N-l-1}+\int_0^T(1+t)^{s+l}\mathcal D_{N-l-1}dt\lesssim E_{N-1}(0)+\sup_{t\in[0,T)}\|a,\mathbf u,\theta,\mathbf h\|^2_{N}.
\end{align*}
This completes the proof of Lemma \ref{decay}.\end{proof}
\subsection{Uniform boundedness estimates of the highest-order energy} Based on Lemmas \ref{priori-basic},  \ref{H-priori-v},  \ref{m-weifenbudengshi} and  \ref{decay},  we here successfully construct  the  uniform boundedness estimates for smooth solutions $(a,\mathbf u,\theta,\mathbf h)$ to  the system $\eqref{heat-mhd}$.
\begin{proposition}\label{Uniform-boundedness}%Let $M\geq4.$
{Let $N\geq5$.} Assume the initial data $(a_0,\mathbf u_0,\theta_0,\mathbf h_0)\in H^{N}(\mathbb T^2)$ satisfying \eqref{a-theta-bound} and \eqref{junzhi}, and smooth solutions $(a,\mathbf u,\theta,\mathbf h)$ to  the system $\eqref{heat-mhd}$ satisfying  $(a,\mathbf u,\theta,\mathbf h)\in C([0,\tilde{T});H^N(\mathbb T^2))$ and $(D\mathbf u,D\theta)\in L^2([0,\tilde{T});H^N(\mathbb T^2))$ for some $\tilde{T}>0$. Given an arbitrarily fixed constant $s\in(0,1)$, if
\begin{equation}\label{2026-10-18-5}\|a_0,\mathbf u_0,\theta_0,\mathbf h_0\|_{N}\leq \epsilon,\end{equation}
then there is a positive constants $\epsilon$ depending only on $N,\mu,\lambda,R,c_\nu,\kappa,|\mathbb T|$ and $s$ such that
\begin{equation}\label{2026-10-18-11}
\sup_{t\in[0,\tilde{T})}\|a,\mathbf u,\theta,\mathbf h\|_{N}\lesssim\epsilon,
\end{equation}
and
\begin{align}
&\sup_{t\in[0,\tilde T)}(1+t)^{s+l}\|\partial_2a,\partial_2\mathbf u,D\theta,\partial_2\mathbf h,\mathrm{div}\mathbf u,\Phi\|_{N-l-1}^2\nonumber\\
&\quad+\int_0^{\tilde T}(1+t)^{s+l}\big(\|\partial_2\mathbf u,\mathrm{div}\mathbf u,D\theta\|^2_{N-l}+\|D\mathbf u,\Phi\|^2_{N-l-1}+\|\partial_2^2a,\partial_2^2\mathbf h\|^2_{N-l-2}\big)dt\lesssim\epsilon^2,\label{2026-10-18-12}
\end{align}
for $l=0,1,2,\cdots, N-3.$
\end{proposition}
\begin{proof}
Under the regularity assumptions of $(a,\mathbf u,\theta,\mathbf h)$, the arguments in the proofs of the preceding lemmas remain justified, and the same a priori estimates can be derived for the corresponding local solution. By \eqref{2026-10-18-5}, choosing $\epsilon<1$ ensures that the following set is nonempty,
$$\mathcal A:=\{T\in(0,\tilde{T});\sup_{t\in[0,T)}\|a,\mathbf u,\theta,\mathbf h\|^2_{N}+\int_0^T\|D\mathbf u,D\theta\|^2_{N}dt\leq\epsilon\}.$$
Let $$T_0:=\sup_{T\in \mathcal A} T.$$ Then $T_0\in(0,\tilde{T}].$ Now we show that, by a suitable choice of $\epsilon$, one can ensure that $T_0=\tilde{T}.$
We choose $\epsilon\leq \epsilon_0$, where $\epsilon_0$ has been defined in Lemma \ref{m-weifenbudengshi},   and then conclude from \eqref{2026-10-18-9} and \eqref{2026-10-17-7} that for $l=0,1,2,\cdots, N-3,$
\begin{align}
&\sup_{t\in[0,T_0)}(1+t)^{s+l}E_{N-l-1}+\int_0^{T_0}(1+t)^{s+l}\mathcal D_{N-l-1}dt\lesssim E_{N-1}(0)+\sup_{t\in[0,T_0)}\|a,\mathbf u,\theta,\mathbf h\|^2_{N},\label{2026-10-18-10}
\end{align}
and for any $t\in[0,T_0),$
\begin{align}
\frac{d}{dt}E_{N-1}+2\mathcal D_{N-1}&\lesssim E_{N-1}\|D\mathbf u\|^2_{N}+(\mathcal D_2)^{\frac12}\|a,\mathbf h\|_{N}(E_{N-1})^\frac{1}{2}\nonumber\\
&\lesssim\sqrt{\epsilon}\|D\mathbf u\|^2_{N}+(\mathcal D_2)^{\frac12}\|a,\mathbf u,\theta,\mathbf h\|^2_{N}.\label{2026-10-17-9}
\end{align}
Here, we have used the the following fact
\begin{equation}\label{2026-10-18-6}E_{N-1}\lesssim\|a,\mathbf u,\theta,\mathbf h\|^2_{N}.\end{equation}
By virtue of \eqref{2026-10-17-3}, we have
\begin{equation}\label{2026-10-18-2}|a|,|\theta|\leq\frac12.\end{equation}
Then, it follows from \eqref{basic-1}, \eqref{1-h-older}, Poincar\'{e}'s inequality and the definition of $T_0$ that, for any $t\in[0,T_0),$
\begin{align}
&\frac{d}{dt}\|\sqrt{R}a,\sqrt{\rho}\mathbf u,\sqrt{c_\nu\rho}\theta,\mathbf h\|_0^2+2\min\{\mu,\kappa\}\|D\mathbf u,D\theta\|_0^2\lesssim(\mathcal D_{N-2})^\frac12\|a,\theta,\mathbf u\|^2_{N},\label{2026-10-18}
\end{align}
and
\begin{align}\frac{d}{dt}&\|\sqrt{\frac{R(1+\theta)}{1+a}}D^{N-1}\partial_1a,\sqrt{\rho}D^{N-1}\partial_1\mathbf u,\sqrt{c_\nu\rho}D^{N-1}\partial_1\theta,D^{N-1}\partial_1\mathbf h\|^2_0\nonumber\\&\quad+2\min\{\mu,\kappa\}\|D^{N}\partial_1\mathbf u,D^{N}\partial_1\theta\|_0^2\nonumber\\
&\lesssim(\mathcal D_{N-2})^\frac12\|a,\mathbf u,\mathbf h\|_{N}^2+\sqrt{\epsilon}\|D^{N+1}\mathbf u,D^{N+1}\theta,{D^{N-2}\partial_2^2a,}D^{N-2}\partial_2^2\mathbf h\|_0^2.\label{2026-10-19}
\end{align}
Adding \eqref{2026-10-17-9}, \eqref{2026-10-18} and \eqref{2026-10-19} together{, noting that $\mathcal D_{N-1}$ contains $\|D^{N}\partial_2\mathbf u\|_0^2$, $\|D^{N+1}\theta\|_0^2$ and $\|D^{N-2}\partial_2^2(a,\mathbf h)\|_0^2$, and that $\mathcal D_2\lesssim\mathcal D_{N-2}$ by Poincar\'{e}'s inequality,} and then   taking $\epsilon$ sufficiently small, we deduce that for any $t\in[0,T_0)$,
\begin{align}
\frac{d}{dt}\mathcal E+\min\{\mu,\kappa{,1}\}\|D\mathbf u,D\theta,D^{N+1}\mathbf u,D^{N+1}\theta\|_0^2+\mathcal D_{N-1}\lesssim(\mathcal D_{N-2})^\frac12\|a,\mathbf u,{\theta,}\mathbf h\|_{N}^2,\label{2026-10-18-1}
\end{align}
where
\begin{align*}\mathcal E(t)&:=E_{N-1}(t)+\|\sqrt{R}a,\sqrt{\rho}\mathbf u,\sqrt{c_\nu\rho}\theta,\mathbf h\|_0^2\\
&\quad+\|\sqrt{\frac{R(1+\theta)}{1+a}}D^{N-1}\partial_1a,\sqrt{\rho}D^{N-1}\partial_1\mathbf u,\sqrt{c_\nu\rho}D^{N-1}\partial_1\theta,D^{N-1}\partial_1\mathbf h\|^2_0.\end{align*}
By \eqref{2026-10-18-6} and \eqref{2026-10-18-2}, we conclude that 
\begin{equation}\label{2026-10-18-3}\|a,\mathbf u,\theta,\mathbf h\|^2_{N}\lesssim\mathcal E\lesssim\|a,\mathbf u,\theta,\mathbf h\|^2_{N}.\end{equation}
This together with \eqref{2026-10-18-1}, yields for any $t\in[0,T_0)$,
\begin{align}\label{2026-10-18-7}
\frac{d}{dt}\mathcal E+c_0\|D\mathbf u,D\theta\|^2_{N}\lesssim(\mathcal D_{N-2})^\frac12\|a,\mathbf u,{\theta,}\mathbf h\|_{N}^2\leq C(\mathcal D_{N-2})^\frac12\mathcal E.
\end{align}
Here, we have used the fact $$c_0\|D(\cdot)\|^2_{N}\leq\min\{\mu,\kappa{,1}\}\|D(\cdot),D^{N+1}(\cdot)\|_0^2,$$
where $c_0$  is a positive constant depending only on $\mu,\kappa,{|\mathbb T|},N$. We deduce  by \eqref{2026-10-18-7}, Gronwall's and H\"{o}lder's inequalities, the definition of $T_0$, \eqref{2026-10-18-10} and \eqref{2026-10-18-6}, that
\begin{align*}
\sup_{t\in[0,T_0)}\mathcal E(t)+c_0\int_0^{T_0}\|D\mathbf u,D\theta\|^2_{N}dt&\leq e^{C\int_0^{T_0}(\mathcal D_{N-2})^\frac12dt}\mathcal E(0)\\
&\leq e^{C\big(\int_0^\infty(1+t)^{-s-1}dt\big)^\frac12\big(\int_0^{T_0}(1+t)^{s+1}\mathcal D_{N-2}dt\big)^\frac12}\mathcal E(0)\\
&\leq e^{C\sqrt{E_{N-1}(0)}+C\sup_{t\in[0,T_0)}\|a,\mathbf u,\theta,\mathbf h\|_{N}}\mathcal E(0)\\&\leq e^{C\sqrt{\epsilon}}\mathcal E(0).
\end{align*}
This together with \eqref{2026-10-18-5} and \eqref{2026-10-18-3}, yields
\begin{align}
\sup_{t\in[0,T_0)}\|a,\mathbf u,\theta,\mathbf h\|^2_{N}+\int_0^{T_0}\|D\mathbf u,D\theta\|^2_{N}dt&\lesssim\sup_{t\in[0,T_0)}\mathcal E(t)+c_0\int_0^{T_0}\|D\mathbf u,D\theta\|^2_{N}dt\nonumber\\
&\lesssim\mathcal E(0)\lesssim\|a_0,\mathbf u_0,\theta_0,\mathbf h_0\|^2_{N}\lesssim\epsilon^2.\label{2026-10-18-4}
\end{align}
Choosing $\epsilon$ small enough, leads to
\begin{equation}\label{2026-10-18-8}\sup_{t\in[0,T_0)}\|a,\mathbf u,\theta,\mathbf h\|^2_{N}+\int_0^{T_0}\|D\mathbf u,D\theta\|^2_{N}dt\leq\frac{\epsilon}{2}.\end{equation}
And then $T_0=\tilde T$. Indeed, if $T_0<\tilde{T}$, then it follows from \eqref{2026-10-18-8} and the regularity assumptions of $(a,\mathbf u,\theta,\mathbf h)$ that one can extend to another time $T_1\in(T_0,\tilde T)$ such that
$$\sup_{t\in[0,T_1)}\|a,\mathbf u,\theta,\mathbf h\|^2_{N}+\int_0^{T_1}\|D\mathbf u,D\theta\|^2_{N}dt\leq\epsilon.$$ This contradicts the definition of $T_0.$ Thus, it must be $T_0=\tilde T.$ \eqref{2026-10-18-11} follows by \eqref{2026-10-18-4}. \eqref{2026-10-18-12} follows by \eqref{dengjia-2}, \eqref{2026-10-18-10}, \eqref{2026-10-18-6}, \eqref{2026-10-18-5}, \eqref{2026-10-18-11}, the definition of $\mathcal D_{N-l-1}$, Poincar\'{e}'s inequality and \eqref{2026-10-18-13}. This completes the proof of Proposition \ref{Uniform-boundedness}.\end{proof}

With Proposition \ref{Uniform-boundedness} in hand, we now  complete the proof of
Theorem \ref{th:main1}

\textbf{The proof of Theorem \ref{th:main1}.} The local well-posedness of system \eqref{heat-mhd} can be established via standard energy methods, based on the solvability of the full compressible Navier-Stokes equations and transport equations. See \cite{KawashimaS} for a similar result. Thus, given initial data $(a_0,\mathbf u_0,\theta_0,\mathbf h_0)\in
H^N(\mathbb T^2)$ satisfying \eqref{a-theta-bound}, \eqref{junzhi}, we may assume that there exists $T>0$ such that the system \eqref{heat-mhd} has a unique
solution $(a,\mathbf u,\theta,\mathbf h)\in C([0,T);H^N(\mathbb T^2))$ with $(D\mathbf u,D\theta)\in L^2([0,T];H^N(\mathbb T^2))$. By iteratively applying the local well-posedness result, one can extend the local solution to the maximal time of existence $T_{\mathrm{max}}.$ Now, under the condition \eqref{initial-date-1}, we claim that
$T_\text{max}=\infty$ and the conclusion follows by Proposition \ref{Uniform-boundedness} . Otherwise, if $T_{\mathrm{max}}<\infty$, then it holds
\begin{equation}\label{2026-5-9-10}
\varlimsup_{t\rightarrow T_\text{max}}\|a,\mathbf u,\theta,\mathbf h\|_{N}=\infty.
\end{equation}
However,  it follows from  Proposition \ref{Uniform-boundedness} that
\begin{align*}
&\sup_{t\in[0,T_{\mathrm{max}})}\|a,\mathbf u,\theta,\mathbf h\|_{N}\lesssim\epsilon<\infty,
\end{align*}
This contradicts  \eqref{2026-5-9-10}. This completes the proof of Theorem \ref{th:main1}.

\bigskip
\section*{Acknowledgement}
%\noindent{\bf Ethical Approval }\

% We certify that this manuscript is original and has not been published and will not be submitted elsewhere for
%publication while being considered by  Calculus of Variations and Partial Differential Equations.

\vskip .1in
\noindent{\bf Competing interests  }\

On behalf of all authors, the corresponding author states that there is no potential conflicts of interest with respect to the research of this article.
\vskip .1in
\noindent{\bf Authors' contributions   }\

Liening  Qiao, Juntao  Sun,  Jiahong Wu and  Fuyi Xu  contributed equally to this work.
\vskip .1in
\noindent{\bf Funding   }\

J. Sun was supported by National Natural Science Foundation of China 12371174. 
J. Wu was supported by the National Science Foundation of the United
States under DMS 2104682 and DMS 2309748. L. Qiao and  F. Xu were supported by the
National Natural Science Foundation of China 12326430 and the Natural Science Foundation of Shandong Province ZR2026MS0022.
\vskip .1in
\noindent{\bf Availability of data and materials  }\

Data and materials  sharing
not applicable to this article as no data and  materials
 were generated or analyzed during the current study.


\vskip .3in
\begin{thebibliography}{10}

%\bibitem{Chen-Tan-Wu-Zou}Chen, Q.; Tan, Z.; Wu, G; Zou, W.: {\it The initial value problem for the compressible Navier-Stokes equations without heat conductivity,} J. Differential Equations {\bf 268} (2020), no.~9, 5469--5490.
%\bibitem{Chen-Zhang-Zhou} Chen, W.; Zhang, Z.; Zhou, J.: {\it Global well-posedness for the 3-D MHD equations with partial diffusion in periodic domain,} Sci China Math. \textbf{65} (2022), 309--318.
%\bibitem{Danchin-Mucha} Danchin, R.; Mucha, P.: {\it The Incompressible Navier-Stokes Equations in Vacuum,} Communications on Pure and Applied Mathematics
    %\textbf{72}(2019), 1351--1385.
\bibitem{AMS} A. Alemany, R. Moreau, P. Sulem and U. Frisch,  {\it Influence of
	an external magnetic field on homogeneous MHD turbulence,}  J. M\'{e}c., \textbf{18} (1979),
	277--313.
	
   \bibitem{Alex} A. Alexakis,  {\it Two-dimensional behavior of three-dimensional
   magnetohydrodynamic flow with a strong guiding field,}   Phys. Rev. E,
  \textbf{84} (2011), Paper No, 056330.

   \bibitem{Alf} H.  Alfv\'{e}n,  {\it Existence of electromagnetic-hydrodynamic
   waves,}  Nature,  \textbf{150} (1942), 405--406.

   \bibitem{Davi0} P.A. Davidson, {\it Magnetic damping of jets and vortices,}
   J. Fluid Mech.,  \textbf{299} (1995), 153--186.

   \bibitem{Davi1} P.A. Davidson,  {\it The role of angular momentum in the magnetic
   damping of turbulence,}
  J. Fluid Mech.,  \textbf{336} (1997), 123--150.

   \bibitem{DA} P.A. Davidson, {\it An Introduction to Magnetohydrodynamics},
   Cambridge University Press, Cambridge, England, 2001.




%\bibitem{Dong-Wu-Zhai} B. Dong, J. Wu and X. Zhai,  {\it Global small solutions to a special 2$\frac12$-D compressible viscous non-resistive MHD system,} J. Nonlinear Sci., \textbf{33} (2023), Paper No. 21, 37 pp.
%\bibitem{Duan-Ma} Duan, R.; Ma, H.: {\it Global existence and convergence rates for the 3-D compressible Navier-Stokes equations without heat conductivity,} Indiana Univ. Math. J. {\bf 57} (2008), no.~5, 2299--2319.
%\bibitem{Fan-Wang-Zhou} Fan, J.; Wang, L.;  Zhou, Y.: {\it Uniform regularity for the full compressible MHD system with zero heat conductivity,} Math. Methods Appl. Sci. {\bf 47} (2024), no.~13, 10694--10701.
\bibitem{Gao-Wu-Xu} N. Gao, J. Wu and F. Xu, {\it Global large solutions to a multi-dimensional compressible magnetohydrodynamic flows with a nonlinear initial
    constraint,} Calc. Var. Partial Differential Equations, {\bf 64} (2025), Paper No. 73, 38 pp.
%\bibitem{Gao-Tan-Wu} Gao, Z.; Tan, Z.; Wu, G.: {\it Global existence and convergence rates of smooth solutions for the 3-D compressible magnetohydrodynamic equations without heat conductivity,} Acta Math. Sci. Ser. B (Engl. Ed.) {\bf 34} (2014), no.~1, 93--106.
\bibitem{HuX}X. Hu,  {\it Global existence for two dimensional compressible magnetohydrodynamic flows with zero magnetic diffusivity,}	arXiv:1405.0274.


\bibitem{Hu-Wang1} X. Hu and D. Wang, {\it Global solutions to the
three-dimensional full compressible magnetohydrodynamic flows,} Comm. Math.
Phys., \textbf{283} (2008), 255--284.

\bibitem{Hu-Wang2}  X. Hu and D. Wang, {\it Global existence and large-time
behavior of solutions to the three-dimensional equations of compressible
magnetohydrodynamic flows}, Arch. Ration. Mech. Anal., \textbf{197} (2010),
203--238.


%\bibitem{Hu-Wang} Hu, X.; Wang, D.: {\it Global solutions to the three-dimensional full compressible magnetohydrodynamic flows,} Comm. Math. Phys. {\bf 283}
    %(2008), no.~1, 255--284.
%\bibitem{Huang-Wang} Huang, X.; Wang, Y.: {\it $L^\infty$ continuation principle to the non-baratropic non-resistive magnetohydrodynamic equations without heat conductivity,} Math. Methods Appl. Sci. {\bf 39} (2016), no.~14, 4234--4245.
%\bibitem{Jiang-Jiang-1} F. Jiang and S. Jiang,  {\it Nonlinear stability and instability in the Rayleigh-Taylor problem of stratified compressible MHD fluids,}
   % Calc. Var. Partial Differential Equations,  {\bf 58} (2019), Paper No. 29, 61 pp.
\bibitem{Jiang-Jang2} F. Jiang and S. Jiang,  {\it On magnetic inhibition theory in 3D non-resistive magnetohydrodynamic fluids: global existence of large solutions,} Arch. Ration. Mech. Anal., {\bf 247} (2023),  Paper No. 96, 35 pp.
\bibitem{Jiu-Liu-Xie}Q. Jiu, J. Liu and Y. Xie, {\it Asymptotic stability for $n$-dimensional isentropic compressible MHD equations without magnetic diffusion,} arXiv:2402.09661.
\bibitem{Kato-Ponce}T.  Kato and G. Ponce, {\it Commutator estimates and the Euler and Navier-Stokes equations,} Commun. Pure Appl. Math., {\bf 41} (1988),  891--907.

\bibitem{KawashimaS}S. Kawashima,  {\it System of a Hyperbolic-Parabolic Composite Type, with Applications to the Equations of Magnetohydrodynamics,} Ph.D. thesis, Kyoto University (1984).
\bibitem{Kenig-Ponce-Vega} C. Kenig, G. Ponce and  L. Vega, {\it Well-posedness of the initial value problem for the Korteweg-de Vries equation,} J. Am. Math. Soc.,
    \textbf4  (1991), 323--347.
\bibitem{Li-Qiao} J. Li and L. Qiao,  {\it Global well-posedness of the inviscid resistive isentropic
compressible MHD system,} Calc. Var. Partial Differential Equations, {\bf 65} (2026),  Paper No. 247.
\bibitem{Li-Li-Wu-Xu} L. Li, Q. Li, J. Wu and X. Xu, {\it Global stability of 3D compressible non-resistive MHD: hidden damping and rational background fields,} arXiv:2609.14244v1


\bibitem{LiY} Y. Li,  {\it Global well-posedness for the three-dimensional full compressible viscous non-resistive MHD system,} J. Math. Fluid Mech., {\bf 24}  (2022), Paper No. 28, 24 pp.


%\bibitem{Li-Sun} Y. Li and Y. Sun,  {\it On global-in-time weak solutions to a two-dimensional full compressible nonresistive MHD system,} SIAM J. Math. Anal., {\bf 53} (2021), 4142--4177.
\bibitem{Li-Xu-Zhai} Y. Li, H. Xu and X. Zhai, {\it Global smooth solutions to the 3D compressible viscous non-isentropic magnetohydrodynamic flows without   magnetic diffusion,} J. Geom. Anal., {\bf 33} (2023),  Paper No. 246, 32 pp.
%\bibitem{LiangJ} Liang, J.: {\it Global well-posedness of strong solutions to the 3D full compressible magnetohydrodynamics equations with zero heat-conduction,} Appl. Anal. {\bf 104} (2025), no.~16, 3117--3133.
%\bibitem{Peng-Zhai} Peng, H.; Zhai, X.: {\it The Cauchy problem for the $n$-dimensional compressible Navier-Stokes equations without heat conductivity,} SIAM J. Math. Anal. {\bf 55} (2023), no.~2, 1439--1463.


%\bibitem{LXZ}F. Lin, L. Xu and  P. Zhang, {\it Global small solutions to 2-D incompressible MHD system,} J. Differential Equations, \textbf{259} (2015), 5440--5485.


%\bibitem{RWXZ}X. Ren, J. Wu, Z. Xiang and Z. Zhang, {\it Global existence and decay of smooth solution for the 2-D MHD equations without magnetic diffusion,} J. Funct. Anal., \textbf{267} (2014), 503--541.



\bibitem{Pu-Guo} X. Pu and B. Guo,  {\it Global existence and convergence rates of smooth solutions for the full compressible MHD equations,} Z. Angew. Math. Phys., \textbf{64} (2013), 519--538.

\bibitem{Qiao-Zhai}L. Qiao,  J. Wu, F. Xu and X. Zhai, {\it  Stabilization by a background magnetic field: global well-posedness of the compressible isentropic ideal MHD equations with velocity damping,} arXiv:2605.04462.

\bibitem{Qiao-Xu}L. Qiao,  J. Sun, J. Wu and F. Xu, {\it Stabilization by a background magnetic field: global well-posedness of the full compressible viscous non-resistive MHD system without heat-conductivity,} arXiv:2609.01488.

\bibitem{A.S} A. Suen and D. Hoff, {\it Global low-energy weak solutions of the equations of three-dimensional compressible magnetohydrodynamics,} Arch. Ration. Mech. Anal.,  {\bf 205} (2012), 27--58.

\bibitem{SuenA} A. Suen,  {\it Existence and uniqueness of low-energy weak solutions to the compressible 3D magnetohydrodynamics equations,}  J. Differ. Equ., {\bf 268} (2020), 2622--2671.
\bibitem{Tan-Wang} Z. Tan and Y. Wang,  {\it Global well-posedness of an initial-boundary value problem for viscous non-resistive MHD systems,} SIAM J. Math. Anal., {\bf 50} (2018), 1432--1470.
%\bibitem{Tan-Xu-Wang} Tan, Z.; Xu. Q.; Wang, H.:{\it Global existence and convergence rates for the compressible magnetohydrodynamic equations without heat conductivity,} Discrete Contin. Dyn. Syst. {\bf 35} (2015), no.~10, 5083--5105.
\bibitem{Wang-Xin}Y.  Wang and Z.  Xin,  {\it Global well-posedness of the inviscid heat-conductive resistive compressible MHD in a strip domain,} Commun. Math. Res., {\bf 38} (2022),  1--27.
%\bibitem{WangY} Wang, Y.: {\it A Beale-Kato-Majda criterion for three dimensional compressible viscous non-isentropic magnetohydrodynamic flows without  heat-conductivity,} J. Differential Equations {\bf 280} (2021), 66--98.
%\bibitem{Wei-Li-Guo} Wei, R.; Li, Y.; Guo, B.: {\it Global existence and convergence rates of solutions for the 3D compressible magnetohydrodynamic equations without heat conductivity,} Appl. Anal. {\bf 99} (2020), no.~10, 1661--1684.
%\bibitem{WuGC} Wu, G.: {\it Global existence and asymptotic behavior for the 3D compressible Navier-Stokes equations without heat conductivity in a bounded  domain,} J. Differential Equations {\bf 262} (2017), no.~2, 844--861.

%\bibitem{Wu-Tan-Zou} Wu, G.; Tan, Z.; Zou, W.: {\it The initial value problem for the compressible magnetohydrodynamic equations without heat conductivity,} J. Math. Anal. Appl. {\bf 484} (2020), no.~2, 123708, 23 pp.
\bibitem{Wu-Wu}J. Wu and Y. Wu,  {\it Global small solutions to the compressible 2D magnetohydrodynamic system without magnetic diffusion,} Adv. Math., \textbf{310} (2017), 759--888.
\bibitem{Wu-Xu-Zhai} J. Wu, F. Xu and X. Zhai,  {\it Magnetic stabilization of compressible flows: global existence in 3D inviscid non-isentropic MHD equations,} arXiv:2507.00888.
\bibitem{Wu-Zhai} J. Wu and X. Zhai, {\it Global small solutions to the 3D compressible viscous non-resistive MHD system,} Math. Models Methods Appl. Sci., {\bf 33} (2023),  2629--2656.
\bibitem{Wu-Zhu}J. Wu and Y. Zhu,  {\it Global well-posedness for 2D non-resistive compressible MHD system in periodic domain,} J. Funct. Anal., \textbf{283} (2022), Paper No. 109602.
%\bibitem{Xi-Guo-Xie-Fang}Xi, X.; Guo, B.; Xie, B.; Fang, S.: {\it Nonlinear thermal instability in the magnetohydrodynamics problem without heat conductivity,} J. Differential Equations, {\bf 263} (2017), no.~10, 6635--6683.
%\bibitem{Xie-Jiu-Liu} Xie, Y.; Jiu, Q.; Liu, J.: {\it Sharp decay estimates and asymptotic stability for incompressible MHD equations without viscosity or magnetic diffusion,} Calc. Var. Partial Differential Equations {\bf 63} (2024),  Paper No. 191, 24pp.

\bibitem{Xu-Zhong} L. Xu and X. Zhong,  {\it Global well-posedness of the three-dimensional non-isentropic compressible magnetohydrodynamic equations under a scaling-invariant smallness condition,} Calc. Var. Partial Differential Equations,  {\bf 65} (2026),  Paper No. 135, 22 pp.
%\bibitem{Yang-Zhao} W. Yang and X.  Zhao, {\it Global well-posedness and asymptotics of full compressible non-resistive magnetohydrodynamics system with large  external potential forces,} Math. Methods Appl. Sci., {\bf 45} (2022), 206--237.
%\bibitem{YangX} Yang, X.: {\it A regularity criterion for the 2D full compressible MHD equations with zero heat conductivity,} Acta Appl. Math. {\bf 169} (2020), 523--531.
%\bibitem{YuHB2} Yu, H.:{\it Global strong solutions with large oscillations to the 3D full compressible Navier-Stokes equations without heat conductivity,} J. Evol. Equ. {\bf 24} (2024), no.~3, Paper No. 71, 19 pp.
%\bibitem{YuHB1} Yu, H.:{\it Global strong solutions and decay to the 3D full compressible Navier-Stokes equations without heat conductivity,} Discrete Contin. Dyn. Syst. {\bf 45} (2025), no.~6, 1911--1927.
%\bibitem{Zhai-Li-Zhao} X. Zhai, Y. Li and Y. Zhao,  {\it Stability for the $2\frac{1}{2}$-D compressible viscous non-resistive and heat-conducting
    %magnetohydrodynamic flow,} NoDEA Nonlinear Differential Equations Appl. {\bf 32} (2025), Paper No. 101, 30 pp.

\bibitem{Zhai-Wu-Xu} X. Zhai, J. Wu and F.  Xu,  {\it Stability for the 3D magneto-micropolar fluids with only velocity dissipation near a background magnetic field,}  J. Differ. Equ.,  {\bf 425} (2025), 596--626.
%\bibitem{Zhai-Zhang}X. Zhai and S. Zhang,  {\it Global well-posedness and large-time behavior for a special $2\frac{1}{2}\rm D$ full compressible viscous non-resistive MHD system,} Anal. Appl., (Singap.) {\bf 23} (2025), 607--633.
%\bibitem{ZhaoY3} Y. Zhao, Global well-posedness and time-decay for the full compressible non-resistive MHD equations in a 3D infinite slab, Nonlinear Anal. RWA., {\bf 68} (2022), Paper No. 103679, 23 pp.
%\bibitem{ZhaoY1} Y. Zhao, {\it Global well-posedness for the compressible non-resistive MHD equations in a 3D infinite slab,} Nonlinear Anal., {\bf 227} (2023),
    %Paper No. 113162, 21 pp.
\bibitem{ZhaoY2}Y. Zhao,  {\it Global solutions and asymptotic behavior for the three-dimensional viscous non-resistive MHD system with some large
    perturbations,} J. Math. Fluid Mech., {\bf 27} (2025),  Paper No. 44, 25 pp.
%\bibitem{ZhongX} Zhong, X.: {\it Singularity formation to the two-dimensional non-barotropic non-resistive magnetohydrodynamic equations with zero heat conduction in a bounded domain,} Discrete Contin. Dyn. Syst. Ser. B {\bf 25} (2020), no.~3, 1083--1096.

\end{thebibliography}
\end{document}